\documentclass[11pt]{article}
\usepackage{upgreek}
\usepackage{graphicx}
\usepackage{amssymb}
\usepackage{epstopdf}
\usepackage{textcomp}
\usepackage{amsmath}
\usepackage{amsthm}
\usepackage{mathdots}
\usepackage[all,cmtip]{xy}
\usepackage[usenames,dvipsnames]{color}
\input diagxy
\usepackage{cancel}
\usepackage{url}
\usepackage{enumerate}
\newcommand{\A}{\mathcal{A}}

\newcommand{\G}{\mathcal{G}}
\newcommand{\Ca}{\mathcal{C}}

\newcommand{\Gl}[1][{}]{\mathcal{G}_{(#1)}}
\newcommand{\Cl}[1][{}]{\mathcal{C}_{(#1)}}
\newcommand{\N}{\mathbb{N}}

\newcommand{\R}{\mathbb{R}}
\newcommand{\C}{\mathbb{C}}
\newcommand{\id}{\mathrm{id}}
\newcommand{\E}{\mathcal{E}}
\newcommand{\End}{\mathrm{End}\,}
\newcommand{\tr}{\mathrm{tr}}
\newcommand{\Hom}{\mathrm{Hom}\,}

\newcommand{\pr}{\mathrm{pr}}

\newcommand{\set}[2]{ \left\{ #1 \,  : #2 \right\} }

\theoremstyle{plain}
\newtheorem{thm}{Theorem}[section]

\newtheorem{prop}[thm]{Proposition}
\newtheorem{lem}[thm]{Lemma}

\theoremstyle{definition}
\newtheorem{defn}[thm]{Definition}
\newtheorem{ex}[thm]{Example}

\theoremstyle{remark}
\newtheorem{rmk}[thm]{Remark}

\title{Cyclic cocycles on twisted convolution algebras}
\author{Eitan Angel\thanks{Email: eitan.angel@colorado.edu} \\ \\ Department of Mathematics, University of Colorado, UCB 395, \\
Boulder, Colorado 80309-0395, USA}
\begin{document}
\maketitle

\abstract{We give a construction of cyclic cocycles on convolution algebras twisted by gerbes over discrete translation groupoids. For proper \'etale groupoids, Tu and Xu in \cite{tu_xu06} provide a map between the periodic cyclic cohomology of a gerbe-twisted convolution algebra and twisted cohomology groups which is similar to the construction of Mathai and Stevenson in \cite{mathai_stevenson06}.  When the groupoid is not proper, we cannot construct an invariant connection on the gerbe; therefore to study this algebra, we instead develop simplicial techniques to construct a simplicial curvature 3-form representing the class of the gerbe.  Then by using a JLO formula we define a morphism from a simplicial complex twisted by this simplicial curvature 3-form to the mixed bicomplex computing the periodic cyclic cohomology of the twisted convolution algebras. }

\tableofcontents

\section{Introduction}

When a manifold $M$ carries the action of a discrete group $\Gamma$ a natural question to ask is what the orbit space $M / \Gamma$ looks like.  In general the orbit space is not a manifold or even a Hausdorff space. The prescription of Connes' noncommutative geometry to deal with such ``bad'' quotients is to instead study the  algebra $C^\infty_c (M \rtimes \Gamma)$ of smooth functions with compact support on $M \times \Gamma$ with the discrete convolution product.  In \cite{connes94}, Chapter III.2, Connes describes the cyclic cohomology of the convolution algebra $C^\infty_c (M \rtimes \Gamma)$ via a morphism $\Phi$ from the Bott complex of \cite{bott76} computing the Borel model of equivariant cohomology $H^*_\Gamma (M;\C) = H^* (M_\Gamma;\C)$, where $M_\Gamma = M \times_\Gamma E \Gamma$, to the $(b,B)$-bicomplex computing the periodic cyclic cohomology of $C^\infty_c (M \rtimes \Gamma)$, $HP^*(C^\infty_c (M \rtimes \Gamma))$.

The main result of this paper is the construction of a morphism analogous to Connes' $\Phi$ map into the $(b,B)$-bicomplex computing the periodic cyclic cohomology of the convolution algebra twisted by a gerbe over a discrete translation groupoid.  Such a morphism, given by a JLO-type formula, was constructed by Tu and Xu in \cite{tu_xu06} for the case of a gerbe over any proper \'etale groupoid.  Our construction removes the properness requirement in the case of a discrete translation groupoid.  A related construction was also given by Mathai and Stevenson in \cite{mathai_stevenson06} for a canonical smooth subalgebra of stable continuous trace $C^*$-algebras having smooth manifolds as their spectrum (see also \cite{mathai_stevenson07}).  In both \cite{tu_xu06} and \cite{mathai_stevenson06} the morphisms constructed are in fact isomorphisms whereas the morphism constructed in this paper is in general not an isomorphism.


A presentation of an $S^1$-gerbe over the discrete translation groupoid $M \rtimes \Gamma$ is given by a line bundle $L \rightarrow M \rtimes \Gamma$ as well as a collection of line bundle isomorphisms $\mu_{g_1 , g_2} : L_{g_1} \otimes (L_{g_2})^{g_1} \xrightarrow{\sim} L_{g_1 g_2}$ for every $g_1 , g_2 \in \Gamma$, where $L_g$ denotes the restriction of $L$ to $M_g = M \times \{g\}$. Gerbes were introduced by J. Giraud \cite{giraud71} and largely developed in a differential geometric framework by J.-L. Brylinski \cite{brylinski93} to study central extensions of loop groups, line bundles on loop spaces, and the Dirac monopole among other applications.
The appearance of gerbes in string theory and quantum field theory has been studied by numerous mathematicians and physicists, e.g., in \cite{bcmms02}, \cite{gawedzki_reis02} and \cite{cmm00}. Sections of $L$ with the convolution product 
\begin{equation*}
(f_1 * f_2) (x,g) = \sum_{g_1 g_2 = g}  \mu_{g_1 , g_2} (f_1|_{M_{g_1}} \otimes (f_2|_{M_{g_2}})^{g_1} ) (x , g).
\end{equation*}
define the twisted convolution algebra $C_c^\infty (M \rtimes \Gamma, L)$.

In the improper case, one cannot construct a $\Gamma$-invariant connection on $L \rightarrow M \rtimes \Gamma$ as is done in \cite{tu_xu06}.  In particular, this means that we cannot in general choose a connection $\nabla$ on $L \rightarrow M \rtimes \Gamma$ that satisfies the cocycle condition
\begin{equation*}
\mu_{g,h}^* (\nabla_{gh}) = \nabla_g \otimes 1 + 1 \otimes (\nabla_h)^g
\end{equation*}
where $\nabla_g$ denotes $\nabla |_{M_g}$.  Instead, there is some discrepancy $\alpha \in \Omega^1 ( (M \rtimes \Gamma)_{(2)} )$
\begin{equation*}
 \mu_{g,h}^* (\nabla_{g.h})  - \nabla_g \otimes 1 + 1 \otimes (\nabla_h)^g  = \alpha (g,h).
\end{equation*}
The curvature forms of $\nabla_g$, $\theta_g \in \Omega^2 (M)$, define a form $\theta \in \Omega^2 ( (M \rtimes \Gamma)_{(1)} )$.  The $\theta_g$ satisfy
\begin{equation*}
\theta_g -\theta_{gh}+ \theta_h^g = - d\alpha(g,h)
\end{equation*}
so that $(\alpha, \theta)$ is a cocycle that represents the Dixmier-Douady class in the cohomology of the Bott complex of $M \rtimes \Gamma$.

To overcome this difficulty, the geometric data which forms the source of our morphism is constructed using simplicial differential forms on the simplicial manifold $(M \rtimes \Gamma)_\bullet$, the nerve of $M \rtimes \Gamma$.
The geometric data we use is a complex termed the \textit{twisted simplicial complex} 
\begin{equation*}
(\underline{\Omega} ( (M \rtimes \Gamma)_\bullet )[u]^* , d_{\Theta_u})
\end{equation*}
where $\underline{\Omega} (M_\bullet)[u]^*$ is a rescaled version of the complex of compatible forms $\Omega^* (M_\bullet)$ described by Dupont in \cite{dupont}, $u$ is a formal variable, and $d_{\Theta_u}$ is the de Rham differential twisted by a rescaled version of the simplicial curvature 3-form $\Theta$ in $\Omega^3 ((M \rtimes \Gamma)_\bullet)$ which is a representative of the Dixmier-Douady class.  The 3-form $\Theta$ is constructed from the gerbe datum $(L, \mu, \nabla)$, .  From $L \rightarrow M \rtimes \Gamma$ we define an infinite dimensional vector bundle $\E \rightarrow M$ as the direct sum bundle $\E = \bigoplus_{g \in \Gamma} L_g$ with the direct sum connection $\nabla^\E$.  Then $\End \E \rightarrow M$ is a vector bundle of finite rank endomorphisms over $M$.


From the gerbe datum $(L, \mu, \nabla)$ we define connections $\nabla^k$ on the vector bundles $\E_k \times \Delta^k \rightarrow (M \rtimes \Gamma)_{(k)} \times \Delta^k$, where $\E_k = p_k^* \E$ and $p_k : M \times \Gamma^k \rightarrow M$ is the map $p_k : (m, g_1 , \ldots , g_k) \mapsto m \in M$.  While $\End \E_k \rightarrow (M \rtimes \Gamma)_{(k)}$ forms a simplicial vector bundle, the sequence of $\End \E_k$-valued 2-forms defined by $(\nabla^k)^2$ do not define a simplicial form.  Instead we are able to define a simplicial 2-form by the formula
\begin{equation*}
\vartheta_{(k)} (g_1 , \ldots , g_k) = (\nabla^k)^2 - \sum_{i=1}^k t_i \theta_{g_1 \cdots g_i}  + \sum_{1 \leq i < j \leq k} \alpha(g_1 \cdots g_i , g_{i+1} \cdots g_j) (t_i dt_j - t_j dt_i).
\end{equation*}
Here the difference between $\vartheta_{(k)}$ and $(\nabla^k)^2$ is a scalar valued form.  The simplicial curvature 3-form $\Theta_u$ is defined by $(\Theta_u)_{(k)} = \nabla^k_u (\vartheta_u)_{(k)}$, where $\nabla^k_u$ and $\vartheta_u$ are rescaled versions of $\nabla^k$ and $\vartheta$.

A simplicial version of the character formula of A. Jaffe, A. Lesniewski, and K. Osterwalder \cite{jlo88}, similar to \cite{mathai_stevenson06} and \cite{gorokhovsky99} 
defines a morphism into group cochains valued in the $(b,B)$ bicomplex of sections of $\End \E$,
\begin{equation*}
\tau_\nabla : (\underline{\Omega} (M_\bullet) [u]^\bullet , d_{\Theta_u})  \to (C^\bullet (\Gamma , CC^\bullet ( \Gamma^\infty_c ( \End \E) )  [u^{-1} , u] ) , b+uB + {\delta_\Gamma}').
\end{equation*}
Following this we define algebraic morphisms
\begin{equation*}
(C^\bullet (\Gamma , CC^\bullet ( \Gamma^\infty_c ( \End \E) )  [u^{-1} , u] ) , b+uB + {\delta_\Gamma}') \to CC^\bullet (C^\infty_c (M \rtimes \Gamma , L)  [u^{-1} , u], b+uB)
\end{equation*}
into the periodic cyclic complex of the twisted convolution algebra $C^\infty_c (M \rtimes \Gamma , L)$.

The rest of this paper is organized as follows.  The content of Section 2 is a review of simplicial notions and the definition of the twisted simplicial complex $(\underline{\Omega} ( (M \rtimes \Gamma)_\bullet )[u]^* , d_{\Theta_u})$. In Section 3 we discuss gerbes on translation groupoids.  Section 4 describes the construction of the simplicial curvature 3-form $\Theta_u$.  Finally, Section 5 contains the main result, a morphism from the twisted simplicial complex of a gerbe datum $(L, \mu, \nabla)$ on $M \rtimes \Gamma$ to the periodic cyclic complex of the twisted convolution algebra $C^\infty_c (M \rtimes \Gamma, L)$.

\textbf{Acknowledgements.} I would like to thank my advisor, Alexander Gorokhovsky, for his advice, support, and patience.  I would also like to thank Arlan Ramsay for his useful suggestions and helpful discussions.

\section{Preliminaries}


\subsection{Categorical notions}
\label{sec:cat}



\begin{defn}
\label{def:simp_cat}
The \textit{simplicial category} $\Delta$ is the small category of objects $[n]$, where $[n]$ is the ordered set of $n+1$ points, $[n] = \{0 < 1 < \cdots < n \}$, $n \in \N$, and arrows the nondecreasing maps $f: [n] \rightarrow [m]$.  For $n,m \in \N$, a map $f : [n] \rightarrow [m]$ is called \textit{nondecreasing} if $f(i) \geq f(j)$ whenever $i > j$.
\end{defn}


\begin{defn}
\label{def:face_deg}
In $\Delta$ the \textit{face maps} are the injections $\updelta^n_i : [n-1] \rightarrow [n]$, $0 \leq i \leq n$, which skip over the $i^\mathrm{th}$ point, i.e.\ such that $\updelta^n_i (i-1) = i-1$, $\updelta^n_i (i) = i+1$.  The \textit{degeneracy maps} are the surjections $\upsigma^n_j : [n+1] \rightarrow [n]$, $0\leq j \leq n$, that sends both $j$ and $j+1$ to $j$, i.e.\ such that $\upsigma^n_j (j) = \upsigma^n_j (j+1) = j$. When the context is clear, we will omit the superscript $n$.
\end{defn}


\begin{defn}
A \textit{simplicial object} (respectively \textit{cosimplicial object}) in a category $\Ca$ is a contravariant functor $X_\bullet : \Delta \rightarrow \Ca$ (respectively a covariant functor $X^\bullet : \Delta \rightarrow \Ca$).  We will often denote the objects $X_n=X([n])$, $n \in \N$, and the morphisms $\delta_i^n= X_\bullet (\updelta^n_i)$ and $\sigma^n_j=X_\bullet (\upsigma^n_j)$, $0 \leq i , j \leq n$. Such a functor is determined entirely by the objects $X_n$ and the morphisms $d^n_i$ and $s^n_j$ (see e.g. \cite{loday98}). When the morphisms $s^n_j$ are not specified, the functor is instead called a \textit{pre-simplicial object}.

A contravariant functor $X_\bullet :\Delta \rightarrow \mathbf{Sets}$ to the category of sets is called a \emph{simplicial set}. A contravariant functor $X_\bullet : \Delta \rightarrow \mathbf{Top}$ to the category of topological spaces is called a \emph{simplicial space}. A contravariant functor $X_\bullet : \Delta \rightarrow \mathbf{Man}$ to the category of smooth manifolds is called a \emph{simplicial manifold} (see Section \ref{sec:Simplicial_manifolds}).

\end{defn}

\begin{defn}
The \textit{geometric $n$-simplex} is the cosimplicial space $\mathbf{\Delta}^\bullet : \Delta \rightarrow \mathbf{Top}$ defined by $\mathbf{\Delta}^\bullet ([n]) = \Delta^n$, where $\Delta^n$ is the standard $n$-simplex
\begin{align}
\Delta^n &= \set{ (t_0 , \ldots , t_n ) \in \R^{n+1} }{ 0 \leq t_i \leq 1, \ 0 \leq i \leq n, \sum_{i =0}^n t_i = 1} & \textrm{(barycentric)} \\
&= \set{(t_1 , \ldots , t_n) \in \R^{n} }{ t_1 , \ldots , t_n \geq 0, \ \sum_{i=1}^n t_i \leq 1 } & \textrm{(Cartesian).}
\end{align}
The Cartesian coordinates may be obtained from the barycentric coordinates eliminating $t_0 = 1 - t_1 - \cdots - t_n$. For barycentric coordinates
\begin{align}
\mathbf{\Delta}^\bullet (\updelta^n_i) (t_0 , \ldots , t_{n-1}) &= (t_0 , \ldots , t_{i-1} , 0 , t_i , \ldots , t_{n-1}) \\
\mathbf{\Delta}^\bullet (\upsigma^n_j) (t_0 , \ldots ,  t_{n+1}) &=  (t_0 , \ldots , t_{j-1} , t_j + t_{j+1}, \ldots , t_{n+1}),
\end{align}
and for Cartesian coordinates
\begin{align}
\mathbf{\Delta}^\bullet (\updelta^n_i) (t_1 , \ldots , t_{n-1}) &= 
\begin{cases}
(1-t_1 - \cdots - t_{n-1}, t_1 , \ldots , t_{n-1}) & i = 0 \\ 
(t_1 , \ldots , t_{i-1} , 0 , t_i , \ldots , t_{n-1}) & 1 \leq i \leq n
\end{cases} \\
\mathbf{\Delta}^\bullet (\upsigma^n_j) (t_1 , \ldots ,  t_{n+1}) &= 
\begin{cases}
( t_2, \ldots , t_{n+1}) & i=0 \\
(t_1 , \ldots , t_j + t_{j+1} , \ldots , t_{n+1}) & 1 \leq i \leq n
\end{cases}
\end{align}
We will denote $\partial^n_i= \mathbf{\Delta}^\bullet (\updelta^n_i)$ and $\varsigma^n_j=\mathbf{\Delta}^\bullet (\upsigma^n_j)$.
\end{defn}

\begin{defn}
Given a simplicial set (or space or manifold) $X_\bullet$, the \textit{fat realization} of $X_\bullet$, denoted $\| X_\bullet \|$, is the space
\begin{equation}
\| X_\bullet \| = \bigcup_{n \geq 0} X_n \times \Delta^n / \sim
\end{equation}
where $\sim$ is the equivalence relation generated by $(\delta^n_i (x) , t) \sim ( x , \partial^n_i (t))$, for $(x,t) \in X_n \times \Delta^{n-1}$. The \textit{geometric realization} of $X$, denoted by $| X |$ is the space
\begin{equation}
| X_\bullet | = \bigcup_{n \geq 0} X_n \times \Delta^n / \approx
\end{equation}
where $\approx$ is the equivalence relation generated by $(f^* (x) , t) \approx (x , f_* (t))$ for $(x,t) \in X_n \times \Delta^{n-1}$ and any $f \in \hom (\Delta)$.
\end{defn}

\begin{defn}
\label{def:nerve}
Let $\Cl[n]$ consist of $n$-tuples of composable arrows in a (small) category $\Ca$, with $\Cl[0]$ the objects of $\Ca$.  The \textit{nerve} of $\Ca$ is the simplicial set $\Ca_\bullet : \Delta \rightarrow \mathbf{Sets}$ given on objects by $\Ca_\bullet ([n]) = \Cl[n]$ and on morphisms by the faces $\Ca_\bullet (\updelta^n_i) = \delta^n_i$ which compose adjacent morphisms in the $i^\mathrm{th}$ place and by the degeneracies $\Ca_\bullet (\upsigma^n_j) = \sigma^n_j$ which insert the identity morphism in the $j^\mathrm{th}$ place, i.e.\
\begin{align}
\label{eq:nerve}
\delta^n_i (f_1 , \ldots , f_n) &= \begin{cases} (f_2 , \ldots , f_n) & \textrm{if } i=0 \\ (f_1 , \ldots , f_{i} f_{i+1} , \ldots , f_n) & \textrm{for } 1 \leq i \leq n-1 \\ (f_1 , \ldots , f_{n-1}) & \textrm{if } i=n \end{cases} \\
\sigma^n_j (f_1 , \ldots , f_n) &= (f_1 , \ldots , f_j , \id , f_{j+1} , \ldots , f_n)
\end{align}
with $\delta^1_0 (f)$ the terminal object of $f$ and $\delta^1_1 (f)$ the initial object of $f$ and where $f_i f_{i+1} = f_{i+1} \circ f_i$.  The \textit{classifying space} $B \Ca$ of the small category $\Ca$ is the geometric realization of the nerve of $\Ca$.
\end{defn}


\subsection{Simplicial manifolds}

\label{sec:Simplicial_manifolds}

We will now recount the ideas of Dupont in \cite{dupont} (and \cite{felisatti_neumann}) to explicitly construct the cohomology of a simplicial manifold.

\begin{defn}
\label{def:compatible_form}
Given a simplicial manifold $M_\bullet : \Delta \rightarrow \mathbf{Man}$, a \textit{simplicial differential $k$-form} on $M_\bullet$ is a sequence of $k$-forms $\{\omega_{(n)} \}$, where $\omega_{(n)} \in \Omega^k (M_n \times \Delta^n)$, and which satisfy the compatibility condition
\begin{equation}
\label{eq:compatibility_cond}
(\id \times \partial_i)^* \omega_{(n)} = (\delta_i \times \id)^* \omega_{(n-1)}
\end{equation}
on $\Omega^k (M_n \times \Delta^{n-1})$ for all $0 \leq i \leq n$ and $n \geq 1$.  The complex $(\Omega^* (M_\bullet) , d)$ of \textit{compatible forms} consists of

\begin{itemize}
\item the differential graded algebra (dga) of simplicial differential forms on $M_\bullet$ form denoted by $\Omega^* (M_\bullet)$
\item an exterior differential 
\begin{equation*}
d : \Omega^k (M_\bullet) \rightarrow \Omega^{k+1} (M_\bullet)
\end{equation*}
induced by the exterior differentials $d : \Omega^k (M_n \times \Delta^n) \rightarrow \Omega^{k+1} (M_n \times \Delta^n)$ for all $n \in \N$.  Explicitly, $d \{ \omega_{(n)} \} = \{ d \omega_{(n)} \}$, which is well defined as the de Rham differential on $M_n \times \Delta^{n-1}$ respects the compatibility condition \eqref{eq:compatibility_cond} for all $n \in \N$.

\item a wedge product
\begin{equation*}
\wedge : \Omega^k (M_\bullet) \otimes \Omega^l (M_\bullet) \rightarrow \Omega^{k+l} (M_\bullet)
\end{equation*}
is induced by the wedge products $\wedge : \Omega^k (M_n \times \Delta^n) \otimes \Omega^l (M_n \times \Delta^n) \rightarrow \Omega^{k+l} (M_n \times \Delta^n)$ for all $n \in \N$.
\end{itemize}
\end{defn}

Notice that $\{\omega_{(n)} \}$ defines a $k$-form on $ \coprod_{n \geq 0} M_n \times \Delta^n$.
In view of this, the compatibility condition \eqref{eq:compatibility_cond} is precisely the condition required for $\{ \omega_{(n)} \}$ to define a form on $\|M_\bullet\|$.

\begin{defn}
We may think of the complex $(\Omega^* (M_\bullet), d)$ as a bicomplex $(\Omega^{r,s} (M_\bullet) , d_{dR} , d_\Delta')$ called the \textit{bicomplex of compatible forms} such that 
\begin{equation}
(\Omega^* (M_\bullet), d) = \mathrm{Tot} \, (\Omega^{r,s} (M_\bullet) , d_{dR} , d_\Delta').
\end{equation}
The bicomplex of compatible forms consists of
\begin{itemize}
\item simplicial $(r+s)$-forms on $M_\bullet$
\begin{equation}
\label{eq:bicompatible_form}
 \Omega^{r,s} (M_\bullet) = \left( \coprod_{n \geq 0} \Omega^r (M_n) \otimes \Omega^s (\Delta^n)  \right) \Big{/} \sim
 \end{equation}
where $\sim$ is essentially the compatibility condition \eqref{eq:compatibility_cond}.  Explicitly, if $\omega_{(n)}^{r,s} \in \Omega^r (M_n) \otimes \Omega^s (\Delta^n)$ we may write $\omega_{(n)}^{r,s}$ as a linear combination of forms $\alpha_{(n)}^r \otimes \beta_{(n)}^s$ where $\alpha_{(n)}^r \in \Omega^r (M_n)$ and $\beta_{(n)}^s \in \Omega^s (\Delta^n)$.  Then $\sim$ is defined by
\begin{equation}
\label{eq:bicompatibility_cond}
 \alpha_{(n)}^r \otimes \beta_{(n)}^s \sim \alpha_{(n-1)}^r \otimes \beta_{(n-1)}^s \Longleftrightarrow \alpha_{(n)}^r \otimes \partial_i^* (\beta_{(n)}^s) = \delta_i^* (\alpha_{(n-1)}^r ) \otimes \beta_{(n-1)}^s,
\end{equation}

\item a differential
\begin{equation}
d_{dR} = d : \Omega^{r,s} (M_\bullet) \rightarrow \Omega^{r+1,s} (M_\bullet)
\end{equation}
induced by the collection of exterior differentials on each $M_n$,
\begin{equation}
d \otimes \id : \Omega^r (M_n) \otimes \Omega^s (\Delta^n) \rightarrow \Omega^{r+1} (M_n) \otimes \Omega^s (\Delta^n),
\end{equation}
for all $n \in \N$.

\item a differential
\begin{equation}
d_\Delta' = (-1)^r d : \Omega^{r,s} (M_\bullet) \rightarrow \Omega^{r,s+1} (M_\bullet)
\end{equation}
induced by the collection of exterior differentials on each $\Delta^n$,
\begin{equation}
\id \otimes d : \Omega^r (M_n) \otimes \Omega^s (\Delta^n) \rightarrow \Omega^r (M_n) \otimes \Omega^{s+1} (\Delta^n),
\end{equation}
for all $n \in \N$.
\end{itemize}
As before, the induced differentials on $\Omega^{r,s} (M_\bullet)$ are well defined as a result of  the compatibility condition \eqref{eq:bicompatibility_cond}.  With this notation the complex of compatible forms may be written $(\Omega^* (M_\bullet) , d_{dR} + d_\Delta')$.

\end{defn}

As the dga of compatible forms carries a bigrading, 
\begin{equation}
\label{eq:bicompatible_sum}
\Omega^k (M_\bullet) = \bigoplus_{r+s=k} \Omega^{r,s} (M_\bullet),
\end{equation}
given an element of the complex of of compatible forms $\omega \in \Omega^k (M_\bullet)$ we will denote the component of $\omega$ in $\Omega^{r,s} (M_\bullet)$ by 
\begin{equation}
\omega^{r,s} = \omega |_{\Omega^{r,s} (M_\bullet)}.
\end{equation}
The $(r+s)$-forms of $\Omega^{r,s} (M_\bullet)$ have a local description as
\begin{equation}
\omega^{r,s} |_{M_n \times \Delta^n} = \sum \omega_{i_1 \cdots i_r j_1 \cdots j_s} dx_{i_1} \wedge \cdots \wedge dx_{i_r} \wedge dt_{j_1} \wedge \cdots \wedge dt_{j_s}
\end{equation}
where $\{x_i\}$ are local coordinates of $M_n$ and $(t_0 , \ldots , t_n)$ are barycentric coordinates of $\Delta^n$.  

\begin{defn}
The \textit{simplicial de Rham complex} $(\A^* (M_\bullet) , \delta + d')$ is defined by 
\begin{equation}
\A^k (M_\bullet) = \bigoplus_{r+s = k} \A^{r,s} (M_\bullet)
\end{equation}
where 
\begin{itemize}
\item $\A^{r,s} (M_\bullet) = \Omega^s (M_r)$ is the collection of differential $s$-forms on the manifold $M_r$,
\item the differential
\begin{equation}
\delta : \A^{r,s} (M_\bullet) \rightarrow \A^{r+1,s} (M_\bullet)
\end{equation}
is the alternating sum of the pullback of the face maps $\delta^{r+1}_i = M_\bullet (\updelta^n_i)$, $0 \leq i \leq r+1$,
\item the differential
\begin{equation}
d' : \A^{r,s} (M_\bullet) \rightarrow \A^{r,s+1} (M_\bullet)
\end{equation}
is the exterior differential $(-1)^r d : \Omega^s (M_r) \rightarrow \Omega^{s+1} (M_r)$.
\end{itemize}
Then we have
\begin{equation}
(\A^* (M_\bullet) , \delta+ d') = \mathrm{Tot} \, (\A^{r,s} (M_\bullet),  \delta,  d').
\end{equation}

\end{defn}

\begin{defn}
Given a commutative ring $R$, the \textit{simplicial singular cochain complex} associated to a simplicial manifold $M_\bullet$ is denoted $(C^*(M_\bullet ; R) , \delta + \partial')$ and consists of
\begin{itemize}
\item the spaces
\begin{equation}
C^k (M_\bullet ; R) = \bigoplus_{r+s=k} C^{r,s} (M_\bullet ; R)
\end{equation}
where $C^{r,s} (M_\bullet ; R) = C^s (M_r ; R)$ is the space of singular cochains of degree $s$ on $M_r$,
\item the differential
\begin{equation}
\delta : C^{r,s} (M_\bullet; R) \rightarrow C^{r+1,s} (M_\bullet; R)
\end{equation}
is the alternating sum of the pullback of the face maps $\delta^{r+1}_i = M_\bullet (\updelta^n_i)$, for all $0 \leq i \leq r+1$,
\item and
\begin{equation}
\partial' : C^{r,s} (M_\bullet; R) \rightarrow C^{r,s+1} (M_\bullet; R)
\end{equation}
is $(-1)^r$ times the usual coboundary on singular cochains.
\end{itemize}
Then we have 
\begin{equation}
(C^* (M_\bullet ; R) , \delta + \partial') = \mathrm{Tot} \, (C^{r,s} (M_\bullet ; R) , \delta , \partial').
\end{equation}
\end{defn}

From \cite{dupont} Proposition 5.15
\begin{thm}
For a simplicial manifold $M_\bullet$, we have the isomorphism
\begin{equation}
H^* (\| M_\bullet \|; R) \cong H(C^* (M_\bullet ; R) , \delta + \partial')
\end{equation}
\end{thm}
and from \cite{dupont} Proposition 6.1,
\begin{thm}[Simplicial de Rham theorem]
The integration map
\begin{equation}
\mathcal{I} : \A^{r,s} (M_\bullet) \rightarrow C^{r,s} (M_\bullet)
\end{equation}
defined by $\mathcal{I} (\omega^{r,s}) (c_{r,s}) = \int_{c_{r,s}} \omega^{r,s}$ for $\omega^{r,s} \in \A^{r,s} (M_\bullet)$ and $c_{r,s} \in C_s (M_r)$, the collection of singular $s$-chains on $M_r$, gives a morphism of double complexes. Furthermore, this integration map induces an isomorphism
\begin{equation}
H(\A^* (M_\bullet) , \delta + d') \cong H(C^* (M_\bullet) , \delta + \partial')
\end{equation}
on cohomology.
\end{thm}

Furthermore, there is another morphism of complexes given by Stokes' theorem. From \cite{dupont} Theorem 6.4,
\begin{thm}
\label{thm:simplex_integration}
Let 
\begin{equation}
\mathcal{I}_\Delta : (\Omega^{r,s} (M_\bullet) , d_{dR} , d_\Delta) \rightarrow (\A^{r,s} (M_\bullet) , \delta , d')
\end{equation}
be the map defined by integration over the standard simplex, i.e. the map defined on $\Omega^s (M_r \times \Delta^r)$ by
\begin{equation}
\mathcal{I}_\Delta : \omega_{(r)} \mapsto \int_{\Delta^r} \omega_{(r)}.
\end{equation}
Then $\mathcal{I}_\Delta$ is a morphism that induces an isomorphism
\begin{equation}
H(\Omega^* (M_\bullet) , d') \cong H((\A^* (M_\bullet) , \delta+d')
\end{equation}

\end{thm}

\subsection{Connections and Curvature}

\begin{defn}
Let $G$ be a Lie group.  A \textit{simplicial $G$-bundle} $\pi_\bullet : E_\bullet \rightarrow M_\bullet$ over a simplicial manifold $M_\bullet$ is a sequence of principal $G$-bundles $\{ \pi_n : E_n \rightarrow M_n\}$ where $E_\bullet$ is itself a simplicial manifold and the diagrams
\begin{equation}
\label{eq:simp_bund_diag}
\xymatrix{
E_{n+1}  \ar[rr]^{E_{n+1} (\updelta_i)} \ar[d]_{\pi_{n+1}} && E_{n} \ar[d]^{\pi_{n}} && E_{n-1} \ar[rr]^{E_{n-1} (\upsigma_j)} \ar[d]_{\pi_{n-1}} && E_{n} \ar[d]^{\pi_{n}} \\
M_{n+1} \ar[rr]^{M_{n+1} (\updelta_i)} && M_{n} & \textrm{and} & M_{n-1} \ar[rr]^{M_{n-1} (\upsigma_j)} && M_{n}
}
\end{equation}
commute.  Given a simplicial $G$-bundle $\pi_\bullet : E_\bullet \rightarrow M_\bullet$, the geometric realization $|\pi_\bullet| :| E_\bullet | \rightarrow | M_\bullet|$ is a principal $G$-bundle with $G$-action induced by
\begin{equation}
E_n \times \Delta^n \times G \rightarrow E_n \times \Delta^n, \quad (x,t,g) \mapsto (xg, t).
\end{equation}
If we only require that the first diagram of \eqref{eq:simp_bund_diag} commute then we may still consider the fat realization $\|\pi_\bullet\| :\| E_\bullet \| \rightarrow \| M_\bullet \|$ which is a principal $G$-bundle.
\end{defn}

\begin{defn}
A \textit{connection} in a simplicial $G$-bundle $\pi_\bullet : E_\bullet \rightarrow M_\bullet$ is a $1$-form $\omega \in \Omega^1 (E_\bullet ; \mathfrak{g})$ on $E_\bullet$ (in the sense of definition \ref{def:compatible_form}) with coefficients in $\mathfrak{g}$, the Lie algebra of $G$, such that $\omega_{(n)} = \omega|_{E_n \times \Delta^n}$ is a connection in the usual sense on the bundle $\pi_n \times \id : E_n \times \Delta^n \rightarrow M_n \times \Delta^n$. The curvature $\theta$ of a connection $\omega$ is the differential form
\begin{equation}
\theta = d\omega + \frac{1}{2} [\omega , \omega] \in \Omega^2 (M_\bullet ; \mathfrak{g}).
\end{equation}

\end{defn}


\subsection{Twisted simplicial cohomology}

In this section we will construct a twisted version of the complex of compatible forms, although rather than using forms, this construction will be based on densities.  We will first define a version of the de Rham complex of a manifold twisted by a $3$-form in preparation for the definition of the de Rham complex of a simplicial manifold twisted by a simplicial $3$-form.

For $\tau$ the orientation bundle of a manifold $M$, let $\Omega_\tau^* (M)$ be the \emph{densities} of $M$. The cohomology of the \emph{$\tau$-twisted de Rham complex}, $(\Omega_\tau^* (M) , d)$, is the \emph{$\tau$-twisted de Rham cohomology}, $H_\tau^* (M)$. The \emph{$\tau$-twisted de Rham cohomology} with compact support, $H_{\tau,c}^* (M)$ is defined similarly. For more details, see \cite{bott_tu} Chapter I.7.

In the following, let $u$ be a formal variable of degree $+2$. 

\begin{defn}

Given a smooth manifold $M$ and a closed $3$-form $\Theta \in \Omega^3 (M)$, the \textit{$\Theta$-twisted de Rham complex of $M$}, denoted by $(\underline\Omega (M)[u]^\bullet , d_\Theta)$, is defined as follows.
\begin{itemize}
\item Define 
\begin{equation}
\label{eq:twisted_complex_grading}
\Omega_k (M) := \Omega_\tau^{\dim M - k } (M).
\end{equation}
  Let $\underline\Omega (M) [u]^k$ denote polynomials in $u$ over $\Omega_* (M)$ such that the degree in $\Omega_* (M)$ plus the degree of powers of $u$ sum to $k$.  In other words elements $\omega \in \underline\Omega (M)[u]^k$ are linear combinations of $u^j \omega_\ell$ for any $\omega_\ell \in \Omega_\ell (M)$ such that $2j + \ell = k$, for all $0 \leq \ell \leq \dim M$.

\item Note the de Rham differential 
\begin{equation} d_{dR} : \Omega_k (M) \rightarrow \Omega_{k-1} (M) \end{equation}
is of degree $-1$ and exterior multiplication by $\Theta$, 
\begin{equation} \Theta \wedge \cdot : \Omega_k (M) \rightarrow \Omega_{k-3} (M)\end{equation}
is of degree $-3$.  We can rescale the de Rham differential by $u$ to obtain a degree $+1$ differential
\begin{equation}
u d_{dR} : \underline\Omega (M) [u]^k \rightarrow \underline\Omega (M) [u]^{k+1}.
\end{equation}

The differential 
\begin{equation}
d_\Theta := ud_{dR} + u^2 \Theta \wedge \cdot : \underline\Omega (M)[u]^k \rightarrow \underline\Omega (M)[u]^{k+1}
\end{equation}
is of degree $+1$.
\end{itemize}
The cohomology of $(\underline\Omega (M) [u]^\bullet , d_\Theta)$ is called the \emph{$\Theta$-twisted cohomology of $M$}, denoted $H^*_\Theta (M)$.

\end{defn}

If $\Theta ' = \Theta + d \eta$ is cohomologous to $\Theta$ then the complexes $(\underline\Omega (M)[u]^\bullet , d_\Theta)$ and $(\underline\Omega (M)[u]^\bullet, d_{\Theta'})$ are isomorphic via the isomorphism
\begin{equation}
I_\eta : \xi \mapsto e^{-u \eta} \wedge \xi
\end{equation}

Now we will similarly describe the de Rham complex of a simplicial manifold $M_\bullet$ in which the twisting of the complex of compatible forms is given by a closed compatible $3$-form $\Theta \in \Omega^3 (M_\bullet)$. This definition involves a modification of the bicomplex of compatible forms \eqref{eq:bicompatible_form}.

\begin{defn}

Given a simplicial manifold $M$ and a closed compatible $3$-form $\Theta \in \Omega^3 (M_\bullet)$, the \textit{$\Theta$-twisted complex of compatible forms of $M_\bullet$}, denoted by $(\underline\Omega (M_\bullet)[u]^\bullet , d_{\Theta_u})$, is defined as follows.


\begin{itemize}

\item First we adjust the grading of $\Omega^{r,s} (M_\bullet)$ to define a bicomplex $\underline\Omega (M_\bullet)[u]^{r , s}$.  Let
\begin{equation}
\underline\Omega (M_\bullet) [u]^{r,s} = \left( \coprod_{n \geq 0} \underline{\Omega} (M_n) [u]^r \otimes \Omega^s (\Delta^n) \right) \Big{/} \sim
\end{equation}
where the equivalence relation $\sim$ is defined in precisely the same way as the compatibility condition in \eqref{eq:bicompatibility_cond}.  For $\alpha^r_{(n)} \in \underline{\Omega} (M_n) [u]^r$ and $\beta^s_{(n)} \in \Omega^s (\Delta^n)$, we have $\alpha_{(n)}^r \otimes \beta_{(n)}^s \sim \alpha_{(n-1)}^r \otimes \beta_{(n-1)}^s$ if and only if the corresponding statement to \eqref{eq:bicompatibility_cond} holds.  As in \eqref{eq:bicompatible_sum}, we define
\begin{equation}
\underline\Omega (M_\bullet) [u]^k = \bigoplus_{r+s=k} \underline\Omega (M_\bullet)[u]^{r,s}.
\end{equation}

Note the degrees are such that for a $(\dim M - \ell)$-form $\omega_\ell \in \Omega_\ell (M_n)$ and $\eta^s \in \Omega^s (\Delta^n)$, 
\begin{equation*}
u^j \omega_\ell \otimes \eta^s \in ( \underline{\Omega} (M_\bullet) [u]^{2j+\ell , s}) |_{M_n \times \Delta^n}
\end{equation*}

\item Furthermore, there is a degree $+1$ differential
\begin{equation}
u\tilde{d}_{dR} : \underline\Omega (M_\bullet) [u]^{r,s} \rightarrow \underline\Omega (M_\bullet) [u]^{r+1,s}
\end{equation}
 induced by the exterior derivative on each $M_n$,
\begin{equation}
(-1)^{\dim M - r+s}ud \otimes \id : \underline{\Omega} (M_n) [u]^r \otimes \Omega^s (\Delta^n) \rightarrow \underline{\Omega} (M_n) [u]^{r+1} (M_n)[u] \otimes \Omega^s (\Delta^n)
\end{equation}
and there is a degree $+1$ differential $d_\Delta$ induced by the exterior derivative on each $\Delta^n$, i.e.\
\begin{equation}
d_{\Delta} =  \id \otimes d : \underline{\Omega} (M_n) [u]^r  \otimes \Omega^s (\Delta^n) \rightarrow \underline{\Omega} (M_n) [u]^r \otimes \Omega^{s+1} (\Delta^n).
\end{equation}


\item Given a compatible differential $3$-form $\Theta \in \Omega^3 (M_\bullet)$ as in Definition \ref{def:compatible_form} we may decompose $\Theta$ with respect to the sum \eqref{eq:bicompatible_sum} as $\Theta = \Theta^{3,0} + \Theta^{2,1} + \Theta^{1,2} + \Theta^{0,3}$ where $\Theta^{j,3-j} \in \Omega^{j,3-j} (M_\bullet)$ for $0 \leq j \leq 3$.  If, further, $\Theta^{0,3} = 0$ then define a rescaling of $\Theta$
\begin{equation}
\label{eq:Theta_rescaled}
\Theta_u = u^2 \Theta^{3,0} + u \Theta^{2,1} + \Theta^{1,2} \in \underline\Omega (M_\bullet)[u]^\bullet.
\end{equation}
Hence, exterior multiplication by $\Theta_u$ is an operator of degree $+1$ on
\begin{equation}
\Theta_u \wedge \cdot : \underline\Omega (M_\bullet)[u]^k \rightarrow \underline\Omega (M_\bullet)[u]^{k+1}.
\end{equation}

\item So for a closed, compatible 3-form $\Theta$ we define
\begin{equation}
d_{\Theta_u} := u\tilde{d}_{dR} + d_\Delta - \Theta_u \wedge \cdot : \underline\Omega (M_\bullet)[u]^k \rightarrow \underline\Omega (M_\bullet)[u]^{k+1}.
\end{equation}
\end{itemize}
The cohomology of $(\underline\Omega (M_\bullet) [u]^\bullet , d_{\Theta_u})$ is called the \emph{$\Theta$-twisted cohomology of $M_\bullet$}, denoted $H^*_\Theta (M_\bullet)$.

\end{defn}

\begin{rmk}
As the nerve of any Lie groupoid $\G$ (see Definition \ref{defn:Lie_etale_groupoid}) defines a simplicial manifold $\G_\bullet$, we may speak of the complex of compatible forms and, given a closed, compatible 3-form $\Theta$, the $\Theta$-twisted complex of compatible forms of $\G$.

\end{rmk}

\begin{prop}
Let $\Theta , \Theta ' \in \Omega^3 (M_\bullet)$ be compatible forms such that the components in $\Omega^{0,3} (M_\bullet)$ satisfy $(\Theta)^{0,3} = (\Theta')^{0,3}  = 0$.  If $\Theta - \Theta ' = d \eta$ for some $\eta \in \Omega^2 (M_\bullet)$ that satisfies $\eta^{0,2}  = 0$ then $(\underline\Omega (M_\bullet)[u]^\bullet , d_{\Theta_u})$ and $(\underline\Omega (M_\bullet)[u]^\bullet , d_{\Theta'_u})$ are isomorphic via the isomorphism
\begin{equation}
I_\eta : \xi \mapsto e^{- \eta_u} \wedge \xi
\end{equation}
where $\eta_u = u \eta^{2,0} + \eta^{1,1}$.
\end{prop}

\begin{lem}
\label{lem:wedge_product_twisted_simplicial_complex}
For $\omega_1 \in \underline\Omega (M_\bullet)[u]^k$ and $\omega_2 \in \Omega^\ell (M_\bullet)[u]$ the induced wedge product satisfies
\begin{equation*}
(-1)^{|\omega_2|} (u \tilde{d}_{dR} + d_\Delta) (\omega_1 \wedge \omega_2) = ( (u \tilde{d}_{dR} + d_\Delta) \omega_1) \wedge \omega_2 +  \omega_1 \wedge ( (u d_{dR} + d_\Delta) \omega_2)
\end{equation*}
\begin{proof}
As the variable $u$ is of even degree and will not alter any signs, let $\tilde{d} = \tilde{d}_{dR} + d_\Delta$ and let $d = d_{dR} + d_\Delta$. Suppose $\omega_1 \in \underline\Omega (M_\bullet)[u]^{r_1,s_1}$ and $\omega_2 \in \Omega^{r_2 , s_2} (M_\bullet)[u]$ so that
\begin{align*}
\tilde{d} (\omega_1 \wedge \omega_2) &= (-1)^{\dim M - r_1 + r_2 + s_1 + s_2} d (\omega_1 \wedge \omega_2) \\
&= (-1)^{\dim M - r_1 + r_2 + s_1 + s_2} ( d \omega_1 \wedge \omega_2 + (-1)^{\dim M - r_1 + s_1} \omega_1 \wedge d \omega_2)\\
&= (-1)^{\dim M - r_1 + r_2 + s_1 + s_2} \left( (-1)^{\dim M - r_1 + s_1} \tilde{d} \omega_1 \wedge \omega_2 \right. \\
&\quad \left. + (-1)^{\dim M - r_1 + s_1} \omega_1 \wedge d \omega_2 \right)\\
&= (-1)^{-r_2 + s_2} \tilde{d} \omega_1 \wedge \omega_2 + (-1)^{-r_2 + s_2} \omega_1 \wedge d\omega_2
\end{align*}
Since $|\omega_2| = r_2+s_2$ which is $s_2 - r_2$ modulo 2,
\begin{equation*}
\tilde{d} (\omega_1 \wedge \omega_2) = (-1)^{|\omega_2|} \tilde{d} \omega_1 \wedge \omega_2 + (-1)^{|\omega_2|} \omega_1 \wedge d\omega_2
\end{equation*}
from which the lemma follows.
\end{proof}
\end{lem}


\section{$S^1$-gerbes}

\subsection{\'Etale groupoids}
\label{sec:gpd}


We will recall some standard material about groupoids.  See \cite{moerdijk02}, \cite{crainic99}, or \cite{moerdijk_mrcun} for example.

\begin{defn}
A \textit{groupoid} is a small category $\G$ in which every arrow is invertible. The set of objects is denoted by $\Gl[0]$ and the set of arrows is denoted by $\Gl[1]$.  The set of arrows will often be denoted simply by $\G$. For each arrow in $\Gl[1]$ there is a source object and a range object given by the \textit{range and source maps}, $r,s : \Gl[1] \rightrightarrows \Gl[0]$.  To denote that $g \in \Gl[1]$ is an arrow with source $s(g) = x$ and range $r(g) = y$ we write either $g : x \rightarrow y$ or $x \xrightarrow{g} y$.

As $\G$ is a category, there is a rule of composition.  Given two arrows $y \xrightarrow{g_1} z, x \xrightarrow{g_2} y \in \Gl[1]$ such that $s(g_1) = r(g_2)$, their composition is an arrow $x \xrightarrow{g_1 g_2} z$ called their \textit{product}.  If we define pairs of composable arrows as
\begin{equation}
\Gl[2] = \Gl[1] \times_{\Gl[0]} \Gl[1] = \{ (g_1 , g_2) \in \Gl[1] \times \Gl[1] : s(g_1) = r(g_2) \}
\end{equation}
then the product defines the \textit{multiplication map}
\begin{equation}
m : \Gl[2] \rightarrow \Gl[1], \qquad m(g_1 , g_2) = g_1g_2.
\end{equation}
As the composition in a category is required to be associative, the groupoid product is associative.

For any object $x \in \Gl[0]$ there is an identity morphism $1_x : x \rightarrow x$ that satisfies $1_x g = g 1_y = g$ for any arrow $x \xrightarrow{g} y \in \Gl[1]$. The \textit{unit map} is then $$u : \Gl[0] \rightarrow \Gl[1], \qquad u(x) = 1_x.$$

Every arrow $g : x \rightarrow y$ in $\Gl[1]$ is invertible so we will denote the inverse of $g$ by $g^{-1} : y \rightarrow x$.  With this we can define the \textit{inverse map} $$i : \Gl[1] \rightarrow \Gl[1], \qquad i(g) = g^{-1}.$$ Then $g^{-1} g = 1_x$ and $g g^{-1} = 1_y$.
\end{defn}

\begin{rmk}
As a groupoid $\G$ is a category, Definition \ref{def:nerve} applies and there is a simplicial set $\G_\bullet$ called the nerve of $\G$.
\end{rmk}

\begin{defn}
As in \cite{bgnt07}, Section 2.3, we will define the following projection maps. For $0 \leq i \leq n$ and $\G_\bullet$ the nerve of $\G$ let $\pr^n_i : \Gl[n] \rightarrow \Gl[0]$ be the final object of the $i^\mathrm{th}$ morphism when $i \neq n$ and the initial object of the last morphism when $i=n$. For $0 \leq j \leq n$, $0 \leq i_j \leq m$, the map $\pr^n_{i_0} \times \cdots \times \pr^n_{i_m} : \Gl[n] \rightarrow (\Gl[0])^m$ factors into maps $\Gl[n] \rightarrow \Gl[m]$ and the canonical projection $\Gl[m] \rightarrow (\Gl[0])^m$. The map $\Gl[n] \rightarrow \Gl[m]$ will be denoted $\pr^n_{i_0 \cdots i_m}$. We may write $\pr^n_{i_0 \cdots i_m}$ explicitly as
\begin{equation*}
\pr^n_{i_0 \cdots i_m} : (g_1 , \ldots, g_n) \mapsto ( g_1 ' , \ldots , g_m ')
\end{equation*}
where, for $1 \leq k \leq m$
\begin{equation*}
g_k ' = \begin{cases}
g_{i_{k-1} +1} \cdots g_{i_{k}} & \text{if $i_{k-1} < i_{k}$} \\
\id_{s(g_{i_k})} & \text{if $i_{k-1} = i_{k}$} \\
(g_{i_{k}} \cdots g_{i_{k-1} - 1})^{-1} & \text{if $i_{k-1} > i_{k}$}
\end{cases}
\end{equation*}
\end{defn}

\begin{rmk}
The maps $\delta^n_i$ of Definition \ref{def:nerve} may be written $\pr^n_{0\cdots \hat{i} \cdots n} : \Gl[n] \rightarrow \Gl[n-1]$ where $\hat{i}$ denotes the omission of $i$.
\end{rmk}

\begin{defn}
\label{defn:Lie_etale_groupoid}
A \textit{Lie groupoid} is a groupoid in which the sets $\Gl[0]$ and $\Gl[1]$ are smooth manifolds and the structure maps $r,s,u,i,m$ are smooth.  Furthermore, the range and source maps $r,s : \Gl[1] \rightrightarrows \Gl[0]$ are required to be submersions so that the domain of the multiplication map, $\Gl[2]$, is a manifold.  An \textit{\'etale groupoid} is a Lie groupoid in which the source map is \'etale, i.e.\ a local diffeomorphism.  In this case the other structure maps are \'etale as well.
\end{defn}

\begin{ex} 
\begin{enumerate}
\item Given a manifold $M$ with an open cover $\mathcal{U} = \{U_\alpha\}$, the \textit{\v{C}ech groupoid} associated to $\mathcal{U}$ has objects $ \coprod_\alpha U_\alpha$ and morphisms $\coprod_{\alpha, \beta} U_{\alpha \beta}$ where $U_{\alpha \beta}  = U_\alpha \cap U_\beta$.  The range and source maps are the embeddings $r : U_{\alpha \beta} \rightarrow U_\beta$ and $s : U_{\alpha \beta} \rightarrow U_\alpha$.
\item Given a manifold $M$ carrying a right action of a group $G$, there is a groupoid called the \textit{translation groupoid} or \textit{action groupoid} $M \rtimes G$ with objects $(M \rtimes G)_{(0)} = M$ and arrows $(M \rtimes G)_{(1)} = M \times G$. The source and range maps are given by $s(x,g) = xg$ and $r(x,g) = x$. Given arrows $xg_1 \xleftarrow{g_1} x$ and $xg_1 g_2 \xleftarrow{g_2} xg_1$ their composition is $m((x,g_1) ,(xg_1, g_2)) = (x, g_1 g_2) = xg_1 g_2 \xleftarrow{g_1 g_2} x$.

\end{enumerate}
\end{ex}

\subsection{$S^1$-gerbes over \'etale groupoids}

\begin{defn}
Given an \'etale groupoid $\G$, a \emph{gerbe datum} is a triple $(L, \mu, \nabla)$ with
\begin{itemize}
\item $L$ a line bundle $L \rightarrow \G$ over $\Gl[1]$;
\item $\mu$ an isomorphism
\begin{equation}
\label{eq:gerbe_multiplication}
\mu : (\pr^2_{01})^* L \otimes (\pr^2_{12})^* L \xrightarrow{\sim} (\pr^2_{02})^* L
\end{equation}
of line bundles over $\Gl[2]$ which satisfies the associativity condition that the diagram
\begin{equation}
\bfig
\Square|aarb|[(\pr^3_{01})^* L \otimes (\pr^3_{12})^* L \otimes (\pr^3_{23})^* L ` (\pr^3_{012})^* L \otimes (\pr^3_{23})^* L ` (\pr^3_{01})^* L \otimes (\pr^3_{123})^* L ` (\pr^3_{0123})^* L; (\pr^3_{012})^* (\mu) \otimes \id ` \id \otimes (\pr^3_{123})^* (\mu) ` (\pr^3_{023})^* (\mu) ` (\pr^3_{013})^* (\mu)]
\efig
\end{equation}
commutes. Such an isomorphism may be written as $\mu_{(g_1 , g_2)} : L|_{g_1} \otimes L|_{g_2} \xrightarrow{\sim} L|_{g_1 g_2}$ on the fiber over $(g_1 , g_2) \in \Gl[2]$.  Hence on fibers this condition means that the diagram
\begin{equation}
\bfig
\Square|aarb|[L|_{g_1} \otimes L|_{g_2} \otimes L|_{g_3} ` L|_{g_1 g_2} \otimes L|_{g_3} ` L|_{g_1} \otimes L|_{g_2 g_3} ` L|_{g_1 g_2 g_3} ; \mu_{(g_1 , g_2)} \otimes \id ` \id \otimes \mu_{(g_2, g_3)} ` \mu_{(g_1 g_2 , g_3)} ` \mu_{(g_1 , g_2 g_3)}]
\efig
\end{equation}
commutes.
\item $\nabla$ a connection on $L \rightarrow \G$.  In other words a linear operator $\nabla : \Gamma^\infty ( L) \rightarrow \Omega^1 (\G , L)$ where $\Gamma^\infty (L)$ denotes sections of $L \rightarrow \G$ and $\Omega^1 (\G, L)$ denotes $L$-valued $1$-forms on $\G$.
\end{itemize}
\end{defn}

Associated to a gerbe datum $(L, \mu, \nabla)$ on an \'etale groupoid $\G$ there are 3-cocycles in $\Omega^* (\G_\bullet)$
\begin{equation}
(\alpha, \theta) \in \Omega^1(\Gl[2]) \oplus \Omega^2 (\Gl[1])
\end{equation}
where
\begin{itemize}
\item $\alpha$ is the \textit{discrepancy} between the connection $(\pr^2_{01})^* \nabla \otimes 1 + 1 \otimes (\pr^2_{12})^* \nabla$ on $(\pr^2_{01})^* L \otimes (\pr^2_{12})^* L$ and the connection $(\mu)^* (\pr^2_{02})^* \nabla$ on $(\mu)^* (\pr^2_{02})^* L$: 
\begin{equation}
\label{eq:gerbe_connection1}
(\mu)^* (\pr^2_{02})^* \nabla - \left( (\pr^2_{01})^* \nabla \otimes 1 + 1 \otimes (\pr^2_{12})^* \nabla  \right) = \alpha \in \Omega^1 (\Gl[2])
\end{equation}
\item $\theta$ is any element of $\Omega^2 (\Gl[1])$ such that 
\begin{equation}
\label{eq:gerbe_connection2}
\delta \theta = -d\alpha.
\end{equation}
In particular, $\theta = \textrm{curv} (\nabla)$, the curvature of $\nabla$, satisfies \eqref{eq:gerbe_connection2}. 
\end{itemize}
Such an $(\alpha, \theta) \in \Omega^1(\Gl[2]) \oplus \Omega^2 (\Gl[1])$ determines a class in $H^3_{dR} (\G)$, which we call the \textit{Dixmier-Douady class} of $(L,\mu,\nabla)$.  Although the Dixmier-Douady class depends on $\nabla$ throughout this article, it is expected this dependence can be removed.  This notion of a gerbe and Dixmier-Douady class agrees with the notion of a gerbe on a manifold (as in \cite{hitchin01}) when $\G$ is the \v{C}ech groupoid.

%

\subsubsection{$S^1$-gerbes over discrete translation groupoids}

As the rest of this paper is concerned with gerbes on translation groupoids arising from a discrete group action on a manifold, we will adopt a bit of new notation for this case.  Given a discrete group $\Gamma$ and a manifold $M$ carrying a (right) action of $\Gamma$, $M \times \Gamma \rightarrow M$, denoted $(x,g) \mapsto xg$, the associated translation groupoid is an \'etale groupoid and is denoted $M \rtimes \Gamma$.  The nerve of $M \rtimes \Gamma$ is a simplicial manifold $(M \rtimes \Gamma)_\bullet$ which is explicitly given on objects by $(M \rtimes \Gamma)_{(k)} = M \times \Gamma^k$ and on morphisms by
\begin{align}
\label{eq:simp_mfld_face_maps}
\delta^k_i (m, g_1 , \ldots , g_k) &= \begin{cases} (mg_1 , g_2 , \ldots , g_k) & \textrm{if } i=0 \\ (m, g_1 , \ldots , g_{i} g_{i+1} , \ldots , g_k) & \textrm{for } 1 \leq i \leq k-1 \\ (m, g_1 , \ldots , g_{k-1}) & \textrm{if } i=k \end{cases} \\
\label{eq:simp_mfld_deg_maps}
\sigma^k_j (m, g_1 , \ldots , g_k) &= (m, g_1 , \ldots , g_j , 1_\Gamma , g_{j+1} , \ldots , g_k)
\end{align}
in the notation of \eqref{eq:nerve}.  In particular, $\delta^1_0 (m, g_1) = mg_1$ and $\delta^1_1 (m , g_1) = m$.

Let $(L,\mu,\nabla)$ be a presentation of a gerbe on $M \rtimes \Gamma$.  As $(M \rtimes \Gamma)_{(1)} = M \times \Gamma$ and $\Gamma$ is a discrete group, we may view $L$ as a collection of line bundles $L_g$ on $M_g = M \times \{ g \}$ for each $g \in \Gamma$. The isomorphism $\mu$ may be restricted in the same way to a collection of isomorphisms
\begin{equation}
\mu_{g_1 , g_2} : L_{g_1} \otimes (L_{g_2})^{g_1} \xrightarrow{\sim} L_{g_1 g_2}
\end{equation}
where $g_1 , g_2 \in \Gamma$ and $(L_{g_2})^{g_1}$ denotes the line bundle $L_{g_2}$ shifted by the action of $g_1$.  In other words $(L_{g_2})^{g_1}$ has sections $s^{g_1}$ where $s^{g_1} (x) = s(xg_1)$ denotes the (left) action of $\Gamma$ for $s \in \Gamma^\infty ( L_{g_2})$.

Also denote the restriction of $\nabla$ to $L_g$ by $\nabla_g = \nabla|_{M_g}$ for each $g \in \Gamma$.   In this case the condition \eqref{eq:gerbe_connection1} may be written
\begin{equation}
\label{eq:gerbe_connection1_group}
\mu_{g_1 , g_2}^* (\nabla_{gh}) - (\nabla_g \otimes 1 + 1 \otimes (\nabla_h)^g) = \alpha (g,h)
\end{equation}
for all $g,h \in \Gamma$ and for some $\alpha (g,h) \in \Omega^1(M)$. Hence $\alpha \in \Omega^1 ( (M \rtimes \Gamma)_{(2)})$. 

Define $\theta_g \in \Omega^2 (M)$ by $[\theta_g , s] = (\nabla_g)^2 s$ for a section $s : M_g \rightarrow L_g$. Then we have a collection $\theta = (\theta_g)_{g \in \Gamma} \in  \Omega^2 ( (M \rtimes \Gamma)_{(1)} )$ and the following relation
\begin{prop}
\label{prop:gerbe_connection2_group}
The 2-forms $\theta_g$ satisfy condition \eqref{eq:gerbe_connection2}, which may be written
\begin{equation}
\label{eq:gerbe_connection2_group}
\theta_g + \theta_h^g - \theta_{gh} = -d\alpha(g,h)
\end{equation}
\begin{proof}
This follows from condition \eqref{eq:gerbe_connection1_group}.  The curvature of the sum of connections $\nabla_g \otimes 1 + 1 \otimes (\nabla_h)^g$ is $\theta_g + \theta_h^g$ and the curvature of the connection minus scalar form $\mu_{g_1 , g_2}^* (\nabla_{gh}) - \alpha (g,h)$ is $ \theta_{gh} - d\alpha(g,h)$.
\end{proof}
\end{prop}


\section{Simplicial Curvature 3-form}
\label{chap:3}

Our goal now will be to construct  a 3-form $\Theta \in \Omega^3 ((M \rtimes \Gamma)_\bullet)$ that is a representative of the Dixmier-Douady class of $(L, \mu, \nabla)$ in the complex of compatible forms on $(M\rtimes \Gamma)_\bullet$ for a fixed manifold $M$ carrying the action of a discrete group $\Gamma$.  Throughout this section, also fix the notation $(L , \mu ,\nabla)$ for a gerbe datum on $M \rtimes \Gamma$.


\subsection{Twisted bundles}

\begin{defn}
An \textit{$(L,\mu,\nabla)$-twisted vector bundle with connection} is a triple $(\E,\varphi,\nabla^\E)$ where
\begin{itemize}
\item $\E$ is a vector bundle on $M$ together with
\item $\varphi = \{\varphi_g\}_{g \in \Gamma}$ is a collection  of vector bundle isomorphisms $\varphi_g : \E \xrightarrow{\sim} \E^g \otimes L_g$ for each $g \in \Gamma$ such that the diagram
\begin{equation}
\label{eq:twisted_bundle}
\bfig
\Square|aarb|[\E^{gh} \otimes (L_h)^g \otimes L_g ` \E^{gh} \otimes L_{gh} ` \E^g \otimes L_g ` \E; \id \otimes \mu_{g,h} `  ( (\varphi_h)^g)^{-1} \otimes \id ` \varphi_{gh}^{-1} ` \varphi_g^{-1}]
\efig
\end{equation}
commutes. Here $\E^g$ denotes the vector bundle $\E$ shifted by the action of $g$.
\item $\nabla^\E : \Gamma^\infty ( \E) \to \Omega^1 (M , \E)$ is a connection on $\E$.

\end{itemize}

\end{defn}

\begin{rmk}
Note that $\varphi_g$ induces an isomorphism $\varphi_g : \End \E \xrightarrow{\sim} \End \E^g$ as well. More precisely, $\varphi_g : \E \xrightarrow{\sim} \E^g \otimes L_g$ induces an isomorphism
\begin{equation*}
\varphi_g : \End \E \xrightarrow{\sim} \End (\E^g \otimes L_g) \cong \End (\E^g ) \otimes \End (L_g).
\end{equation*}
As $\End (L_g) \cong \C$ we see $\End \E \cong \End \E^g$ via $\varphi_g$.
\end{rmk}

Let $\theta = \textrm{curv} (\nabla)$ and denote the discrepancy of $(L, \mu , \nabla)$ by $\alpha$.  Then there is a connection on $\E^g \otimes L_g$ given by ${(\nabla^\E)}^g \otimes 1 + 1 \otimes \nabla_g$ for every $g \in \Gamma$.  By ${(\nabla^\E)}^g$ we mean the connection ${(\nabla^\E)}^g s^g = {(\nabla^\E s)}^g$ where $s^g$ is a section of $\E^g$.  This defines a connection on $\E^g$ which we denote $\nabla^{\E^g} = {(\nabla^\E)}^g$.  The discrepancy between $\nabla^\E$ and $\nabla^{\E^g} \otimes 1 + 1 \otimes \nabla_g$ is an $\End \E$-valued 1-form which we denote by
\begin{equation}
 \varphi^*_g ( \nabla^{\E^g} \otimes 1 + 1 \otimes \nabla_g)  - \nabla^\E = A(g) \in \Omega^1 (M , \End \E)
\end{equation}
for each $g \in \Gamma$. We will call $A \in \Omega^1 (M \times \Gamma , \End \E)$ the \textit{discrepancy} of $\nabla^\E$. Notice that we then have the identity
\begin{align}
\label{eq:twisted_bundle_curvature}
(\nabla^\E + A(g))^2 &= \varphi_g^* ((\nabla^{\E^g}) \otimes 1 + 1 \otimes \nabla_g)^2 \notag \\
 (\nabla^\E)^2 + A(g)^2 + [\nabla^\E , A(g)] &= \varphi_g^* ((\nabla^{\E^g})^2 ) + \theta_g \notag \\
 \theta^\E + A(g)^2 + [\nabla^\E , A(g)] &= \varphi_g^* (\theta^\E)^g + \theta_g 
\end{align}
on $\Omega^2 (M , \End \E)$ where we define $\theta^\E = (\nabla^\E)^2$.


\begin{ex}[Direct Sum Bundle]
\label{ex:direct_sum_bundle}

There is a particular $(L, \mu, \nabla)$-twisted vector bundle called the \textit{direct sum bundle} of $(L, \mu ,\nabla)$ which is the most relevant example of a twisted vector bundle for our purposes. The direct sum bundle $(\E , \varphi, \nabla^\E)$ of $(L, \mu ,\nabla)$ consists of
\begin{itemize}
\item the (possibly infinite dimensional) vector bundle $\pi : \E \rightarrow M$ where $\E = \bigoplus_{g \in \Gamma} L_g$,
\item a collection of isomorphisms $\varphi= \{\varphi_g\}_{g \in \Gamma}$ with $\varphi_g : \E \to \E^g \otimes L_g$ for each $g \in \Gamma$ defined by
\begin{equation}
\varphi_g : \bigoplus_{g' \in \Gamma} L_{g'} \to \bigoplus_{g' \in \Gamma} \mu_{g , g^{-1}g'}^{-1} (L_{g'}) = \bigoplus_{g' \in \Gamma} (L_{g^{-1}g'})^g \otimes L_g
\end{equation}
where the factor $L_{g'}$ is mapped to $\mu_{g , g^{-1}g'}^{-1} (L_{g'})$ under $\varphi_g$.
\item the direct sum connection $\nabla^\E = \bigoplus_{g \in \Gamma} \nabla_g$ on $\E \rightarrow M$ induced by $\nabla$.
\end{itemize}

\begin{lem}
\label{lem:discrepancy_identity}
Let $(\E , \varphi , \nabla)$ be the direct sum bundle of $(L,\mu, \nabla)$.  Then for any $g, h \in \Gamma$ we have the identity
\begin{equation}
A(g) + (\varphi_g)^* (A(h)^g) - A(gh) =  \alpha (g,h) \in \Omega^1 (M , \End \E),
\end{equation}
where $\alpha$ is the discrepancy of $\nabla$ and $A$ is the discrepancy of $\nabla^\E$.
\begin{proof}
This follows from the commutativity of \eqref{eq:twisted_bundle}.
\end{proof}
\end{lem}

\begin{lem}
\label{cor:discrepancy_identity2}
Let $(\E , \varphi , \nabla)$ be the direct sum bundle of $(L,\mu, \nabla)$. For any $g_1 , \ldots , g_j \in \Gamma$ and $1 < i < j$,
\begin{equation}
\alpha (g_1 , g_2 \cdots g_i) - \alpha (g_1 , g_2 \cdots g_j) + \alpha (g_1 \cdots g_i , g_{i+1} \cdots g_j) \\=  \varphi_{g_1}^* \alpha (g_2 \cdots g_i , g_{i+1} \cdots g_j)^{g_1} ,
\end{equation}
where $\alpha$ is the discrepancy of $\nabla$.
\begin{proof}
This follows from a direct calculation using Lemma \ref{lem:discrepancy_identity}.
\end{proof}
\end{lem}

The bundle $\End \E \rightarrow M$ is vector bundle of finite rank endomorphisms when $(\E, \varphi, \nabla^\E)$ is a direct sum bundle of $(L, \mu ,\nabla)$.  Note that any section of $\Gamma^\infty ( \End \E)$ may be decomposed as follows.  Since $\E = \bigoplus_{g \in \Gamma} L_g$, $\End \E$ consists of morphisms $ \phi : \bigoplus_{g \in \Gamma} L_g \rightarrow \bigoplus_{g \in \Gamma} L_g$.  Let $E_{g , h} (\phi) \in \Hom (L_g , L_h)$ denote the restriction of $\phi$ to the component $E_{g,h} (\phi): L_g \rightarrow L_h$, so that $\phi = \bigoplus_{g,h \in \Gamma} E_{g,h} (\phi)$.  So given a section $f \in \Gamma^\infty ( \End \E)$ we may decompose $f$ in the same way, with $E_{g,h} (f)$ a section on $M$ with values in morphisms $L_g \rightarrow L_h$. In other words, $E_{g,h} (f) \in \Gamma^\infty ( \Hom (L_g , L_h))$.

Given such a decomposition we can describe the action of $g \in \Gamma$ on sections $\Gamma^\infty (\End \E)$ as follows.  Let $f \in \Gamma^\infty (\Hom(L_{g_1} , L_{g_2}) )$, i.e.\ $f=E_{g_1 , g_2} (f)$.  For any $g \in \Gamma$, the section $f^g$ is a section of $\Hom (L_{g_1}^g , L_{g_2}^g)$ and $\id \otimes f^g$ is a section of $\Hom (L_{g} \otimes (L_{g_1})^{g} , L_{g} \otimes (L_{g_2})^{g})$.  We may denote the corresponding section by $f^{g}$ as well.  As there are isomorphisms $\mu_{g, g_1} : L_{g} \otimes (L_{g_1})^{g} \rightarrow L_{g g_1}$ and $\mu_{g , g_2} : L_{g} \otimes (L_{g_2})^{g} \rightarrow L_{g g_2}$, under these isomorphisms $f^{g}$ corresponds to a section in $\Gamma^\infty (\Hom (L_{g g_1} , L_{g g_2}) )$ which is denoted by $E_{g g_1 , g g_2} (f^g)$.  In other words
\begin{equation*}
E_{g g_1 , g g_2} (f^g) = \mu_{g , g_2} \circ \id \otimes f^g \circ \mu_{g , g_1}^{-1}
\end{equation*}
Hence the (left) action of $g \in \Gamma$ on $E_{g_1 , g_2} (f)$ is $g \cdot E_{g_1 , g_2} (f) = E_{g g_1 , g g_2} (f^g)$.

For a section $a \in \Gamma^\infty ( L_g)$, we may consider $a$ as the section $E_{1, g} ( a) : L_1 \rightarrow L_g$ of $\Gamma^\infty ( \End \E)$ since $a : M_g \rightarrow L_g$ defines a section $\tilde{a} : M \times \C \rightarrow L_g$ by $\tilde{a} ( x, \lambda) = a(x)$ for all $\lambda \in \C$ and $L_1 \cong M \times \C$. We have $L_1 \cong M \times \C$ as $\mu_{1 , 1} : L_1 \otimes L_1 \xrightarrow{\sim} L_1$.  
We will not make a distinction in notation between a section $a : M \rightarrow L_g$ and the corresponding section of $ \Gamma^\infty (\Hom (L_1 , L_g ) )$.

Let $f_1 , f_2 \in \Gamma^\infty ( \End \E)$ be sections of $\End \E$. The sections $E_{g_1 , g_2} (f_1)$ and $E_{g_2 , g_3} (f_2)$ may be composed to produce a section $E_{g_2 , g_3} (f_2) \circ E_{g_1, g_2} (f_1)$. We will denote this operation of composition by
\begin{equation}
E_{g_1, g_2} (f_1) E_{g_2 , g_3} (f_2) = E_{g_1 , g_3} ( f_1 f_2 )  \in \Gamma^\infty (\Hom (L_{g_1},  L_{g_3}))
\end{equation}
Thus we may obtain a matrix product on $\Gamma^\infty ( \End \E)$.  Let $f_1 f_2 \in \Gamma^\infty ( \End \E)$ be defined by
\begin{equation}
E_{h , k} (f_1 f_2) = \sum_{g \in \Gamma} E_{h , g} (f_1) E_{g , k} (f_2)
\end{equation}
for $h,k \in \Gamma$ so that $f_1 f_2 = \bigoplus_{h,k \in \Gamma} E_{h,k} (f_1 f_2)$.

\begin{lem}

Given sections $a_g \in \Gamma^\infty ( L_g)$ and $a_h \in \Gamma^\infty ( L_h)$
\begin{equation*}
E_{1,g} (a_g) E_{g,h} (a_h^g) = E_{1, gh} (\mu_{g,h} (a_g \otimes a_h^g) )
\end{equation*}
where we consider  $a_g \in \Gamma^\infty (\Hom (L_1 , L_g))$, $a_h \in \Gamma^\infty (\Hom (L_1 , L_h))$, as well as $\mu_{g,h} (a_g \otimes a_h^g) \in \Gamma^\infty ( \Hom (L_1 , L_{gh}) )$.
\begin{proof}
The left hand side is defined by
\begin{equation*}
L_1 \xrightarrow{a_g} L_g \xrightarrow{\mu_{g,1}^{-1} }  L_g \otimes L_1^g \xrightarrow{\id \otimes a_h^g} L_g \otimes L_h^g  \xrightarrow{\mu_{g,h}} L_{gh}
\end{equation*}
while the right hand side is defined by
\begin{equation*}
L_1 \xrightarrow{ a_g \otimes a_h^g } L_g \otimes L_h^g \xrightarrow{ \mu_{g,h} } L_{gh}
\end{equation*}
which are corresponding sections of $\Hom (L_1 , L_{g,h})$ under the isomorphism $\mu$.
\end{proof}
\end{lem}

Let $(\E , \varphi , \nabla)$ be the direct sum bundle of $(L,\mu, \nabla)$.  For $f \in \Gamma^\infty ( \End \E)$ there is a trace map $\tr : \Gamma^\infty ( \End \E) \rightarrow C^\infty (M)$ defined by
\begin{equation}
\tr(f) = \sum_{g \in \Gamma} E_{g , g} (f)
\end{equation}
which is an element of $C^\infty (M)$ as $E_{g, g} (f) \in \End (L_{g})$ for all $g \in \Gamma$.  In fact for $\Omega^* (M, \End \E) = \Omega^* (M) \otimes \Gamma^\infty ( \End \E)$ this trace map induces a trace
\begin{equation}
\tr : \Omega^* (M , \End \E) \rightarrow \Omega^* (M)
\end{equation}
as for $\omega \otimes f \in \Omega^*(M) \otimes \Gamma^\infty (\End \E)$ we may define $\tr (\omega \otimes f) = \omega \tr(f)$.

\begin{lem}
\label{lem:trace_commutator}
Let $(\E , \varphi , \nabla)$ be the direct sum bundle of $(L,\mu, \nabla)$. For $\omega_1, \omega_2 \in \Omega^* (M , \End \E)$
\begin{equation}
\tr (\omega_1 \omega_2) = (-1)^{ |\omega_1| |\omega_2|} \tr(\omega_2 \omega_1),
\end{equation}
and furthermore,
\begin{equation}
\tr ( [\omega_1 , \omega_2]) = 0,
\end{equation}
where $[\omega_1 , \omega_2]$ is the graded commutator of $\omega_1$ and $\omega_2$.
\begin{proof}
We will first mention that $\tr (f_1 f_2) = \tr (f_2 f_1)$ for any $f_1 , f_2 \in \Gamma^\infty (\End \E)$. This follows as usual
\begin{align*}
\tr(f_1 f_2) &= \sum_{g \in \Gamma} E_{g,g} (f_1 f_2) \\
&= \sum_{g \in \Gamma}  \left( \sum_{h \in \Gamma} E_{g , h} (f_1) E_{h,g} (f_2) \right) \\
&= \sum_{h \in \Gamma} \left( \sum_{g \in \Gamma}  E_{h,g} (f_2) E_{g,h} (f_1) \right) \\
&= \sum_{h \in \Gamma} E_{h,h} (f_2 f_1) \\
&= \tr(f_2 f_1)
\end{align*}
by changing the summation order.

Let $\omega_1 = \alpha_1 \otimes f_1 $, $\omega_2 = \alpha_2 \otimes f_2$ where $\alpha_1 , \alpha_2 \in \Omega^* (M)$ and $f_1, f_2 \in \Gamma^\infty (\End \E)$.  Then $| \alpha_1 | = | \omega_1|$ and $| \alpha_2 | = | \omega_2 |$ and
\begin{align*}
\tr (\omega_1 \omega_2) &= \tr ((\alpha_1 \wedge \alpha_2) \otimes (f_1 f_2) ) \\
&= (\alpha_1 \wedge \alpha_2) \tr ( (f_1 f_2) ) \\
&= (-1)^{|\alpha_1| |\alpha_2|} (\alpha_2 \wedge \alpha_1) \tr(  (f_2 f_1) ) \\
&= (-1)^{|\omega_1| |\omega_2|} \tr(\omega_2 \omega_1)
\end{align*}
From this it follows that
\begin{equation*}
\tr ([\omega_1, \omega_2]) = \tr(\omega_1 \omega_2 - (-1)^{|\omega_1| |\omega_2|} \omega_2 \omega_1) =  \tr(\omega_1 \omega_2) - (-1)^{|\omega_1| |\omega_2|} \tr( \omega_2 \omega_1) = 0.
\end{equation*}
\end{proof}
\end{lem}

\end{ex}


\subsection{Derivations $\nabla^k$}

We will again denote the discrepancy of  $(L, \mu, \nabla)$ by $\alpha$.  For the rest of this section, fix the notation $(\E , \varphi , \nabla^\E)$ for the direct sum bundle of $(L , \mu , \nabla)$ as in Example \ref{ex:direct_sum_bundle}.

Let $p_k : (M\rtimes \Gamma)_{(k)} = M \times \Gamma^k \rightarrow M$ be defined by $p_k : (m, g_1 , \ldots , g_k) \mapsto m \in M$.  The pullback bundle $\E_k = p_k^* \E$ is a vector bundle over $(M\rtimes \Gamma)_{(k)}$.  There is also a connection on $\E_k$ defined by $\nabla^{\E_k} = p_k^* \nabla^{\E}$.  We will denote the connection $\nabla^{\E_k} |_{M \times \{g_1 \} \times \cdots \times \{g_k\}}$ by $\nabla^{\E_k} (g_1 , \ldots , g_k)$.  Note with this definition that $\nabla^{\E_k} (g_1 , \ldots , g_k) = \nabla^{\E}$ for all $g_1 , \ldots , g_k$.

\begin{rmk}
In general the collection $\{ \pi_k : \E_k \rightarrow (M \rtimes \Gamma)_{(k)} \}$ cannot be made into simplicial vector bundle using the maps \eqref{eq:simp_mfld_face_maps} and \eqref{eq:simp_mfld_deg_maps}.  In order for the diagram \eqref{eq:simp_bund_diag} of face maps to commute, we must have that $\delta_i^* \E_{k-1} \cong \E_{k}$ for $0 \leq i \leq k$ which does not hold in general.  Specifically, consider $\delta_0^* \E_{k-1}$.  For any $g_1 , \ldots , g_k \in \Gamma$ and $m \in M$,
\begin{align*}
\delta_0^* \E_{k-1} (m , g_1 , \ldots ,g_k) &= \delta_0^* (p_{k-1}^* \E) (m , g_1 ,\ldots , g_k) \\
&= (p_{k-1}^* \E) (mg_1 , g_2 ,\ldots , g_k) \\
&= \E^{g_1} \\
&\cong \E \otimes L_{g_1^{-1}},
\end{align*}
where the last isomorphism is given by $\varphi_{g_1^{-1}}$.
\end{rmk}

On the other hand, the vector bundles $\End \E_k \rightarrow (M \rtimes \Gamma)_{(k)}$ do form a simplicial vector bundle.

\begin{prop}
The vector bundles $\pi_{k} : \End \E_{k} \rightarrow (M \rtimes \Gamma)_{(k)}$ form a simplicial vector bundle $\pi_\bullet : \End \E_\bullet \rightarrow (M \rtimes \Gamma)_\bullet$.
\begin{proof}
This follows as $\End (\E \otimes L_g) \cong \End \E$.  So $\delta_i^* \End \E_{k-1} \cong \End \E_k$ and the diagram \eqref{eq:simp_bund_diag} of face maps commutes for $0 \leq i \leq k$.
\end{proof}
\end{prop}

\begin{defn}
Denote $(M \rtimes \Gamma)_{(k)}^\Delta = (M \rtimes \Gamma)_{(k)} \times \Delta^k$ and $\E_k^\Delta = \E_k \times \Delta^k$.  Consider the bundle $ \pi_k \times \id :   \E_k^\Delta \rightarrow  (M \rtimes \Gamma)_{(k)}^\Delta $. Then define a connection 
\begin{equation}
\nabla^{k} : \Gamma^\infty (\E_k^\Delta) \rightarrow \Omega^1 ((M \rtimes \Gamma)_{(k)}^\Delta , \E_k^\Delta)
\end{equation}
by the formula 
\begin{multline}
\label{eq:derivation_def}
((\nabla^{k} (g_1 , \ldots , g_k) ) s) ( x_1 , \ldots , x_m, t_1 , \ldots , t_k ) \\
= ((\nabla^{\E_k} (g_1 , \ldots , g_k) + t_1 A(g_1) + \cdots + t_k A(g_1 \cdots g_k)) s ) ( x_1 , \ldots , x_m, t_1 , \ldots , t_k )
\end{multline}
where $s \in \Gamma^\infty (\E_k^\Delta)$, $x_1 , \ldots , x_m$ are coordinates on $M$ (with $m = \dim M$), and $t_1 , \ldots , t_k$ are Cartesian coordinates on $\Delta^k$.

Here we use $\nabla^{\E_k}$ to denote the derivation on $\E_k^\Delta \rightarrow  (M \rtimes \Gamma)_{(k)}^\Delta$ that is the derivation $\nabla^{\E_k}$ on $\E_k \rightarrow (M \rtimes \Gamma)_{(k)}$ at every point of $\Delta^k$.  In other words, if $\bar{p}_k : (M \rtimes \Gamma)_{(k)}^\Delta \rightarrow (M \rtimes \Gamma)_{(k)}$ is the projection $(x,t) \mapsto x$ then we denote $\bar{p}_k^* \nabla^{\E_k}$ merely by $\nabla^{\E_k}$.  Note that, while $\nabla^k$ is a connection on $\E_k^\Delta \rightarrow  (M \rtimes \Gamma)_{(k)}^\Delta$, $\nabla^k$ is not a simplicial connection as $\E_k \rightarrow (M \rtimes \Gamma)_{(k)}$ does not form a simplicial bundle.

We may use the formal variable $u$ to rescale the connection $\nabla^k$ by considering $(\nabla^k)^{1,0}$ the component of the operator $\nabla^k$ 
\begin{equation*}
(\nabla^k)^{1,0} : \Gamma^\infty (\E_k^\Delta) \rightarrow \Omega^1 ((M \rtimes \Gamma)_{(k)} , \E_k) \otimes C^\infty (\Delta^k , \Delta^k)
\end{equation*}
and $(\nabla^k)^{0,1} = d_\Delta$ the component that maps
\begin{equation*}
(\nabla^k)^{0,1} : \Gamma^\infty (\E_k^\Delta) \rightarrow \Gamma^\infty (\E_k) \otimes \Omega^1 (\Delta^k, \Delta^k).
\end{equation*}
Let $\nabla^k_u = (\nabla^k)^{1,0} + u^{-1} d _\Delta$ so that
\begin{equation}
\nabla^k_u :  \Gamma^\infty (\E_k^\Delta) [u] \to \Omega^1 ((M \rtimes \Gamma)_{(k)}^\Delta , \E_k^\Delta) [u , u^{-1}]
\end{equation}

We may extend $\nabla^k$ to an operator
\begin{equation}
\nabla^k : \Omega^\ell ((M \rtimes \Gamma)_{(k)}^\Delta , \E_k^\Delta) \rightarrow \Omega^{\ell + 1} ((M \rtimes \Gamma)_{(k)}^\Delta , \E_k^\Delta)
\end{equation}
in the usual way.  Given a section $s \in \Gamma^\infty ( \E_k^\Delta)$ and a form $\omega \in \Omega^\ell ((M \rtimes \Gamma)_{(k)}^\Delta)$ we define
\begin{equation}
\nabla^k (\omega \otimes s) = d \omega \otimes s + (-1)^\ell \omega \wedge \nabla^k s.
\end{equation}

\end{defn}

\begin{rmk}
The operator $\nabla^k$ may be defined to act on the algebra of $\End \E$-valued forms $\Omega^* (M , \End \E)$ and in particular sections $\Gamma^\infty (\End \E)$ as follows.  Given $\eta \in \Omega^* (M , \End \E)$, define
\begin{equation}
\nabla^k \eta = [\nabla^k , \bar{p}_k^* p_k^* \eta]
\end{equation}
where 
\begin{equation}
p_k : (M \rtimes \Gamma)_{(k)} \rightarrow M
\end{equation}
is the projection $(x, g_1 , \ldots , g_k) \mapsto x$, 
\begin{equation}
\bar{p}_k : (M \rtimes \Gamma)_{(k)} \times \Delta^k \rightarrow (M \rtimes \Gamma)_{(k)}
\end{equation}
is the natural projection $(x,t) \mapsto x$ and we consider $\nabla^k$ to be degree $+1$ with respect to the graded commutator. In other words, $\nabla^k \eta$ operates on sections $s \in \Gamma^\infty (\E_k^\Delta)$ by
\begin{equation}
(\nabla^k \eta) s = [\nabla^k , \bar{p}_k^* p_k^* \eta] s = \nabla^k ((\bar{p}_k^* p_k^* \eta) s) - (-1)^{|\eta| + 1} (\bar{p}_k^* p_k^* \eta) (\nabla^k s).
\end{equation}
\end{rmk}

With this extension the rescaling $\nabla^k_u =  (\nabla^k)^{1,0} + u^{-1} d_\Delta$ still makes sense.

\begin{lem}
\label{lem:trace_connection}
For any $\eta \in \Omega^* (M, \End \E)$,
\begin{equation}
(u d_{dR} + d_\Delta) \tr (\eta) = u\tr (\nabla^k_u (\eta))
\end{equation}
where by $\tr (\eta)$ we mean $\tr(\bar{p}_k^* p_k^* \eta)$.
\begin{proof}
We will omit the pullback maps $\bar{p}_k^* p_k^*$ as they should be clear from the context.  As $\nabla^k_u =  (\nabla^k)^{1,0} + u^{-1}d_\Delta$ and $(\nabla^k)^{1,0} = d_{dR} + \omega$ locally for some $\omega \in \Omega^1 (M, \End \E)$.  For a local section $s \in \Gamma^\infty (\E_k^\Delta)$
\begin{align*}
((\nabla^k)^{1,0} \eta) (s) &= [(\nabla^k)^{1,0} , \eta] s \\
&= [d_{dR}  + \omega , \eta] s \\
&= d_{dR} (\eta s) - (-1)^{|\eta|} \eta d_{dR} s + [\omega , \eta] s \\
&= (d_{dR} \eta + [\omega , \eta]) s
\end{align*}
So by Lemma \ref{lem:trace_commutator}
\begin{align*}
\tr ((\nabla^k)^{1,0} \eta) &= \tr ( (d_{dR} \eta +  [\omega , \eta])) \\
&= \tr (d_{dR} \eta) + \tr([\omega , \eta]) \\
&= \tr (d_{dR} \eta) \\
&= d_{dR} \tr (\eta).
\end{align*}
Hence
\begin{equation*}
u\tr ( ( (\nabla^k)^{1,0} + u^{-1}d_\Delta) \eta) = u \tr ( (\nabla^k)^{1,0} \eta) + \tr (d_\Delta \eta) = u d_{dR} \tr(\eta) + d_{\Delta} \tr (\eta),
\end{equation*}
as required. 
\end{proof}
\end{lem}


\subsection{Simplicial 2-form}

The collection of $\End \E_{k}$-valued $2$-forms on $(M \rtimes \Gamma)_{(k)} \times \Delta^k$ given by $(\nabla^k)^2$ for each $k \geq 0$ does not define a simplicial differential form on $(M \rtimes \Gamma)_\bullet$.  In the theorem below we adjust $(\nabla^k)^2$ by a scalar-valued form to obtain a simplicial differential $2$-form on $M_\bullet$.

\begin{thm}
\label{thm:simplicial_2_form}
For each $k \geq 0$, let $\vartheta_{(k)} \in \Omega^2 ( (M \rtimes \Gamma)_{(k)}^\Delta, \End \E_{k}^\Delta)$ be defined by the formula
\begin{multline}
\label{eq:theta_def}
\vartheta_{(k)} (g_1 , \ldots , g_k) = (\nabla^k (g_1, \ldots , g_k))^2  \\ - \sum_{i=1}^k t_i \theta_{g_1 \cdots g_i}  - \sum_{1 \leq i < j \leq k} \alpha(g_1 \cdots g_i , g_{i+1} \cdots g_j) (t_i dt_j - t_j dt_i).
\end{multline}
where $\alpha$ denotes the discrepancy of $(L,\mu,\nabla)$. Then $\vartheta = \{ \vartheta_{(k)} \}$ is a compatible $2$-form on $(M \rtimes \Gamma)_\bullet$ with values in $\End \E_\bullet$, which means that $\vartheta \in \Omega^2 ((M \rtimes \Gamma)_\bullet , \End \E_\bullet)$. We call $\vartheta$ the \emph{simplicial $2$-form associated to $(L , \mu , \nabla)$}.
\begin{proof}
Our goal is to check the compatibility conditions of \eqref{eq:compatibility_cond}. First it will be helpful to rewrite $\vartheta_{(k)} (g_1 , \ldots , g_k)$. By expanding $(\nabla^k)^2$ using the definition of $\nabla^k$ in \eqref{eq:derivation_def},
\begin{align*}
(\nabla^k (g_1, \ldots , g_k))^2 &= (\nabla^{\E})^2 + \sum_{i=1}^k t_i^2 A(g_1 \cdots g_i)^2 +  [\nabla^{\E}  , t_i A(g_1 \cdots g_i) ]  \\
& \quad + \sum_{1 \leq i < j \leq k} [t_i A(g_1 \cdots g_i) , t_j A(g_1 \cdots g_j)] \\
\intertext{and by expanding $[ \nabla^{\E}  , t_i A(g_1 \cdots g_i) ] = dt_i A(g_1 \cdots g_i) +  t_i [ \nabla^{\E}  , A(g_1 \cdots g_i) ]$,}
&= \theta^\E  + \sum_{i=1}^k t_i^2 A(g_1 \cdots g_i)^2 + dt_i A(g_1 \cdots g_i) + t_i [ \nabla^{\E}  , A(g_1 \cdots g_i) ] \\
&\quad + \sum_{1 \leq i < j \leq k} [t_i A(g_1 \cdots g_i) , t_j A(g_1 \cdots g_j)] \\
\intertext{which, by substituting into \eqref{eq:theta_def}, gives}
 \vartheta_{(k)}(g_1 , \ldots , g_k) &=  \theta^\E  + \sum_{i=1}^k t_i^2 A(g_1 \cdots g_i)^2 + dt_i A(g_1 \cdots g_i) + t_i [ \nabla^{\E}  , A(g_1 \cdots g_i) ]  - t_i \theta_{g_1 \cdots g_i} \displaybreak[0]\\
& \quad + \sum_{1 \leq i < j \leq k} [t_i A(g_1 \cdots g_i) , t_j A(g_1 \cdots g_j)] - \alpha(g_1 \cdots g_i , g_{i+1} \cdots g_j) (t_i dt_j - t_j dt_i) .\displaybreak[2]\\
\intertext{By adding and subtracting $\displaystyle \sum_{i = 1}^k t_i \theta^\E $, factoring $t_i$ from some terms, and using $t_i^2 = t_i - (1-t_i)t_i$, }
 \vartheta_{(k)}(g_1 , \ldots , g_k) &= \left(1 - \sum_{i = 1}^k t_i  \right) \theta^\E  + \sum_{i=1}^k  t_i \left( \theta^\E +  A(g_1 \cdots g_i)^2 + [ \nabla^{\E}  , A(g_1 \cdots g_i) ]  - \theta_{g_1 \cdots g_i} \right) \displaybreak[0]\\
&\qquad \qquad + dt_i A(g_1 \cdots g_i) - (1-t_i)t_i A(g_1 \cdots g_i)^2 \displaybreak[0]\\
& \quad + \sum_{1 \leq i < j \leq k} [t_i A(g_1 \cdots g_i) , t_j A(g_1 \cdots g_j)] - \alpha(g_1 \cdots g_i , g_{i+1} \cdots g_j) (t_i dt_j - t_j dt_i) 
\end{align*}
So by Equation \eqref{eq:twisted_bundle_curvature} we have obtained
\begin{multline}
\label{eq:theta_rewritten}
\vartheta_{(k)}(g_1 , \ldots , g_k) \\ =  \left(1 - \sum_{i = 1}^k t_i  \right) \theta^\E  +  \sum_{i=1}^k  t_i \varphi_{g_1\cdots g_i}^* \theta^{\E^{g_1 \cdots g_i}} + dt_i A(g_1 \cdots g_i) - (1-t_i)t_i A(g_1 \cdots g_i)^2 \\+ \sum_{1 \leq i < j \leq k} [t_i A(g_1 \cdots g_i) , t_j A(g_1 \cdots g_j)] - \alpha(g_1 \cdots g_i , g_{i+1} \cdots g_j) (t_i dt_j - t_j dt_i) 
\end{multline}

As we are using Cartesian coordinates, when $1 \leq \ell \leq k$ we have that 
\begin{equation}
( \id \times \partial_\ell )^* \vartheta_{(k)} (g_1 , \ldots , g_k) = \vartheta_{(k)} (g_1 , \ldots , g_k) |_{t_\ell = 0}
\end{equation}
and when $\ell = 0$ we have that 
\begin{equation}
( \id \times \partial_0 )^* \vartheta_{(k)} (g_1 , \ldots , g_k) = \vartheta_{(k)} (g_1 ,\ldots , g_k) |_{t_1+ \cdots + t_k = 1}
\end{equation}
Therefore, to show that $( \id \times \partial_\ell )^* \vartheta_{(k)} = (\delta_\ell \times \id)^* \vartheta_{(k-1)}$ for $1 \leq \ell \leq k-1$, we explicitly show
\begin{equation}
 \vartheta_{(k)} (g_1 , \ldots , g_k) |_{t_\ell = 0} = \vartheta_{(k-1)} (g_1 , \ldots , g_\ell g_{\ell +1} , \ldots , g_k)
\end{equation}
Note that when $t_\ell = 0$, $dt_\ell = 0$.  From the above
\begin{align*}
 \vartheta_{(k)} (g_1 , \ldots , g_k) |_{t_\ell = 0} &= \theta^\E + \sum_{\substack{ i=1 \\ i \neq \ell} }^k t_i^2 A(g_1 \cdots g_i)^2 + [ \nabla^{\E}  , t_i A(g_1 \cdots g_i) ] - t_i \theta_{g_1 \cdots g_i} \displaybreak[0]\\
& \quad + \sum_{\substack{ 1 \leq i < j \leq k \\ i,j \neq \ell} } [t_i A(g_1 \cdots g_i) , t_j A(g_1 \cdots g_j)] \\&\qquad \qquad - \alpha(g_1 \cdots g_i , g_{i+1} \cdots g_j) (t_i dt_j - t_j dt_i) \displaybreak[2]
\end{align*}
which may be written as $\vartheta_{(k-1)} (g_1 , \ldots , g_\ell g_{\ell +1} , \ldots , g_k)$ upon reindexing $t_{i+1}$ to $t_i$ when $\ell \leq i+1 \leq k$.

Now for the $\ell = 0$ case, we explicitly show that
\begin{equation}
\vartheta_{(k)} (g_1 ,\ldots , g_k) |_{t_1+ \cdots + t_k = 1} = \varphi_{g_1}^* \vartheta_{(k-1)} (g_2 , \ldots , g_k)^{g_1}
\end{equation}
Note that when $t_1 + \cdots + t_k = 1$ we have $dt_1 + \cdots + dt_k = 0$. Before handling the $\ell = 0$ case, it will be helpful to do a couple side calculations.
\begin{lem}
\label{lem:2-form_lemma_1}
With the hypothesis of the theorem and assuming $t_1 + \cdots + t_k = 1$,
\begin{multline}
 \sum_{1 \leq i < j \leq k} (t_i dt_j - t_j dt_i)  \alpha(g_1 \cdots g_i , g_{i+1} \cdots g_j) = \sum_{j=2}^k dt_j \alpha(g_1 , g_2 \cdots g_j) \\
 \quad + \sum_{2 \leq i < j \leq k} (t_i dt_j - t_j dt_i)  (\alpha(g_1 , g_2 \cdots g_i) - \alpha( g_1 , g_2 \cdots g_j) + \alpha(g_1 \cdots g_i , g_{i+1} \cdots g_j)) 
\end{multline}
\begin{proof}
This is just a rearrangement of the summation:
\begin{align*}
& \sum_{1 \leq i < j \leq k} (t_i dt_j - t_j dt_i)  \alpha(g_1 \cdots g_i , g_{i+1} \cdots g_j) =\\
&= \sum_{2 \leq i < j \leq k} (t_i dt_j - t_j dt_i)  \alpha(g_1 \cdots g_i , g_{i+1} \cdots g_j) \\
& \quad + \sum_{j=2}^k (t_1 dt_j - t_j dt_1) \alpha (g_1 , g_2 \cdots g_j) \\
&= \sum_{2 \leq i < j \leq k} (t_i dt_j - t_j dt_i)  \alpha(g_1 \cdots g_i , g_{i+1} \cdots g_j) \\
& \quad + \sum_{j=2}^k \left( \left(1-\sum_{\ell=2}^k t_\ell \right) dt_j - t_j \left( \sum_{\ell=2}^k (-dt_\ell) \right) \right) \alpha (g_1 , g_2 \cdots g_j) \\
&= \sum_{2 \leq i < j \leq k} (t_i dt_j - t_j dt_i)  \alpha(g_1 \cdots g_i , g_{i+1} \cdots g_j) \\
& \quad + \sum_{j,\ell=2, j \neq \ell}^k (t_j dt_\ell - t_\ell dt_j) \alpha (g_1 , g_2 \cdots g_j) + \sum_{j=2}^k dt_j \alpha(g_1 , g_2 \cdots g_j) \\
&= \sum_{j=2}^k dt_j \alpha(g_1 , g_2 \cdots g_j) \\
& \quad + \sum_{2 \leq i < j \leq k} (t_i dt_j - t_j dt_i)  (\alpha(g_1 , g_2 \cdots g_i) - \alpha( g_1 , g_2 \cdots g_j) + \alpha(g_1 \cdots g_i , g_{i+1} \cdots g_j)) 
\end{align*}
\end{proof}
\end{lem}

\begin{lem}
\label{lem:2-form_lemma_2}
With the hypothesis of the theorem and assuming $t_1 + \cdots + t_k = 1$,
\begin{multline}
\sum_{1 \leq i < j \leq k} - t_i t_j \left( \varphi_{g_1 \cdots g_i}^* A(g_{i+1} \cdots g_j)^{g_1 \cdots g_i} \right)^2 \\
=\sum_{i=2}^k - \left( 1 -  t_i \right) t_i (\varphi_{g_1}^* A(g_2 \cdots g_i)^{g_1})^2 + \sum_{2 \leq i < j \leq k} - t_i t_j  [\varphi_{g_1}^* A(g_2 \cdots g_i)^{g_1} , \varphi_{g_1}^* A(g_2 \cdots g_j)^{g_1} ] 
\end{multline}
\begin{proof}
We extract the $i=1$ summands into a separate summation and rearrange
\begin{align*}
& \sum_{1 \leq i < j \leq k} - t_i t_j \left( \varphi_{g_1 \cdots g_i}^* A(g_{i+1} \cdots g_j)^{g_1 \cdots g_i} \right)^2 \\
&= \sum_{j=2}^k - \left( 1 - t_2 - \cdots - t_k \right) t_j (\varphi_{g_1}^* A(g_2 \cdots g_j)^{g_1})^2  + \sum_{2 \leq i < j \leq k} - t_i t_j \left( \varphi_{g_1 \cdots g_i}^* A(g_{i+1} \cdots g_j)^{g_1 \cdots g_i} \right)^2 \\
&= \sum_{j=2}^k - \left( 1 -  t_j \right) t_j (\varphi_{g_1}^* A(g_2 \cdots g_j)^{g_1})^2 \displaybreak[0]\\
& \quad + \sum_{2 \leq i < j \leq k} - t_i t_j \left( \varphi_{g_1 \cdots g_i}^* (A(g_{i+1} \cdots g_j)^{g_1 \cdots g_i})^2  - (\varphi_{g_1}^* A(g_2 \cdots g_i)^{g_1})^2 - (\varphi_{g_1}^* A(g_2 \cdots g_j)^{g_1})^2 \right).  \displaybreak[0]\\
\intertext{By adding and subtracting the same term in the second summation}
& = \sum_{j=2}^k - \left( 1 -  t_j \right) t_j (\varphi_{g_1}^* A(g_2 \cdots g_j)^{g_1})^2 \displaybreak[0]\\
& \quad + \sum_{2 \leq i < j \leq k} - t_i t_j \left( \varphi_{g_1 \cdots g_i}^* (A(g_{i+1} \cdots g_j)^{g_1 \cdots g_i})^2 - (\varphi_{g_1}^* A(g_2 \cdots g_i)^{g_1})^2 - (\varphi_{g_1}^* A(g_2 \cdots g_j)^{g_1})^2 \right.  \displaybreak[0]\\
& \qquad \qquad +  \left. [\varphi_{g_1}^* A(g_2 \cdots g_i)^{g_1} , \varphi_{g_1}^* A(g_2 \cdots g_j)^{g_1} ]  - [\varphi_{g_1}^* A(g_2 \cdots g_i)^{g_1} , \varphi_{g_1}^* A(g_2 \cdots g_j)^{g_1} ] \right) \displaybreak[0]\\
\intertext{we can rewrite the middle three terms in the second summation as a perfect square}
&=   \sum_{i=2}^k - \left( 1 -  t_i \right) t_i (\varphi_{g_1}^* A(g_2 \cdots g_i)^{g_1})^2 \displaybreak[0] \\
& \quad + \sum_{2 \leq i < j \leq k} - t_i t_j \left( \varphi_{g_1 \cdots g_i}^* (A(g_{i+1} \cdots g_j)^{g_1 \cdots g_i})^2  - (\varphi_{g_1}^* A(g_2 \cdots g_i)^{g_1} - \varphi_{g_1}^* A(g_2 \cdots g_j)^{g_1})^2 \right. \displaybreak[0] \\
& \qquad \qquad +   \left. [\varphi_{g_1}^* A(g_2 \cdots g_i)^{g_1} , \varphi_{g_1}^* A(g_2 \cdots g_j)^{g_1} ] \right)  \displaybreak[0] \\
\intertext{and using Lemma \ref{lem:discrepancy_identity} we have}
&=  \sum_{i=2}^k  - \left( 1 -  t_i \right) t_i (\varphi_{g_1}^* A(g_2 \cdots g_i)^{g_1})^2 \displaybreak[0]\\
& \quad + \sum_{2 \leq i < j \leq k} - t_i t_j \left( \varphi_{g_1 \cdots g_i}^* (A(g_{i+1} \cdots g_j)^{g_1 \cdots g_i})^2  - (\alpha (g_2 \cdots g_i , g_{i+1} \cdots g_j) -\varphi_{g_1 \cdots g_i}^* A(g_{i+1} \cdots g_j)^{g_1 \cdots g_i})^2 \right. \\
& \qquad \qquad \left.   +  [\varphi_{g_1}^* A(g_2 \cdots g_i)^{g_1} , \varphi_{g_1}^* A(g_2 \cdots g_j)^{g_1} ] \right).  \displaybreak[2]\\
\intertext{Since $\alpha$ is a scalar-valued 1-form, the first two terms in the second summation cancel so that}
&=  \sum_{i=2}^k - \left( 1 -  t_i \right) t_i (\varphi_{g_1}^* A(g_2 \cdots g_i)^{g_1})^2 + \sum_{2 \leq i < j \leq k} - t_i t_j  [\varphi_{g_1}^* A(g_2 \cdots g_i)^{g_1} , \varphi_{g_1}^* A(g_2 \cdots g_j)^{g_1} ] 
\end{align*}
\end{proof}
\end{lem}

Now we will explicitly show that $( \id \times \partial_0 )^* \vartheta_{(k)} = (\delta_0 \times \id)^* \vartheta_{(k-1)}$, i.e.
\begin{equation}
\vartheta_{(k)} (g_1 ,\ldots , g_k) |_{t_1+ \cdots + t_k = 1} = \varphi_{g_1}^* \vartheta_{(k-1)} (g_2 , \ldots , g_k)^{g_1}.
\end{equation}
Note that when $t_1 + \cdots + t_k = 1$ we have $dt_1 + \cdots + dt_k = 0$.  We set $t_1 + \cdots + t_k =1$ in the expression we obtained for $\vartheta_{(k)} (g_1 , \ldots , g_k)$ in \eqref{eq:theta_rewritten} to obtain
\begin{multline}
 \vartheta(g_1 , \ldots , g_k)|_{\sum_{i=1}^k t_i = 1} =  \sum_{i=1}^k  t_i \varphi_{g_1\cdots g_i}^* \theta^{\E^{g_1 \cdots g_i}}  + dt_i A(g_1 \cdots g_i) - \left( \sum_{ j=1,  j\neq i}^k t_j \right)t_i A(g_1 \cdots g_i)^2 \\
 \quad + \sum_{1 \leq i < j \leq k} t_i t_j [ A(g_1 \cdots g_i) , A(g_1 \cdots g_j)] -  (t_i dt_j - t_j dt_i)  \alpha(g_1 \cdots g_i , g_{i+1} \cdots g_j).
\end{multline}
By rewriting
$$\sum_{i=1}^k \left( \sum_{ j=1,  j\neq i}^k t_j \right)t_i A(g_1 \cdots g_i)^2 = \sum_{1 \leq i < j \leq k} t_i t_j \left( A(g_1 \cdots g_i)^2 +A(g_1 \cdots g_j)^2\right), $$
we can factor the resulting terms appearing in the summation $\displaystyle \sum_{1 \leq i \leq j \leq k}$ as a difference of squares
$$
-A(g_1 \cdots g_i)^2 - A(g_1 \cdots g_j)^2 + [ A(g_1 \cdots g_i) , A(g_1 \cdots g_j)] = - \left( A(g_1 \cdots g_i) - A(g_1 \cdots g_j) \right)^2.
$$
Applying Lemma \ref{lem:discrepancy_identity} this term equals
$$
-\left(\alpha (g_1 \cdots g_i , g_{i+1} \cdots g_j) - \varphi_{g_1 \cdots g_i}^* A(g_{i+1} \cdots g_j)^{g_1 \cdots g_i} \right)^2
$$
and because $\alpha$ is a scalar-valued 1-form this term simplifies to
$$
 \left( \varphi_{g_1 \cdots g_i}^* A(g_{i+1} \cdots g_j)^{g_1 \cdots g_i} \right)^2.
$$
Hence
\begin{multline}
\label{2-form_restricted_intermediate1}
 \vartheta(g_1 , \ldots , g_k)|_{\sum_{i=1}^k t_i = 1}  =  \sum_{i=1}^k  t_i \varphi_{g_1\cdots g_i}^* \theta^{\E^{g_1 \cdots g_i}}  + dt_i A(g_1 \cdots g_i) \\
 + \sum_{1 \leq i < j \leq k} - t_i t_j \left( \varphi_{g_1 \cdots g_i}^* A(g_{i+1} \cdots g_j)^{g_1 \cdots g_i} \right)^2  -  (t_i dt_j - t_j dt_i)  \alpha(g_1 \cdots g_i , g_{i+1} \cdots g_j) .
\end{multline}
Using Lemma \ref{lem:2-form_lemma_1} to rewrite the second term in the second summation as well as rewriting the second term in the first summation of \eqref{2-form_restricted_intermediate1} using $dt_1 + \cdots + dt_k = 0$,
$$
dt_i A(g_1 \cdots g_i) = \sum_{i=2}^k dt_i (A(g_1 \cdots g_i) - A(g_1)),
$$
we obtain, 
\begin{multline*}
\vartheta(g_1 , \ldots , g_k)|_{\sum_{i=1}^k t_i = 1} = \sum_{i=1}^k  t_i \varphi_{g_1\cdots g_i}^* \theta^{\E^{g_1 \cdots g_i}} + \sum_{i=2}^k dt_i (A(g_1 \cdots g_i) - A(g_1) + \alpha (g_1 , g_2 \cdots g_i)) \\ 
 + \sum_{1 \leq i < j \leq k} - t_i t_j \left( \varphi_{g_1 \cdots g_i}^* A(g_{i+1} \cdots g_j)^{g_1 \cdots g_i} \right)^2  \\
 -  \sum_{2 \leq i < j \leq k} (t_i dt_j - t_j dt_i)  (\alpha(g_1 , g_2 \cdots g_i) - \alpha( g_1 , g_2 \cdots g_j) + \alpha(g_1 \cdots g_i , g_{i+1} \cdots g_j)) 
\end{multline*}
and by Lemma \ref{lem:discrepancy_identity} and Corollary \ref{cor:discrepancy_identity2} we can simplify to
\begin{multline*}
\label{2-form_restricted_intermediate2}
\vartheta(g_1 , \ldots , g_k)|_{\sum_{i=1}^k t_i = 1} =  \sum_{i=1}^k  t_i \varphi_{g_1\cdots g_i}^* \theta^{\E^{g_1 \cdots g_i}}  - \sum_{i=2}^k dt_i \varphi_{g_1}^* A(g_2 \cdots g_i)^{g_1} \\+ \sum_{1 \leq i < j \leq k} - t_i t_j \left( \varphi_{g_1 \cdots g_i}^* A(g_{i+1} \cdots g_j)^{g_1 \cdots g_i} \right)^2
-  \sum_{2 \leq i < j \leq k} (t_i dt_j - t_j dt_i)  \varphi_{g_1}^* \alpha(g_2 \cdots g_i , g_{i+1} \cdots g_j)^{g_1} .
\end{multline*}

Now focusing on the third summation of the previous expression, we apply Lemma \ref{lem:2-form_lemma_2} to obtain
\begin{multline*}
\vartheta(g_1 , \ldots , g_k)|_{\sum_{i=1}^k t_i = 1} =  \sum_{i=1}^k  t_i \varphi_{g_1\cdots g_i}^* \theta^{\E^{g_1 \cdots g_i}}  \\
 + \sum_{i=2}^k -  dt_i \varphi_{g_1}^* A(g_2 \cdots g_i)^{g_1}  - \left( 1 -  t_i \right) t_i (\varphi_{g_1}^* A(g_2 \cdots g_i)^{g_1})^2 \\
 - \sum_{2 \leq i < j \leq k}  t_i t_j  [\varphi_{g_1}^* A(g_2 \cdots g_i)^{g_1} , \varphi_{g_1}^* A(g_2 \cdots g_j)^{g_1} ]  + (t_i dt_j - t_j dt_i)  \varphi_{g_1}^* \alpha(g_2 \cdots g_i , g_{i+1} \cdots g_j)^{g_1} 
\end{multline*}
Finally, as $t_1 = 1-t_2 - \cdots - t_k$, we rewrite first term of the first sum to obtain
\begin{multline*}
\vartheta(g_1 , \ldots , g_k)|_{\sum_{i=1}^k t_i = 1} =  (1-t_2 - \cdots - t_k) \varphi_{g_1}^* \theta^{\E^{g_1}} + \sum_{i=2}^k  t_i \varphi_{g_1\cdots g_i}^* \theta^{\E^{g_1 \cdots g_i}}  \\
 + \sum_{i=2}^k -  dt_i \varphi_{g_1}^* A(g_2 \cdots g_i)^{g_1}  - \left( 1 -  t_i \right) t_i (\varphi_{g_1}^* A(g_2 \cdots g_i)^{g_1})^2 \\
 - \sum_{2 \leq i < j \leq k}  t_i t_j  [\varphi_{g_1}^* A(g_2 \cdots g_i)^{g_1} , \varphi_{g_1}^* A(g_2 \cdots g_j)^{g_1} ]  + (t_i dt_j - t_j dt_i)  \varphi_{g_1}^* \alpha(g_2 \cdots g_i , g_{i+1} \cdots g_j)^{g_1} 
\end{multline*}
 which is the form of $\varphi_{g_1}^* \vartheta_{(k-1)} (g_2 , \ldots , g_k)^{g_1}$ that we obtained in \eqref{eq:theta_rewritten} (with a relabeling $t_i \rightarrow t_{i-1}$ for $2 \leq i \leq k$), thus completing the calculation.

\end{proof}
\end{thm}

\begin{rmk}
If we define $\vartheta_{(k)}'$ by $(\nabla^k)^2 = [\vartheta_{(k)} ' , \cdot]$ then the discrepancy between $\vartheta_{(k)}'$ and $\vartheta_{(k)}$, $(\vartheta_{(k)}' - \vartheta_{(k)}) (g_1 , \ldots , g_k) $ is given by
\begin{equation}
\sum_{i=1}^k t_i \theta_{g_1 \cdots g_i}  + \sum_{1 \leq i < j \leq k} \alpha(g_1 \cdots g_i , g_{i+1} \cdots g_j) (t_i dt_j - t_j dt_i)
\end{equation}
which is a scalar-valued form.  In other words, $\vartheta_{(k)} ' - \vartheta_{(k)} \in \Omega^2 (M_{(k)}^\Delta )$.  Since a scalar valued form in $\Omega^2 (M_{(k)}^\Delta)$ is central in $\Omega^2 (M_{(k)}^\Delta , \End \E_{k}^\Delta)$, $[\vartheta_{(k)}' , \cdot] = [\vartheta_{(k)}, \cdot]$.  Hence for any $a \in \Gamma^\infty (\End \E_{k}^\Delta)$,
\begin{equation}
\label{eq:nabla_squared_equals_theta}
(\nabla^k)^2 (a) = [\vartheta_{(k)} , a].
\end{equation}

\end{rmk}

\begin{defn}
Notice from the previous formula for $\vartheta$ that $\vartheta^{0,2} = 0$. We will then define $\vartheta_u = u \vartheta^{2,0} + \vartheta^{1,1}$.  With this rescaling we can view $\vartheta_u$ as an element of $\Omega^0_- (M_\bullet)[u]$.  
\end{defn}


\subsection{Simplicial curvature 3-form}

\begin{thm}
Let $\vartheta$ be the simplicial $2$-form associated to $(L, \mu, \nabla)$ as in Theorem \ref{thm:simplicial_2_form}.  Let $\Theta_{(k)} = - \nabla^k \vartheta_{(k)}$. The collection of $3$-forms $\Theta_{(k)}$ define a simplicial $3$-form $\Theta = \{ \Theta_{(k)} \} \in \Omega^3 (M_\bullet)$.  Furthermore, $\Theta_{(k)}^{0,3} = 0$ as an element of $\Omega^{r,s} (M_\bullet)$. The simplicial $3$-form $\Theta = \{- \nabla^k \vartheta_{(k)} \}$ is called the \emph{simplicial curvature 3-form} associated to $(L , \mu, \nabla)$.
\begin{proof}
First let us provide a formula for $\Theta_{(k)}$. Using the Bianchi identity
\begin{align*}
& \Theta_{(k)} (g_1, \ldots , g_k) = \\
&= -(\nabla^k \vartheta_{(k)} ) (g_1 , \ldots , g_k) \displaybreak[2]\\
&= -(\nabla^k)(g_1 , \ldots , g_k) \left( (\nabla^k)^2 (g_1 , \ldots , g_k) - \sum_{i=1}^k t_i \theta_{g_1 \cdots g_i}  - \sum_{1 \leq i < j \leq k} \alpha ( g_1 \cdots g_i , g_{i+1} \cdots g_j ) (t_i dt_j - t_j dt_i ) \right) \displaybreak[2]\\
&=  -(\nabla^k)^3 (g_1 , \ldots , g_k) +  (\nabla^k (g_1 , \ldots , g_k) ) \left(  \sum_{i=1}^k t_i \theta_{g_1 \cdots g_i} + \sum_{1 \leq i < j \leq k} \alpha ( g_1 \cdots g_i , g_{i+1} \cdots g_j ) (t_i dt_j - t_j dt_i ) \right) \displaybreak[0]\\
&= \sum_{i=1}^k dt_i \theta_{g_1 \cdots g_i}  + \sum_{1 \leq i < j \leq k} d\alpha (g_1 \cdots g_i, g_{i+1} \cdots g_j) (t_i dt_j - t_j dt_i)  + 2 \alpha (g_1 \cdots g_i , g_{i+1} \cdots g_j) dt_i dt_j,
\end{align*}
since $\theta_g$ and $\alpha(g , g')$ are scalar-valued forms.  Since there are no terms containing $dt_{i_1} dt_{i_2} dt_{i_3}$ where $1 \leq i_1 , i_2 , i_3 \leq k$, we explicitly see that $\Theta_{(k)}^{0,3} = 0$.

Using this formula, we can see $\Theta$ is a compatible form by a direct calculation. First we check that $(\id \times \partial_\ell)^* \Theta_{(k)}  = (\delta_\ell \times \id)^* \Theta_{(k-1)} $ for $1 \leq \ell \leq k$.  Explicitly, when $t_\ell = 0$ we have
\begin{multline*}
\Theta_{(k)} (g_1 , \ldots  , g_k) |_{t_\ell = 0} = \sum_{\substack{ i=1 \\ i \neq \ell } }^k dt_i \theta_{g_1 \cdots g_i} \displaybreak[0]\\
+ \sum_{\substack{1 \leq i < j \leq k \\ i,j \neq \ell} } d\alpha (g_1 \cdots g_i, g_{i+1} \cdots g_j) (t_i dt_j - t_j dt_i) + 2 \alpha (g_1 \cdots g_i , g_{i+1} \cdots g_j) dt_i dt_j \displaybreak[1]\\
= \Theta_{(k-1)} (g_1 , \ldots , g_\ell g_{\ell + 1} , \ldots , g_k)
\end{multline*}
by substituting $t_i$ for $t_{i+1}$ whenever $i > \ell$.  Now when $\ell=0$, on the face $t_1 + \cdots + t_k = 1$ we have $dt_1 + \cdots + dt_k = 0$ so that by replacing $t_1$
\begin{multline*}
\Theta_{(k)} (g_1 , \ldots  , g_k) |_{t_1 + \cdots + t_k = 1} = \sum_{i=2}^k dt_i \theta_{g_1 \cdots g_i} - \sum_{i=2}^k dt_i \theta_{g_1} \displaybreak[0]\\
+ \sum_{2 \leq i < j \leq k} d\alpha (g_1 \cdots g_i, g_{i+1} \cdots g_j) (t_i dt_j - t_j dt_i) \displaybreak[0]\\
+ \sum_{i=2}^k d \alpha (g_1 , g_2 \cdots g_i) ( (1- t_2 - \cdots - t_k) dt_i - t_i (- dt_2 - \cdots - dt_k))\displaybreak[0] \\
 + \sum_{2 \leq i < j \leq k} 2 \alpha (g_1 \cdots g_i , g_{i+1} \cdots g_j) dt_i dt_j \displaybreak[0]
 + \sum_{i=2}^k 2 \alpha (g_1 , g_2 \cdots g_i) (-dt_2 -\cdots -dt_k) dt_i \displaybreak[2]
\end{multline*}
\begin{align*}
\intertext{By Proposition \ref{prop:gerbe_connection2_group} applied in the first two summations, Corollary \ref{cor:discrepancy_identity2}, and $dt_i dt_j = - dt_j dt_i$,}
&= \sum_{i=2}^k dt_i ( d \alpha (g_1 , g_2 \cdots g_i) + \theta_{g_2 \cdots g_i}^{g_1}) \displaybreak[0]\\
&\quad + \sum_{2 \leq i < j \leq k} (d \alpha (g_1 , g_2 \cdots g_i) - d\alpha (g_1 , g_2 \cdots g_j) + d \alpha (g_1 \cdots g_i , g_{i+1} \cdots g_j) ) (t_i dt_j - t_j dt_i) \displaybreak[0]\\
&\quad + \sum_{i=2}^k d\alpha (g_1 , g_2 \cdots g_i) dt_i + \sum_{2 \leq i < j \leq k} 2 (\alpha (g_1 , g_2 \cdots g_i ) - \alpha (g_1 , g_2 \cdots g_j) + \alpha ( g_1 \cdots g_i , g_{i+1} \cdots g_j)) dt_i dt_j \\
&=  \sum_{i=2}^k dt_i \theta_{g_2 \cdots g_i}^{g_1} + \sum_{2 \leq i < j \leq k} d\alpha (g_2 \cdots g_i , g_{i+1} \cdots g_j)^{g_1} (t_i dt_j -t_j dt_i)  +  2\alpha (g_2 \cdots g_i , g_{i+1} \cdots g_j)^{g_1} dt_i dt_j \displaybreak[1]\\
&= \Theta_{(k-1)} (g_2 ,\ldots , g_k)^{g_1}
\end{align*}
by a cancellation and application of Corollary \ref{cor:discrepancy_identity2} again.
\end{proof}
\end{thm}

\begin{thm}
Let $\alpha$ be the discrepancy of $( L , \mu , \nabla)$, $\vartheta$ the associated simplicial 2-form, and $\Theta$ the associated simplicial curvature 3-form. Then $\Theta \in \Omega^3 (M_\bullet)$ is a representative of the Dixmier-Douady class of $(L, \mu, \nabla)$.
\begin{proof}

We will apply the quasi-isomorphism of integration along simplices $\mathcal{I}_\Delta$ described in Theorem \ref{thm:simplex_integration} to see $\mathcal{I}_\Delta (\Theta) = (\alpha , \theta) \in \Omega^1 ( M_{(2)}) \oplus \Omega^2 ( M_{(1)} )$.  For any $g \in \Gamma$,
\begin{equation*}
\Theta_{(1)} (g) = dt_1 \theta_g
\end{equation*}
so that restricted to $M \times \{g\}$
\begin{equation*}
\mathcal{I}_\Delta : \Theta_{(1)} (g) \mapsto \int_{\Delta^1} \Theta_{(1)} (g) = \int_0^1  \theta_g dt_1 = \theta_g.
\end{equation*}
For any $g_1 , g_2 \in \Gamma$,
\begin{equation*}
\Theta_{(2)} (g_1 , g_2) = dt_1 \theta_{g_1}  dt_2 \theta_{g_2} + d\alpha (g_1 , g_2) (t_1 dt_2 - t_2 dt_1) + 2 \alpha (g_1 , g_2) dt_1 dt_2
\end{equation*}
so that restricted to $M \times \{g_1\} \times \{g_2\}$ we have
\begin{align*}
\mathcal{I}_\Delta : \Theta_{(2)} (g_1 , g_2) &\mapsto \int_{\Delta^2} \Theta_{(2)} (g_1 , g_2)\\
&= \int_{\Delta^2} dt_1 \theta_{g_1} dt_2 \theta_{g_2} + d\alpha (g_1 , g_2) (t_1 dt_2 - t_2 dt_1) + 2 \alpha (g_1 , g_2) dt_1 dt_2 \\
&= \int_{\Delta^2} 2 \alpha (g_1 , g_2) dt_1 dt_2 \\
&= \int_0^1 \int_0^{1-t_1} 2 \alpha (g_1 , g_2) dt_1 dt_2 \\
&= \alpha (g_1 , g_2).
\end{align*}
Furthermore, $\Theta_{(k)}^{r,s} = 0$ for $s \geq 3$, i.e.\ $\Theta_{(k)}$ has at most two $dt$'s in any term, so that 
\begin{equation*}
\mathcal{I}_\Delta (\Theta_{(k)}) = \int_{\Delta^k} \Theta_{(k)}= 0
\end{equation*}
whenever $k \geq 3$. Hence $\mathcal{I}_\Delta : \Theta \mapsto \alpha + \theta$ so that $\Theta$ is a representative of the Dixmier-Douady class defined by $(\alpha, \theta)$ in the bicomplex of compatible forms.
\end{proof}
\end{thm}

We may also introduce a rescaled version of the simplicial Dixmier-Douady form $\Theta_u$ in the complex $\underline{\Omega} (M_\bullet) [u]^*$ as in \eqref{eq:Theta_rescaled}.  With the rescaling $\nabla^k_u =  (\nabla^k)^{1,0} + u^{-1}d_\Delta$ we have $\nabla^k_u (\vartheta_u)_{(k)} = u^{-1} (\Theta_u)_{(k)}$ as
\begin{align*}
( (\nabla^k)^{1,0} + u^{-1} d_\Delta) ( u (\vartheta_{(k)})^{2,0} + (\vartheta_{(k)})^{1,1}) &= u ((\nabla^k)^{1,0} (\vartheta_{(k)})^{2,0})^{3,0} \\
&\quad +  (( \nabla^k)^{1,0} (\vartheta_{(k)})^{1,1} + d_\Delta (\vartheta_{(k)})^{2,0} )^{2,1} \\
&\quad + u^{-1} (d_\Delta (\vartheta_{(k)}))^{1,2} \\
&=u (\Theta_{(k)})^{3,0} +  (\Theta_{(k)})^{2,1} + u^{-1} (\Theta_{(k)})^{1,2}
\end{align*}
Also note that under this rescaling, $u (\nabla^k_u)^2 = [ (\vartheta_u)_{(k)} , \cdot]$.


\section{Cocycles on the twisted convolution algebra}
\label{chap:4}

Now we will prove the main theorem of this paper:

\begin{thm}
\label{thm:main_thm}
Given a gerbe datum $(L , \mu, \nabla)$ on $M \rtimes \Gamma$, the map
\begin{equation}
\Psi_2 \circ \Psi_1 \circ \Psi_0 \circ \tau_\nabla : (\underline{\Omega} ((M\rtimes \Gamma) _\bullet) [u]^\bullet , d_{\Theta_u})  \to CC^\bullet (C^\infty_c (M \rtimes \Gamma , L)  [u^{-1} , u], b+uB)
\end{equation}
is a morphism from the twisted simplicial complex of $M \rtimes \Gamma$ to the $(b,B)$ complex of the twisted convolution algebra $C^\infty_c ( M \rtimes \Gamma , L)$.
\end{thm}

In Section \ref{sec:twisted_convolution_algebra} we define $C^\infty_c ( M \rtimes \Gamma , L)$.  The JLO morphism, 
\begin{equation*}
\tau_\nabla : (\underline{\Omega} (M_\bullet) [u]^\bullet , d_{\Theta_u})  \to (C^\bullet (\Gamma , CC^\bullet ( \Gamma^\infty_c ( \End \E) )  [u^{-1} , u] ) , b+uB + {\delta_\Gamma}'),
\end{equation*}
is described in Section \ref{sec:JLO_morphism}.  The morphisms $\Psi_0$, $\Psi_1$, and $\Psi_2$ are algebraic morphisms which compose to give a morphism
\begin{equation*}
(C^\bullet (\Gamma , CC^\bullet ( \Gamma^\infty_c ( \End \E) )  [u^{-1} , u] ) , b+uB + {\delta_\Gamma}') \to CC^\bullet (C^\infty_c (M \rtimes \Gamma , L)  [u^{-1} , u], b+uB)
\end{equation*}
and are described in Section \ref{sec:algebraic_morphisms}.

\subsection{Twisted Convolution Algebra}
\label{sec:twisted_convolution_algebra}

\begin{defn}
For a gerbe datum $(L, \mu, \nabla)$ on $\G$, the \textit{twisted convolution algebra} $C^\infty_c (\G, L)$ is defined as follows. The convolution product on $C^\infty_c (\G , L)$ is defined by the (finite) sum
\begin{equation}
 (f_1 * f_2 ) (g) = \sum_{\{(g_1, g_2) \in \G^{(2)} | g_1 g_2 = g \} } f_1 (g_1) \cdot f_2 (g_2)
 \end{equation}
for sections $f_1 , f_2 \in C^\infty_c (\G , L)$, where $f_1 (g_1) \cdot f_2 (g_2)$ is defined by applying the isomorphism $\mu_{(g_1,g_2)} : L|_{g_1} \otimes L|_{g_2} \xrightarrow{\sim} L|_{g}$.  This gives $C^\infty_c (\G, L)$ the structure of an algebra.

In particular, for $\G = M \rtimes \Gamma$, the convolution product on $C^\infty_c (M \rtimes \Gamma, L)$ may be written
\begin{align}
(f_1 * f_2) (x,g) &= \sum_{g_1 g_2 = g}  \mu_{g_1 , g_2} (f_1|_{M_{g_1}} \otimes (f_2|_{M_{g_2}})^{g_1} ) (x , g) \\
&= \sum_{g_1 g_2 = g} f_1 (x,g_1) \cdot f_2 ( xg_1 , g_2)
\end{align}
for sections $f_1 , f_2 \in C^\infty_c (M \times \Gamma , L)$, where $f_1 (x,g_1) \cdot f_2 (xg_1 , g_2)$ is computed using the product $\mu_{g_1 , g_2} : L_{g_1} \otimes (L_{g_2})^{g_1} \xrightarrow{\sim} L_g$.
\end{defn}

\subsection{JLO morphism}
\label{sec:JLO_morphism}

Following \cite{jlo88}, \cite{gorokhovsky99}, \cite{mathai_stevenson06}, and \cite{tu_xu06} we now construct a JLO-type morphism from the twisted simplicial complex of $M \rtimes \Gamma$.

For an algebra $\A$, let $\widetilde{\A} = \A \oplus \C$ denote the unitization of $\A$.  Define 
\begin{equation}
\label{eq:cochains}
CC^n (\A) = \Hom (\widetilde{\A} \otimes \A^{\otimes n} , \C).
\end{equation}
In the following, we will be considering compactly supported sections which are topologized as the inductive limit of Fr\'echet spaces.  Therefore we will use the inductive tensor product.
The periodic cyclic cohomology of $\A$ is the cohomology of the complex 
\begin{equation}
(CC^n (\A) [u^{-1},u] , b+uB),
\end{equation}
where $b$ is a differential of degree $+1$ and $B$ is a differential of degree $-1$ defined as follows.  Let
\begin{align}
d^i f (\tilde{a}_0 , a_1 , \ldots , a_n) &= \begin{cases} f(\tilde{a}_0 , a_1 , \ldots , a_i a_{i+1} , \ldots , a_n) & \textrm{if } 0 \leq i \leq n-1 \\ f(a_n \tilde{a}_0 , a_1 , \ldots , a_{n-1} ) & \textrm{if } i=n \end{cases} \\
\intertext{for $f \in CC^{n-1} (\A)$ and}
\bar{s}^i_n f (\tilde{a}_0 , \ldots , a_n) &= (-1)^{ni} f(1 , a_i , \ldots, a_n , a_0 , \ldots , a_{i-1}), \ 0 \leq i \leq n,
\end{align}
for $f \in CC^{n+1} (\A)$. Then
\begin{equation}
b^n = \sum_{i=0}^n d^i_n, \label{eq:Hochschild_coboundary} \qquad 
B^n = \sum_{i=0}^n \bar{s}^i_n, 
\end{equation}
and, dropping the superscript $n$ denoting degree, $b^2 = B^2 = bB+Bb=0$ so that $b+uB$ is a differential of degree $+1$ on $CC^n (\A)[u^{-1} , u]$.

For the discrete group $\Gamma$ and a $\Gamma$-module $K$, let $C^n (\Gamma , K)$ denote the space of degree $n$ $\Gamma$-cochains, i.e.\ the space of maps $\Gamma^n \rightarrow K$.  For any $f \in C^n (\Gamma, K)$ define $\delta_\Gamma : C^n (\Gamma , K) \rightarrow C^{n+1} (\Gamma , K)$ to be the group coboundary $\delta_\Gamma = \sum_{i=0}^n (-1)^i {\delta_\Gamma}_i$ with
\begin{equation*}
({\delta_\Gamma}_i f) (g_1 , \ldots , g_{n+1}) = 
\begin{cases}
f(g_2 , \ldots , g_{n+1})^{g_1} & \text{for $i=0$}\\
 (-1)^i f( g_1 , \ldots , g_{i-1} , g_i g_{i+1} , g_{i+2} , \ldots , g_{n+1}) & \text{for $1 \leq i \leq n$}\\
 (-1)^{n+1} f(g_1 , \ldots , g_n) & \text{for $i=n+1$}.
\end{cases}
\end{equation*}
The superscript ${}^{g_1}$ denotes the action of $g_1 \in \Gamma$ on $K$.  The cohomology of the complex $(C^\bullet (\Gamma, K) , {\delta_\Gamma})$ is group cohomology of $\Gamma$ with coefficients in $K$.  In particular we will consider the complex 
\begin{equation}
(C^k (\Gamma, CC^n( \Gamma^\infty_c (\End \E) ) [u^{-1} , u] ) , (b+uB) + (-1)^n {\delta_\Gamma})
\end{equation}
of $\Gamma$-cochains with values in the periodic cyclic complex of $\Gamma^\infty_c ( \End \E)$.  By the operators $b$ and $uB$ on such $\Gamma$-cochains we mean the operators defined by 
\begin{equation}
( b f) (g_1 , \ldots , g_n) = b (f(g_1 , \ldots , g_n))
\end{equation}
 and 
\begin{equation}
(uB f) (g_1 , \ldots , g_n) = uB (f (g_1 , \ldots , g_n))
\end{equation}
for any $f \in C^n (\Gamma, CC^\bullet ( \Gamma^\infty_c ( \End \E) ) [u^{-1} , u] )$. Let us also denote ${\delta_\Gamma}' = (-1)^n {\delta_\Gamma}$ when $\delta_\Gamma$ is of degree $n$.

\begin{lem}
\label{lem:Stokes_twisted_simplicial_complex}
For $\alpha_{(n)} \otimes \beta_{(n)} \in \underline{\Omega} ( (M \rtimes \Gamma)_{(n)} ) [u]^0 \otimes \Omega^{n-1} (\Delta^n)$ we have the following
\begin{equation*}
\int_M \int_{\Delta^n} (u\tilde{d}_{dR} + d_\Delta) (\alpha_{(n)} \otimes \beta_{(n)}) = \int_M \alpha_{(n)} \int_{\partial (\Delta^n)} \beta_{(n)}.
\end{equation*}
Also, for $\beta_{(n)}$ not of degree $n-1$ or $\alpha_{(n)}$ not of degree $0$ as elements of $( \underline{\Omega} (M_\bullet)[u]^* ) |_{ (M \rtimes \Gamma)_{(n)} \times \Delta^n  }$ then the integral over $M \times \Delta^n$ is zero.
\begin{proof}
Note that $\alpha_{(n)} \in \underline{\Omega} ( (M \rtimes \Gamma)_{(n)} ) [u]^0$ must be a $(\dim M)$-form on $(M \rtimes \Gamma)_{(n)}$.   For any $g_1 , \ldots , g_n \in \Gamma$ denote $\alpha = \alpha_{(n)} (g_1 , \ldots , g_n)$. Then 
\begin{align*}
\int_M \int_{\Delta^n} (u\tilde{d}_{dR} + d_\Delta) (\alpha \otimes \beta_{(n)}) &= \int_M \int_{\Delta^n} (u \tilde{d}_{dR} \alpha) \otimes \beta_{(n)} + \int_M \int_{\Delta^n} \alpha \otimes (d_\Delta \beta_{(n)}) \\
&=  \int_M \alpha  \int_{\Delta^n}  d \beta_{(n)} \\
&=  \int_M \alpha  \int_{\partial (\Delta^n)} \beta_{(n)}
\end{align*}
by Stokes' theorem.
\end{proof}
\end{lem}

\begin{lem}
\label{lem:Stokes_wedge_twisted_simplicial_complex}
For any $\omega_1 \in \underline{\Omega} ( (M\rtimes \Gamma)_\bullet) [u]^*$ and $ \omega_2 \in \Omega^* ( (M\rtimes \Gamma)_\bullet)[u]$ 
\begin{equation*}
\int_M \int_{\Delta^n} ((u \tilde{d}_{dR} + d_\Delta) \omega_1) \wedge \omega_2  = - \int_M \int_{\Delta^n} \omega_1 \wedge ( (ud_{dR} + d_\Delta) \omega_2) + (-1)^{|\omega_2|} \int_M \int_{\partial (\Delta^n)} \omega_1 \wedge \omega_2
\end{equation*}
where $\omega_2$ is of fixed degree.
\begin{proof}
We have
\begin{align*}
&\int_M \int_{\Delta^n} ((u \tilde{d}_{dR} + d_\Delta) \omega_1) \wedge \omega_2 \\
&= -  \int_M \int_{\Delta^n} \omega_1 \wedge ( (u d_{dR} + d_\Delta) \omega_2) + (-1)^{|\omega_2|}  \int_M \int_{\Delta^n} (u \tilde{d}_{dR} + d_\Delta) (\omega_1 \wedge \omega_2) \\
&= -  \int_M \int_{\Delta^n} \omega_1 \wedge ( (u d_{dR} + d_\Delta) \omega_2) + (-1)^{|\omega_2|} \int_M \int_{\partial(\Delta^n)} \omega_1 \wedge \omega_2
\end{align*}
where the first inequality follows from Lemma \ref{lem:wedge_product_twisted_simplicial_complex} and the second from Lemma \ref{lem:Stokes_twisted_simplicial_complex}.
\end{proof}
\end{lem}

\begin{thm}
Given a gerbe datum $(L , \mu, \nabla)$ on $M \rtimes \Gamma$ with associated simplicial $2$-form $\vartheta$, simplicial curvature 3-form $\Theta$, and direct sum bundle $\E$, there is a morphism 
\begin{equation}
\tau_\nabla : (\underline{\Omega} (( M \rtimes \Gamma)_\bullet) [u]^\bullet , d_{\Theta_u} ) \to (C^\bullet (\Gamma , CC^\bullet ( \Gamma^\infty_c ( \End \E) )  [u^{-1} , u] ) , b+uB + {\delta_\Gamma}')
\end{equation}
defined by the JLO-type formula
\begin{multline}
\tau_\nabla (  \omega ) (\tilde{a}_0 , a_1 , \ldots , a_n) \\ = \sum_k \int_M \int_{\Delta^k}  \omega_{(k)} \wedge \left( \int_{\Delta^n} \tr (\tilde{a}_0 e^{-\sigma_0 \left( \vartheta_u \right)_{(k)}} \nabla_u^k (a_1) e^{-\sigma_1  \left( \vartheta_u \right)_{(k)}} \cdots \nabla_u^k (a_n) e^{-\sigma_n  \left( \vartheta_u \right)_{(k)}} ) d \sigma_1 \cdots d \sigma_n \right) 
\end{multline}
where $\omega =  \{ \omega_{(k)} \}  \in \underline{\Omega} ( (M \rtimes \Gamma)_\bullet))[u]^*$, $a_1 , \ldots , a_n \in \Gamma^\infty_c (\End \E)$, $\tilde{a}_0 \in \widetilde{\Gamma^\infty_c ( \End \E)}$, and $\sigma_0 , \ldots , \sigma_n$ are barycentric coordinates on $\Delta^n$.
\begin{proof}
Let us first discuss the reason the sum is finite.  We must check in particular that the degrees of the forms to be integrated on the simplices $\Delta^k$ are bounded as $k$ increases.  Given $\omega \in \underline{\Omega} ( (M \rtimes \Gamma)_\bullet) [u]^\ell$, and $a_0 , \ldots , a_n \in \Gamma^\infty_c (\End \E) $, consider
\begin{equation}
\omega_{(k)} \wedge \left( \int_{\Delta^n} \tr (\tilde{a}_0 e^{-\sigma_0 \left( \vartheta_u \right)_{(k)}} \nabla_u^k (a_1) e^{-\sigma_1  \left( \vartheta_u \right)_{(k)}} \cdots \nabla_u^k (a_n) e^{-\sigma_n  \left( \vartheta_u \right)_{(k)}} ) d \sigma_1 \cdots d \sigma_n \right) 
\end{equation}
The maximum number of $dt_i$'s due to $\omega$ is $\ell$, where $t_i$ are coordinates on the simplex $\Delta^k$, as the number of $dt_i$'s cannot exceed the degree of $\omega$.  Within the trace, the terms $\nabla^k_u (a_i)$ do not contribute any $dt_i$'s because $a_j$ is extended trivially to $\Delta^n$ and $(\nabla^k)^{0,1} = d_\Delta$.  Since $\vartheta_{(k)} = \vartheta_{(k)}^{2,0} + \vartheta_{(k)}^{1,1}$ and in particular $\vartheta_{(k)}^{0,2} = 0$, the contribution to the degree of the form from all the $(\vartheta_u)_{(k)}$ appearing must be at most $\dim M$, i.e.\ the number of $dx_i$'s, coordinates on the manifold $M$, is greater than or equal to the number of $dt_j$'s. Hence the simplex degree of the form above is at most $\ell  + \dim M$ so the sum is finite.

Without loss of generality we will compute $\tau_\nabla ( d_{\Theta_u} \omega) (a_0 , a_1 , \ldots , a_n)$ since $$\tilde{a}_0 e^{-\sigma_0 \left( \vartheta_u \right)_{(k)}} \nabla_u^k (a_1) \cdots  e^{-\sigma_n  \left( \vartheta_u \right)_{(k)}} $$ in the integrand above means
\begin{equation*}
a_0 e^{-\sigma_0 \left( \vartheta_u \right)_{(k)}} \nabla_u^k (a_1)  \cdots  e^{-\sigma_n  \left( \vartheta_u \right)_{(k)}} + \lambda e^{-\sigma_0 \left( \vartheta_u \right)_{(k)}} \nabla_u^k (a_1)  \cdots e^{-\sigma_n  \left( \vartheta_u \right)_{(k)}} 
\end{equation*}
for $\tilde{a}_0 = (a_0 , \lambda)$.

Following \cite{mathai_stevenson06}, we will compute the effect of $(u\tilde{d}_{dR} + d_{\Delta})$.  This gives
\begin{align}
&\tau_\nabla ( (u\tilde{d}_{dR} + d_{\Delta}) \omega) (a_0 , a_1 , \ldots , a_n)  \notag \\
&= \sum_k  \int_M \int_{\Delta^k}  (u\tilde{d}_{dR} + d_{\Delta}) \omega_{(k)} \wedge \bigg( \int_{\Delta^n} \tr (a_0 e^{-\sigma_0 \left( \vartheta_u \right)_{(k)}} \nabla_u^k (a_1) \cdots  e^{-\sigma_n  \left( \vartheta_u \right)_{(k)}} ) d \sigma_1 \cdots d \sigma_n \bigg)  \displaybreak[3] \notag\\
&=- \sum_k  \int_M \int_{\Delta^k}  \omega_{(k)} \wedge \bigg( \int_{\Delta^n} (ud_{dR} + d_{\Delta}) \tr (a_0 e^{-\sigma_0 \left( \vartheta_u \right)_{(k)}} \nabla_u^k (a_1) \cdots e^{-\sigma_n  \left( \vartheta_u \right)_{(k)}} ) d \sigma_1 \cdots d \sigma_n \bigg)  \displaybreak[2]\notag\\
&\quad + (-1)^{n} \int_M \int_{\partial(\Delta^k)}  \omega_{(k)} \wedge \bigg( \int_{\Delta^n} \tr (a_0 e^{-\sigma_0 \left( \vartheta_u \right)_{(k)}} \nabla_u^k (a_1) \cdots e^{-\sigma_n  \left( \vartheta_u \right)_{(k)}} ) d \sigma_1 \cdots d \sigma_n \bigg)  \notag\\
&= -\sum_k  \int_M \int_{\Delta^k}  \omega_{(k)} \wedge \left( \int_{\Delta^n} u \tr (\nabla^k_u (a_0 e^{-\sigma_0 \left( \vartheta_u \right)_{(k)}}   \cdots e^{-\sigma_n  \left( \vartheta_u \right)_{(k)}} ) ) d \sigma_1 \cdots d \sigma_n \right)  \displaybreak[2] \label{eq:JLO_formula_1}\\
&\quad + (-1)^{n } \int_M \int_{\partial(\Delta^k)}  \omega_{(k)} \wedge \left( \int_{\Delta^n} \tr (a_0 e^{-\sigma_0 \left( \vartheta_u \right)_{(k)}}   \cdots  e^{-\sigma_n  \left( \vartheta_u \right)_{(k)}} ) d \sigma_1 \cdots d \sigma_n \right) , \label{eq:JLO_formula_2}
\end{align}
where the second equality follows from Lemma \ref{lem:Stokes_wedge_twisted_simplicial_complex}, as the terms inside the trace $\nabla_u^k (a_i)$ are of degree 1 in $\Omega^* ( (M \rtimes \Gamma)_\bullet) [u]$ for $1 \leq i \leq n$ and the terms $e^{-\sigma_i  \left( \vartheta_u \right)_{(k)}} $ are of even degree.  The last equality is from Lemma \ref{lem:trace_connection}. By commuting forms inside the trace within the first term \eqref{eq:JLO_formula_1} we have
\begin{align*}
&u \tr \left( \nabla^k_u (a_0 e^{-\sigma_0 \left( \vartheta_u \right)_{(k)}} \nabla_u^k (a_1) \cdots e^{-\sigma_n  \left( \vartheta_u \right)_{(k)}} ) \right) \\
&=  u \tr \left( \nabla^k_u  (a_0) e^{-\sigma_0 \left( \vartheta_u \right)_{(k)}} \nabla_u^k (a_1) \cdots e^{-\sigma_n  \left( \vartheta_u \right)_{(k)}} \right) \\
&\quad + \sum_{i=0}^{n-1} (-1)^i \left( u \tr \left(a_0 e^{-\sigma_0 \left( \vartheta_u \right)_{(k)}} \cdots \nabla^k_u (a_i) \nabla^k_u \left(e^{-\sigma_i \left( \vartheta_u \right)_{(k)}} \right) \nabla^k_u (a_{i+1}) \cdots e^{-\sigma_n  \left( \vartheta_u \right)_{(k)}}  \right) \right. \\
&\qquad \qquad + \left. u \tr \left(  a_0 e^{-\sigma_0 \left( \vartheta_u \right)_{(k)}} \cdots \nabla^k_u (a_i) e^{-\sigma_i \left( \vartheta_u \right)_{(k)}} (\nabla^k_u)^2 (a_{i+1}) \cdots e^{-\sigma_n  \left( \vartheta_u \right)_{(k)}} \right) \right) \\
&\quad + (-1)^n u \tr \left( a_0 e^{-\sigma_0 \left( \vartheta_u \right)_{(k)}} \nabla_u^k (a_1) \cdots \nabla^k_u (a_n) \nabla^k_u (e^{-\sigma_n  \left( \vartheta_u \right)_{(k)}} ) \right) \displaybreak[2]\\
&=  u\tr \left( \nabla^k_u  (a_0) e^{-\sigma_0 \left( \vartheta_u \right)_{(k)}} \nabla_u^k (a_1) \cdots e^{-\sigma_n  \left( \vartheta_u \right)_{(k)}} \right) \\
&\quad + \sum_{i=0}^{n}  (-1)^{i+i(n+1-i)} u \tr \bigg( \nabla^k_u \left(e^{-\sigma_i \left( \vartheta_u \right)_{(k)}} \right) \nabla^k_u (a_{i+1}) e^{-\sigma_n  \left( \vartheta_u \right)_{(k)}} a_0 e^{-\sigma_0 \left( \vartheta_u \right)_{(k)}} \cdots \nabla^k_u (a_i)  \bigg) \\
&\quad + \sum_{i=0}^{n-1}  (-1)^i  \tr \left(  a_0 e^{-\sigma_0 \left( \vartheta_u \right)_{(k)}} \cdots \nabla^k_u (a_i) e^{-\sigma_i \left( \vartheta_u \right)_{(k)}} [(\vartheta_u)_{(k)} , a_{i+1}] \cdots e^{-\sigma_n  \left( \vartheta_u \right)_{(k)}} \right) \displaybreak[2]\\
&=  u\tr \left( \nabla^k_u  (a_0) e^{-\sigma_0 \left( \vartheta_u \right)_{(k)}} \nabla_u^k (a_1) \cdots e^{-\sigma_n  \left( \vartheta_u \right)_{(k)}} \right) \\
&\quad + \sum_{i=0}^{n}  (-1)^{i+i(n+1-i)+i(n-i)}  \sigma_i (\Theta_u)_{(k)} \wedge \tr \left(  a_0 e^{-\sigma_0 \left( \vartheta_u \right)_{(k)}} \nabla_u^k (a_1) \cdots e^{-\sigma_n  \left( \vartheta_u \right)_{(k)}}\right) \\
&\quad + \sum_{i=0}^{n-1}  (-1)^i  \tr \left(  a_0 e^{-\sigma_0 \left( \vartheta_u \right)_{(k)}} \cdots \nabla^k_u (a_i) e^{-\sigma_i \left( \vartheta_u \right)_{(k)}} [(\vartheta_u)_{(k)} , a_{i+1}] \cdots e^{-\sigma_n  \left( \vartheta_u \right)_{(k)}} \right)\displaybreak[2] \\
\intertext{and as $i+i(n+1-i)+i(n-i)$ is even,}
&=  u\tr \left( \nabla^k_u  (a_0) e^{-\sigma_0 \left( \vartheta_u \right)_{(k)}} \nabla_u^k (a_1) \cdots e^{-\sigma_n  \left( \vartheta_u \right)_{(k)}} \right) \\
&\quad + (\Theta_u)_{(k)} \wedge \tr \left(  a_0 e^{-\sigma_0 \left( \vartheta_u \right)_{(k)}} \nabla_u^k (a_1) \cdots e^{-\sigma_n  \left( \vartheta_u \right)_{(k)}}\right) \\
&\quad + \sum_{i=0}^{n-1}  (-1)^i  \tr \left(  a_0 e^{-\sigma_0 \left( \vartheta_u \right)_{(k)}} \cdots \nabla^k_u (a_i) e^{-\sigma_i \left( \vartheta_u \right)_{(k)}} [(\vartheta_u)_{(k)} , a_{i+1}] \cdots e^{-\sigma_n  \left( \vartheta_u \right)_{(k)}} \right) 
\end{align*}
where we have used $u (\nabla_u^k)^2 (a) = [(\vartheta_u)_{(k)} , a]$, $u \nabla^k_u \left( \left( \vartheta_u \right)_{(k)} \right) =-\left( \Theta_u \right)_{(k)}$ and
\begin{equation}
u \nabla^k_u \left( e^{-\sigma_i \left( \vartheta_u \right)_{(k)}} \right) = - \left(d \sigma_i \left( \vartheta_u \right)_{(k)} - \sigma_i  \left( \Theta_u \right)_{(k)} \right) e^{-\sigma_i \left( \vartheta_u \right)_{(k)}}
\end{equation}
as well as $\sum_{i=0}^n d \sigma_i = 0$ and $\sum_{i=0}^n \sigma_i = 1$ on $\Delta^n$.  Therefore, by subtracting $\tau_\nabla (  \Theta_u \wedge \omega ) $ from $\tau_\nabla ( (u \tilde{d}_{dR} + d_\Delta))$, we obtain a formula for $\tau_\nabla (d_{\Theta_u} \omega) (a_0 , \cdots , a_n)$
\begin{align}
& - \sum_k  \int_M \int_{\Delta^k}  \omega_{(k)}  \wedge \left( \int_{\Delta^n} u \tr (\nabla^k_u (a_0) e^{-\sigma_0 \left( \vartheta_u \right)_{(k)}}  \cdots e^{-\sigma_n  \left( \vartheta_u \right)_{(k)}} )  d \sigma_1 \cdots d \sigma_n \right) \label{eq:JLO_formula_B1} \\
&\quad - \int_M \int_{\Delta^k}  \omega_{(k)}  \wedge \bigg( \int_{\Delta^n} \sum_{i=0}^{n-1}  (-1)^i  \tr \Big(  a_0 e^{-\sigma_0 \left( \vartheta_u \right)_{(k)}} \cdots \notag\\
&\qquad \qquad \cdots \nabla^k_u (a_i) e^{-\sigma_i \left( \vartheta_u \right)_{(k)}} [(\vartheta_u)_{(k)} , a_{i+1}] \cdots e^{-\sigma_n  \left( \vartheta_u \right)_{(k)}} \Big) \bigg) \label{eq:JLO_formula_b1} \\
&\quad + (-1)^n \int_M \int_{\partial(\Delta^k)} \omega_{(k)} \wedge \bigg( \int_{\Delta^n} \tr \big(a_0 e^{-\sigma_0 \left( \vartheta_u \right)_{(k)}} \nabla_u^k (a_1)  \cdots e^{-\sigma_n  \left( \vartheta_u \right)_{(k)}} \big) d \sigma_1 \cdots d \sigma_n \bigg) \label{eq:JLO_formula_G1} 
\end{align}
Next we compute the effect of the algebraic operators.  The first term above, \eqref{eq:JLO_formula_B1}, is equivalent to $-(uB) \tau_\nabla (\omega) (a_0 , \ldots , a_n)$ as
\begin{align*}
&(uB) \tau_\nabla (\omega) (a_0 , \ldots , a_n) \displaybreak[2]\\
&= \sum_{i=0}^n (-1)^{ni} u \tau_\nabla (\omega) (1, a_i , \ldots , a_n , a_0 , \ldots a_{i-1}) \displaybreak[2] \\
&= \sum_{i=0}^n (-1)^{ni} u \sum_k \int_M \int_{\Delta^k} \omega_{(k)} \wedge \bigg( \int_{\Delta^{n+1}} \tr \Big( e^{-\sigma_0 \left( \vartheta_u \right)_{(k)}} \nabla_u^k (a_i) e^{-\sigma_1  \left( \vartheta_u \right)_{(k)}} \cdots  \displaybreak[0]\\ 
&\quad \cdots \nabla_u^k (a_n) e^{-\sigma_{n-i+1}  \left( \vartheta_u \right)_{(k)}} \nabla_u^k (a_0) e^{-\sigma_{n-i+2}  \left( \vartheta_u \right)_{(k)}} \cdots \nabla_u^k (a_{i-1}) e^{-\sigma_{n+1}  \left( \vartheta_u \right)_{(k)}} \Big) d \sigma_1 \cdots d \sigma_{n+1} \bigg) \displaybreak[2]\\
&= \sum_{i=0}^n (-1)^{ni} u \sum_k \int_M \int_{\Delta^k} \omega_{(k)} \wedge \bigg( \int_{\Delta^{n+1}} \tr \Big(  \nabla_u^k (a_i) e^{-\sigma_1  \left( \vartheta_u \right)_{(k)}} \cdots  \\ 
&\quad \cdots  \nabla_u^k (a_n) e^{-\sigma_{n-i+1}  \left( \vartheta_u \right)_{(k)}} \nabla_u^k (a_0) e^{-\sigma_{n-i+2}  \left( \vartheta_u \right)_{(k)}} \cdots  \nabla_u^k (a_{i-1}) e^{-(\sigma_0+\sigma_{n+1})  \left( \vartheta_u \right)_{(k)}} \Big) d \sigma_1 \cdots d \sigma_{n+1} \bigg) \displaybreak[2]\\
&= \sum_{i=0}^n (-1)^{ni+(n-i+1)i} u \sum_k \int_M \int_{\Delta^k} \omega_{(k)} \wedge \bigg( \int_{\Delta^{n+1}} \tr \Big( \nabla_u^k (a_0) e^{-\sigma_{n-i+2}  \left( \vartheta_u \right)_{(k)}} \cdots  \\ 
&\quad  \cdots \nabla_u^k (a_{i-1}) e^{-(\sigma_0+\sigma_{n+1})  \left( \vartheta_u \right)_{(k)}} \nabla_u^k (a_i) e^{-\sigma_1  \left( \vartheta_u \right)_{(k)}} \cdots  \nabla_u^k (a_n) e^{-\sigma_{n-i+1}  \left( \vartheta_u \right)_{(k)}}  \Big) d \sigma_1 \cdots d \sigma_{n+1} \bigg) \displaybreak[2]\\
&= \sum_{i=0}^n  u \sum_k \int_M \int_{\Delta^k} \omega_{(k)} \wedge \bigg( \int_{\Delta^{n}} \sigma_i \tr \Big( \nabla_u^k (a_0) e^{-\sigma_{0}  \left( \vartheta_u \right)_{(k)}} \nabla^k_u (a_1) \cdots \\
&\quad \cdots \nabla_u^k (a_n) e^{-\sigma_{n}  \left( \vartheta_u \right)_{(k)}}  \Big) d \sigma_1 \cdots d \sigma_{n} \bigg) \displaybreak[2]\\
&= u \sum_k \int_M \int_{\Delta^k} \omega_{(k)} \wedge \bigg( \int_{\Delta^{n}} \tr \big( \nabla_u^k (a_0) e^{-\sigma_{0}  \left( \vartheta_u \right)_{(k)}} \nabla^k_u (a_1) \cdots \nabla_u^k (a_n) e^{-\sigma_{n}  \left( \vartheta_u \right)_{(k)}}  \Big) d \sigma_1 \cdots d \sigma_{n} \bigg)
\end{align*}
where we have changed variables, performed a partial integration, and then used that $\sum \sigma_i = 1$.

The second term in $\tau_\nabla (d_{\Theta_u} \omega) (a_0 , \ldots , a_n)$, \eqref{eq:JLO_formula_b1}, is $-b \tau_\nabla (\omega) (a_0 , \ldots , a_n)$ as
\begin{align*}
& b \tau_\nabla (\omega) (a_0 , \ldots , a_n) \\
&= \sum_{i=0}^{n-1} (-1)^i \tau_\nabla (a_0 , \ldots , a_i a_{i+1} , \ldots , a_n) + (-1)^{n} \tau_\nabla (a_n a_0 , a_1 , \ldots , a_{n-1}) \displaybreak[2]\\
&=  \sum_k \int_M \int_{\Delta^k} \omega_{(k)} \wedge  \bigg( \int_{\Delta^{n-1}} \tr \Big(\nabla^k_u (a_0 a_1) e^{-\sigma_0 \left( \vartheta_u \right)_{(k)}} \nabla_u^k (a_2)  \cdots \\
&\quad \cdots \nabla_u^k (a_n) e^{-\sigma_{n-1}  \left( \vartheta_u \right)_{(k)}} \Big) d \sigma_1 \cdots d \sigma_{n-1} \bigg) \displaybreak[0]\\
&\quad + \sum_{i=1}^{n-1} (-1)^i \sum_k \int_M \int_{\Delta^k} \omega_{(k)} \wedge  \Big( \int_{\Delta^{n-1}} \big( \tr (a_0 e^{-\sigma_0 \left( \vartheta_u \right)_{(k)}} \nabla_u^k (a_1) \cdots \\
&\qquad \qquad  \cdots e^{-\sigma_{i-1} \left( \vartheta_u \right)_{(k)}} \nabla^k_u (a_i a_{i+1}) e^{-\sigma_{i} \left( \vartheta_u \right)_{(k)}} \cdots \nabla^k_u (a_n) e^{-\sigma_{n-1}  \left( \vartheta_u \right)_{(k)}} \big) d \sigma_1 \cdots d \sigma_{n-1} \Big) \displaybreak[0]\\
& \quad + (-1)^n  \sum_k \int_M \int_{\Delta^k} \omega_{(k)} \wedge  \Big( \int_{\Delta^{n-1}} \tr \big(\nabla^k_u (a_n a_0) e^{-\sigma_0 \left( \vartheta_u \right)_{(k)}} \nabla_u^k (a_1)  \cdots   \\
& \qquad \qquad \cdots  \nabla_u^k (a_{n-1}) e^{-\sigma_{n-1}  \left( \vartheta_u \right)_{(k)}} \big) d \sigma_1 \cdots d \sigma_{n-1} \Big) \displaybreak[2]\\
\intertext{and because $\nabla^k_u (a_i a_{i+1}) = \nabla^k_u (a_i) a_{i+1} + a_i \nabla^k (a_{i+1})$ we have}
&= \sum_k \int_M \int_{\Delta^k} \omega_{(k)} \wedge  \Big( \int_{\Delta^{n-1}} \tr \big(\nabla^k_u (a_0) a_1 e^{-\sigma_0 \left( \vartheta_u \right)_{(k)}} \nabla_u^k (a_2)  \cdots \\
&\cdots \nabla_u^k (a_n) e^{-\sigma_{n-1}  \left( \vartheta_u \right)_{(k)}} \big) d \sigma_1 \cdots d \sigma_{n-1} \Big) \displaybreak[1]\\
&\quad + \sum_k \int_M \int_{\Delta^k} \omega_{(k)} \wedge  \Big( \int_{\Delta^{n-1}} \tr \big(a_0 \nabla^k_u ( a_1) e^{-\sigma_0 \left( \vartheta_u \right)_{(k)}} \nabla_u^k (a_2)  \cdots \\
&\cdots \nabla_u^k (a_n) e^{-\sigma_{n-1}  \left( \vartheta_u \right)_{(k)}} \big) d \sigma_1 \cdots d \sigma_{n-1} \Big) \displaybreak[1]\\
&\quad + \sum_{i=1}^{n-1} (-1)^i \sum_k \int_M \int_{\Delta^k} \omega_{(k)} \wedge  \Big( \int_{\Delta^{n-1}} \tr \big(a_0 e^{-\sigma_0 \left( \vartheta_u \right)_{(k)}} \nabla_u^k (a_1) \cdots  \\
&\qquad \qquad \cdots e^{-\sigma_{i-1} \left( \vartheta_u \right)_{(k)}} \nabla^k_u (a_i) a_{i+1} e^{-\sigma_{i} \left( \vartheta_u \right)_{(k)}} \cdots \nabla^k_u (a_n) e^{-\sigma_{n-1}  \left( \vartheta_u \right)_{(k)}} \big) d \sigma_1 \cdots d \sigma_{n-1} \Big) \displaybreak[1]\\
&\quad + \sum_{i=1}^{n-1} (-1)^i \sum_k \int_M \int_{\Delta^k} \omega_{(k)} \wedge  \Big( \int_{\Delta^{n-1}} \tr \big(a_0 e^{-\sigma_0 \left( \vartheta_u \right)_{(k)}} \nabla_u^k (a_1) \cdots  \\
&\qquad \qquad \cdots e^{-\sigma_{i-1} \left( \vartheta_u \right)_{(k)}} a_i \nabla^k_u (a_{i+1})  e^{-\sigma_{i} \left( \vartheta_u \right)_{(k)}} \cdots \nabla^k_u (a_n) e^{-\sigma_{n-1}  \left( \vartheta_u \right)_{(k)}} \big) d \sigma_1 \cdots d \sigma_{n-1} \Big) \displaybreak[1]\\
& \quad + (-1)^n  \sum_k \int_M \int_{\Delta^k} \omega_{(k)} \wedge  \Big( \int_{\Delta^{n-1}} \tr \big(\nabla^k_u (a_n) a_0 e^{-\sigma_0 \left( \vartheta_u \right)_{(k)}} \nabla_u^k (a_1)  \cdots  \\
& \qquad \qquad \cdots \nabla_u^k (a_{n-1}) e^{-\sigma_{n-1}  \left( \vartheta_u \right)_{(k)}} \big) d \sigma_1 \cdots d \sigma_{n-1} \Big)\displaybreak[1] \\
& \quad + (-1)^n  \sum_k \int_M \int_{\Delta^k} \omega_{(k)} \wedge  \Big( \int_{\Delta^{n-1}} \tr \big(a_n\nabla^k_u ( a_0) e^{-\sigma_0 \left( \vartheta_u \right)_{(k)}} \nabla_u^k (a_1)  \cdots \\
& \qquad \qquad \cdots  \nabla_u^k (a_{n-1}) e^{-\sigma_{n-1}  \left( \vartheta_u \right)_{(k)}} \big) d \sigma_1 \cdots d \sigma_{n-1} \Big) \displaybreak[2] \\
\intertext{and by combining terms as well as commuting terms in the trace of the last two terms}
&= \sum_{i=1}^{n} (-1)^{i-1} \sum_k \int_M \int_{\Delta^k} \omega_{(k)} \wedge  \Big( \int_{\Delta^{n-1}} \tr \big(a_0 e^{-\sigma_0 \left( \vartheta_u \right)_{(k)}} \nabla_u^k (a_1) \cdots \\
&\qquad \qquad \cdots e^{-\sigma_{i-2} \left( \vartheta_u \right)_{(k)}} \nabla^k_u (a_{i-1}) [a_{i}, e^{-\sigma_{i-1} \left( \vartheta_u \right)_{(k)}} ] \nabla^k_u (a_{i+1}) \cdots \\
&\qquad \qquad \cdots \nabla^k_u (a_n) e^{-\sigma_{n-1}  \left( \vartheta_u \right)_{(k)}} \big) d \sigma_1 \cdots d \sigma_{n-1} \Big) 
\end{align*}
As $[a_i , \cdot]$ is a derivation, we can use the integration formula of \cite{quillen}, Equation 7.2, to obtain
\begin{align*}
[a_i , e^{- \sigma_i  \left( \vartheta_u \right)_{(k)}} ] &= \int_0^1 e^{- (1-\sigma) \sigma_i  \left( \vartheta_u \right)_{(k)}} [a_i , - \sigma_i  \left( \vartheta_u \right)_{(k)} ] e^{- \sigma \sigma_i  \left( \vartheta_u \right)_{(k)}} d \sigma \\
&= \sigma_i \int_0^1 e^{- (1-\sigma) \sigma_i  \left( \vartheta_u \right)_{(k)}} [   \left( \vartheta_u \right)_{(k)} , a_i ] e^{- \sigma \sigma_i  \left( \vartheta_u \right)_{(k)}} d \sigma 
\end{align*}
so that we have
\begin{align*}
&= \sum_{i=1}^{n} (-1)^{i-1} \sum_k \int_M \int_{\Delta^k} \omega_{(k)} \wedge  \left( \int_{\Delta^{n-1}} \int_0^1 \sigma_i \tr (a_0 e^{-\sigma_0 \left( \vartheta_u \right)_{(k)}} \nabla_u^k (a_1) \cdots \right. \\
&\qquad \qquad \left. \cdots e^{-\sigma_{i-2} \left( \vartheta_u \right)_{(k)}} \nabla^k_u (a_{i-1}) e^{- (1-\sigma) \sigma_i  \left( \vartheta_u \right)_{(k)}} [   \left( \vartheta_u \right)_{(k)} , a_i ] e^{- \sigma \sigma_i  \left( \vartheta_u \right)_{(k)}} \nabla^k_u (a_{i+1}) \cdots \right. \\
&\qquad \qquad \left. \cdots \nabla^k_u (a_n) e^{-\sigma_{n-1}  \left( \vartheta_u \right)_{(k)}} ) d \sigma d \sigma_1 \cdots d \sigma_{n-1} \right) 
\end{align*}
Now by the change of variables in the $i$-th term, $(1-\sigma)\sigma_i = \sigma_{i-1}'$ and $\sigma_i' = \sigma_{i-1} \sigma$ (as well as $\sigma_0 = \sigma_0', \ldots , \sigma_{i-2} = \sigma_{i-2}'$ and $\sigma_i = \sigma_{i+1}' , \ldots , \sigma_{n-1} = \sigma_n '$) we have $\sigma_{i-1} = \sigma_{i-1}' + \sigma_i'$, and in fact $\sigma_0 ' + \cdots + \sigma_n ' =1$ and $d\sigma_0' + \cdots + d\sigma_n' = 0$, i.e.\ barycentric coordinates on $\Delta^n$.  So with this change of variables the above expression is equivalent to \eqref{eq:JLO_formula_b1}.

The third term in $\tau_\nabla (d_{\Theta_u} \omega)$, \eqref{eq:JLO_formula_G1}, equals ${\delta_\Gamma}' \tau_\nabla (\omega) (a_0 , \ldots , a_n)$ because $\{ \vartheta_{(k)} \}$ and $\{ \omega_{(k)} \}$ are simplicial forms.  In particular, if we define Cartesian coordinates on $\Delta^k$ by $t_1 , \ldots , t_k$ then the restriction of $\vartheta_{(k)}$ to the face $t_i = 0$ is a form on $\Delta^{k-1}$ given by
\begin{equation*}
\vartheta_{(k)} (g_1 , \ldots , g_k) |_{t_i=0} = \vartheta_{(k-1)} (g_1 , \ldots , g_i g_{i+1} , \ldots , g_k)
\end{equation*}
for $1 \leq i \leq k-1$ and 
\begin{equation*}
\vartheta_{(k)} (g_1 , \ldots , g_k) |_{t_k=0} = \vartheta_{(k-1)} (g_1 , \ldots , g_{k-1})
\end{equation*}
The restriction to the face $\sum_{i=1}^k t_i = 1$ is a form on $\Delta^{k-1}$ given by
\begin{equation*}
\vartheta_{(k)} (g_1 , \ldots , g_k) |_{\sum t_i=1} = \vartheta_{(k-1)} (g_2 , \ldots , g_k)^{g_1}
\end{equation*}
as shown in Theorem \ref{thm:simplicial_2_form}.  Hence the corresponding identities hold for $\vartheta_u$. We also have the same sort of identities for $\{ \omega_{(k)} \}$ because $\{ \omega_{(k)} \}$ is a simplicial form. In addition $\nabla^k (g_1 , \ldots , g_k)$ satisfies the following identities for $1 \leq i \leq k-1$
\begin{align*}
\nabla^k (g_1 , \ldots , g_k) |_{t_i =0} &= \nabla^\E + t_1 A(g_1) + \cdots + t_{i-1} A(g_1 \cdots g_{i-1}) + t_{i+1} A(g_1 \cdots g_{i+1}) + \cdots \\
&\qquad \qquad  \cdots + t_k A(g_1 \cdots g_k) \\
&= \nabla^{k-1} (g_1 , \ldots , g_i g_{i+1} , \ldots , g_k)
\end{align*}
with a change of variables.  Similarly $\nabla^k (g_1 , \ldots , g_k)|_{t_k=0} = \nabla^{k-1} (g_1 , \ldots , g_{k-1})$. Furthermore,
\begin{align*}
\nabla^k (g_1 , \ldots , g_k) |_{\sum t_i = 1} &= \nabla^\E + A ( g_1) + t_2 (A (g_1 g_2) - A(g_1)) + \cdots  + t_k ( A(g_1 \cdots g_k) - A(g_1)) \\
&= \varphi_{g_1}^* (\nabla^{\E^g} \otimes 1 + 1 \otimes \nabla_g) + t_2 (\alpha( g_1 , g_2) + \varphi_{g_1}^* A(g_2)^{g_1}) + \cdots \\
&\qquad \qquad \cdots + t_k (\alpha (g_1 , g_2 \cdots g_k) + \varphi_{g_1}^* A(g_2 \cdots g_k) ) \\
&= \varphi_{g_1}^* (\nabla^{\E^g} + t_2 A(g_2) + \cdots + t_k A(g_2 \cdots g_k))^{g_1}  \\
&\quad + \varphi_{g_1}^* \nabla_g + t_2 \alpha (g_1 , g_2) + \cdots t_k \alpha(g_1 , g_2 \cdots g_k)\\
&= \varphi_{g_1}^* (\nabla^{\E^g} + t_2 A(g_2) + \cdots + t_k A(g_2 \cdots g_k))^{g_1}  \\
&\cong \nabla^{k-1} (g_2 , \ldots , g_k)^{g_1}
\end{align*}
by Lemma \ref{lem:discrepancy_identity} and the definition of $A(g_1)$. Since $\nabla_g + t_2 \alpha (g_1 , g_2) + \cdots t_k \alpha(g_1 , g_2 \cdots g_k)$ is scalar-valued,
\begin{equation*}
[\nabla^{k-1} (g_2 , \ldots , g_k)^{g_1} + \varphi_{g_1}^* \nabla_g + t_2 \alpha (g_1 , g_2) + \cdots t_k \alpha(g_1 , g_2 \cdots g_k) , \eta ]  = [\nabla^{k-1} (g_2 , \ldots , g_k)^{g_1} , \eta]
\end{equation*}
for any $\eta \in C^\infty_c (M , \End \E)$.  The corresponding identities hold for $\nabla^k_u$ as well.

Taking orientation into account we may write 
\begin{equation*}
\partial (\Delta^k) = \Delta^k |_{\sum t_i = 1} - \Delta^k |_{t_1 = 0} + \Delta^k |_{t_2 = 0}  \cdots + (-1)^{k+1} \Delta^k |_{t_k = 0}
\end{equation*}
where a negative sign denotes the orientation opposite the orientation induced by $\Delta^k$. With this the third term of $\tau_\nabla (d_{\Theta_u} \omega) (a_0 , \ldots , a_n)$, \eqref{eq:JLO_formula_G1}, may be written
\begin{align*}
&(-1)^n \sum_k  \int_M \int_{\Delta^k |_{\sum t_i = 1}} \omega_{(k)} \wedge \left( \int_{\Delta^n} \tr (a_0 e^{-\sigma_0 \left( \vartheta_u \right)_{(k)}} \nabla_u^k (a_1)  \cdots  e^{-\sigma_n  \left( \vartheta_u \right)_{(k)}} ) d \sigma_1 \cdots d \sigma_n \right) \\
&+ (-1)^n \sum_{i=1}^k (-1)^i \int_M \int_{\Delta^k |_{t_i = 0}} \omega_{(k)} \wedge \left( \int_{\Delta^n} \tr (a_0 e^{-\sigma_0 \left( \vartheta_u \right)_{(k)}}   \cdots  e^{-\sigma_n  \left( \vartheta_u \right)_{(k)}} ) d \sigma_1 \cdots d \sigma_n \right).
\end{align*}
And evaluation on $g_1 , \ldots , g_k$ gives
\begin{align*}
&(-1)^n \sum_k  \int_M \int_{\Delta^{k-1}} \omega_{(k-1)} (g_2 , \ldots , g_k)^{g_1} \wedge \left( \int_{\Delta^n} \tr (a_0 e^{-\sigma_0 \left( \left( \vartheta_u \right)_{(k-1)} (g_2 , \ldots , g_k) \right)^{g_1} } \right. \\
&\qquad \qquad \left. (\nabla_u^{k-1} (g_2 , \ldots , g_k) (a_1))^{g_1} \cdots  e^{-\sigma_n \left( \left( \vartheta_u \right)_{(k-1)} (g_2 , \ldots , g_k) \right)^{g_1}}  ) d \sigma_1 \cdots d \sigma_n \right)  \displaybreak[2]\\
&+ (-1)^n \sum_{i=1}^{k-1} (-1)^i \int_M \int_{\Delta^{k-1}} \omega_{(k-1)} (g_1 , \ldots , g_i g_{i+1} , \ldots , g_k) \wedge  \\
&\qquad \qquad  \wedge \left( \int_{\Delta^n} \tr (a_0 e^{-\sigma_0 \left( \vartheta_u \right)_{(k-1)} (g_1 , \ldots , g_i g_{i+1} , \ldots , g_k) } (\nabla_u^{k-1} (g_1 , \ldots , g_i g_{i+1} , \ldots , g_k)) (a_1)  \cdots \right. \\
&\qquad \qquad \cdots \left. e^{-\sigma_n  \left( \vartheta_u \right)_{(k-1)} (g_1 , \ldots , g_i g_{i+1} , \ldots , g_k) } ) d \sigma_1 \cdots d \sigma_n \right)  \displaybreak[2]\\
&\quad + (-1)^n (-1)^{k} \int_M \int_{\Delta^{k-1}} \omega_{(k-1)} (g_1 , \ldots , g_{k-1}) \wedge  \\
&\qquad \qquad  \wedge \left( \int_{\Delta^n} \tr (a_0 e^{-\sigma_0 \left( \vartheta_u \right)_{(k-1)} (g_1 , \ldots , g_{k-1}) } (\nabla_u^{k-1} (g_1 , \ldots , g_{k-1})) (a_1)  \cdots \right. \\
&\qquad \qquad \cdots \left. e^{-\sigma_n  \left( \vartheta_u \right)_{(k-1)} (g_1 , \ldots , g_{k-1}) } ) d \sigma_1 \cdots d \sigma_n \right)  \displaybreak[2]\\
&= (-1)^n \tau_\nabla (\omega) (g_2 , \ldots , g_k)^{g_1} (a_0 , \ldots , a_n) \\
&\quad + (-1)^n \sum_{i=1}^{k-1} \tau_\nabla (\omega) (g_1 , \ldots , g_i g_{i+1}, \ldots , g_k) (a_0 , \ldots , a_n)\\
&\quad + (-1)^n (-1)^k \tau_\nabla (\omega) (g_1 , \ldots , g_{k-1}) (a_0 , \ldots , a_n) \\
&={\delta_\Gamma}' \tau_\nabla (\omega) (g_1 , \ldots , g_k) (a_0 , \ldots , a_n)
\end{align*}
This proves that $\tau_\nabla \circ d_{\Theta_u} = (-(b+ uB) + {\delta_\Gamma}') \circ \tau_\nabla$.
\end{proof}
\end{thm}


\subsection{Algebraic morphisms}
\label{sec:algebraic_morphisms}

Now we will need a different cochain complex computing group cohomology.  For a discrete group $\Gamma$ and a $\Gamma$-module $K$, let $\widetilde{C}^n (\Gamma , K)$ denote the space of \textit{homogenous $n$-cochains of $\Gamma$ with values in $K$}.  This is the space of degree $n+1$ $\Gamma$-cochains, i.e.\ the space of maps $\Gamma^{n+1} \rightarrow K$, that satisfy the condition
\begin{equation*}
g f (g_0 , \ldots , g_n) = f (g g_0 , \ldots , g g_n).
\end{equation*}
For any $f \in \widetilde{C}^n (\Gamma , K)$ define $\widetilde{\delta}_\Gamma : \widetilde{C}^n (\Gamma , K) \to \widetilde{C}^{n+1} (\Gamma , K)$ by 
\begin{equation}
({\widetilde{\delta}_\Gamma} f) (g_0 , \ldots , g_n) = \sum_{i=0}^n (-1)^i f(g_0 , \ldots , \hat{g_i} , \ldots , g_n)
\end{equation}
and let ${\widetilde{\delta}_\Gamma}' = (-1)^n \widetilde{\delta}_\Gamma$.  There is a bijection from $C^n (\Gamma , K) \to \widetilde{C}^n (\Gamma ,K)$ given by $f \mapsto \widetilde{f}$ where $\widetilde{f}$ is defined by
\begin{equation}
\widetilde{f} (g_0 , \ldots , g_n) = f(g_0^{-1} g_1 , (g_0 g_1)^{-1} g_2 , \ldots , (g_0 \cdots g_{n-1})^{-1} g_n )
\end{equation}
for any $f \in C^n (\Gamma , K)$ and $g_0 , \ldots , g_n \in \Gamma$.  The complex $(\widetilde{C} (\Gamma ,K) , {\widetilde{\delta}_\Gamma})$ also computes group cohomology. 

\begin{prop}

The above bijection induces a map
\begin{multline}
\Psi_0 : (C^\bullet (\Gamma , CC^\bullet ( \Gamma^\infty_c (\End \E) )  [u^{-1} , u]) , b+uB + {\delta_\Gamma}') \\ \to (\widetilde{C}^\bullet (\Gamma , CC^\bullet ( \Gamma^\infty_c ( \End \E) ) [u^{-1} , u] ) , b+uB + {\widetilde{\delta}_\Gamma}') 
\end{multline}
by the formula
\begin{equation}
\Psi_0 (c) (g_0 , \ldots , g_n) = c (g_0^{-1} g_1 , (g_0 g_1)^{-1} g_2 , \ldots , (g_0 \cdots g_{n-1})^{-1} g_n)
\end{equation}
for any $c \in C^\bullet (\Gamma , CC^\bullet ( \Gamma^\infty_c (\End \E) ) [u^{-1} , u] )$.

\end{prop}

\begin{lem}
There is a contracting homotopy
\begin{multline*}
h : (\widetilde{C}^k (\Gamma , CC^n ( \Gamma^\infty_c ( \End \E) )  [u^{-1} , u] ) , {\widetilde{\delta}_\Gamma} , (-1)^k (b+uB) ) \\ \to (\widetilde{C}^{k-1} (\Gamma , CC^n ( \Gamma^\infty_c ( \End \E) )  [u^{-1} , u] ) , {\widetilde{\delta}_\Gamma} , (-1)^k (b+uB) )
\end{multline*}

\begin{proof}

First we will define a map $\gamma: \Gamma^\infty_c (\End \E)^{\otimes n} \to \C \Gamma$ by
\begin{equation}
\gamma : E_{g_0 , g_1} (a_0) \otimes E_{g_1 , g_2} (a_1) \otimes \ldots \otimes E_{g_n , g_0} (a_n) \mapsto g_0
\end{equation}
where $E_{g_i , g_{i+1}} (a_i), E_{g_n , g_0} (a_n) \in \Gamma^\infty_c ( \End \E)$ for $0 \leq i \leq n-1$.  Note that by multilinearity,
\begin{equation}
\gamma ( E_{g_0 , g_1} (a_0) + E_{g_0' , g_1} (a_0') , E_{g_1 , g_2} (a_1) , \ldots , E_{g_n , g_0} (a_n) + E_{g_n , g_0'} (a_n') ) = g_0 + g_0' .
\end{equation}
Furthermore, the map $\gamma$ is equivariant:
\begin{align*}
\gamma ( g \cdot (  E_{g_0 , g_1} (a_0)  \otimes \ldots \otimes E_{g_n , g_0} (a_n)  ) ) &= \gamma ( E_{gg_0 , gg_1} (a_0^g)  \otimes \ldots \otimes E_{gg_n , gg_0} (a_n^g) ) \\
&= g g_0 \\
&= g \gamma ( E_{g_0 , g_1} (a_0)  \otimes \ldots \otimes E_{g_n , g_0} (a_n) )
\end{align*}

Then we can define a map
\begin{equation}
h : \widetilde{C}^{k+1} (\Gamma , CC^n ( \Gamma^\infty_c ( \End \E) ) ) \to \widetilde{C}^{k} (\Gamma , CC^n ( \Gamma^\infty_c ( \End \E) ) )
\end{equation}
by
\begin{equation}
(h \varphi) (g_0 , \ldots , g_k) ( s ) = (-1)^k \varphi (g_0 , \ldots , g_k , \gamma (s)) (s)
\end{equation}
for any $s \in \Gamma^\infty_c ( \End \E)^{\otimes n}$.

The map $h$ is a contracting homotopy for group cochains as we have
\begin{align*}
({\widetilde{\delta}_\Gamma} h) f (g_0 , \ldots , g_k) (s) &= \sum_{i=0}^k (-1)^i (hf) (g_0 , \ldots , \hat{g_i} , \ldots , g_k) (s) \\
&= (-1)^{k-1} \sum_{i=0}^k (-1)^i f (g_0 , \ldots , \hat{g_i} , \ldots , g_k , \gamma (s) ) (s) \displaybreak[2]\\
(h {\widetilde{\delta}_\Gamma} ) f (g_0 , \ldots , g_k) (s) &= (-1)^k ({\widetilde{\delta}_\Gamma} f) (g_0 , \ldots , g_k , \gamma (s) ) (s) \\
&= (-1)^k \sum_{i=0}^k (-1)^i f( g_0 , \ldots , \hat{g_i} , \ldots , g_k , \gamma(s) ) (s)  \\
&\quad+ (-1)^{k+(k+1)} f( g_0 , \ldots , g_k) (s) 
\end{align*}
so that
\begin{equation}
({\widetilde{\delta}_\Gamma} h + h {\widetilde{\delta}_\Gamma})  = 1
\end{equation}
is the identity operator.
\end{proof}
\end{lem}

\begin{thm}
Given a gerbe datum $(L , \mu, \nabla)$ on $M \rtimes \Gamma$ and an $(L, \mu, \nabla)$-twisted vector bundle with connection $(\E, \varphi, \nabla^\E)$, there is a morphism
\begin{equation*}
\Psi_1 : (\widetilde{C}^\bullet (\Gamma , CC^\bullet ( \Gamma^\infty_c (\End \E) )  [u^{-1} , u] ) , b+uB + {\widetilde{\delta}_\Gamma}')  \to  ( CC^\bullet (\Gamma^\infty_c ( \End \E)) [u^{-1} , u], b+uB)^\Gamma
\end{equation*}
from the complex of cochains on $\Gamma$ taking values in the periodic cyclic complex of $\Gamma^\infty_c ( \End \E)$ to the periodic cyclic complex of $\Gamma$-invariant cochains on $\Gamma^\infty_c ( \End \E)$.

\begin{proof}

Now we will view $$(\widetilde{C}^\bullet (\Gamma , CC^\bullet ( \Gamma^\infty_c ( \End \E) )  [u^{-1} , u]) , b+uB + {\widetilde{\delta}_\Gamma}')$$ as a bicomplex $$(\widetilde{C}^k (\Gamma , CC^n ( \Gamma^\infty_c (\End \E) )  [u^{-1} , u] ) , {\widetilde{\delta}_\Gamma} , (-1)^k (b+uB) )$$ with horizontal differential $\widetilde{\delta}_\Gamma$ and vertical differential $b+uB$.  More explicitly this is done via a morphism
\begin{multline*}
\Psi_{1/2} : (\widetilde{C}^\bullet (\Gamma , CC^\bullet ( \Gamma^\infty_c ( \End \E) ) [u^{-1} , u] ) , b+uB + {\widetilde{\delta}_\Gamma}') \\ \to (\widetilde{C}^k (\Gamma , CC^n ( \Gamma^\infty_c ( \End \E) )  [u^{-1} , u] ) , {\widetilde{\delta}_\Gamma} , (-1)^k (b+uB) )
\end{multline*}
which is defined by exchanging the vertical and horizontal differentials and making the appropriate sign adjustments ${\widetilde{\delta}_\Gamma}' \mapsto {\widetilde{\delta}_\Gamma}$ and $(b+ uB) \mapsto (-1)^k (b+uB)$.

Since the map $h$ is a contracting homotopy for the rows of this bicomplex, i.e.\ the rows of the augmented double complex are exact, the cohomology of the bicomplex is the cohomology of the initial column.  As described in \cite{bott_tu} in the case of the \v{C}ech-de Rham complex, an explicit morphism
\begin{multline*}
\Psi_1^k : (\widetilde{C}^k (\Gamma , CC^n ( \Gamma^\infty_c (\End \E) )  [u^{-1} , u]) , {\widetilde{\delta}_\Gamma} , (-1)^k (b+uB) ) \\ \to  (CC^{n+k} (\Gamma^\infty_c ( \End \E))  [u^{-1} , u] , b+uB)^\Gamma
\end{multline*}
is given by the formula
\begin{equation*}
\Psi_1^k (c) = ( -(b+uB)h)^k c - h (-(b+uB) h)^{k-1}  ((b+uB) c) - h (-(b+uB) h)^k ( {\widetilde{\delta}_\Gamma} c ).
\end{equation*}
For cochains in column $k$ denote $D= (-1)^k (b+uB)$ for the sake of brevity.  Then this formula is
\begin{equation*}
\Psi_1^k (c) = ( -Dh)^k c - h (-D h)^{k-1}  (D c) - h (-D h)^k ( {\widetilde{\delta}_\Gamma} c ).
\end{equation*}
Then let $\Psi_1 = \sum_k \Psi^k$.  We will also denote $\Psi_1 \circ \Psi_{1/2}$ merely by $\Psi_1$, as in the theorem statement.  We can see immediately that if $c$ is a group $k$-cochain then 
\begin{equation*}
\Psi_1 ( c) \in  C^0 ( \Gamma , CC^{n+k} (\Gamma^\infty_c (\End \E))  [u^{-1} , u])
\end{equation*}
because 
\begin{equation*}
(Dh) : C^i ( \Gamma , CC^{j} (\Gamma^\infty_c ( \End \E))  [u^{-1} , u] ) \to C^{i-1} ( \Gamma , CC^{j+1} (\Gamma^\infty_c ( \End \E))  [u^{-1} , u]).
\end{equation*}

At this point let us recount a lemma from \cite{bott_tu}.
\begin{lem}
For $i \geq 1$
\begin{equation*}
{\widetilde{\delta}_\Gamma} (D h)^i = (D h)^i {\widetilde{\delta}_\Gamma} - (D h)^{i-1} D
\end{equation*}

\begin{proof}[Proof of lemma]
Since ${\widetilde{\delta}_\Gamma}$ anticommutes with $D$ and ${\widetilde{\delta}_\Gamma} h + h {\widetilde{\delta}_\Gamma} = 1$,
\begin{align*}
{\widetilde{\delta}_\Gamma} (D h)(D h)^{i-1} &= -D {\widetilde{\delta}_\Gamma} h (D h)^{i-1} \\
&= - D (1- h{\widetilde{\delta}_\Gamma}) (D h)^{i-1} \\
&= (D h) {\widetilde{\delta}_\Gamma} (D h)^{i-1}
\end{align*}
Therefore we can commute ${\widetilde{\delta}_\Gamma}$ and $(D h)$ until we reach $(D h)^{i-1} {\widetilde{\delta}_\Gamma} (Dh)$.  We obtain the formula by
\begin{align*}
{\widetilde{\delta}_\Gamma} (D h)(D h)^{i-1} &= (D h)^{i-1} {\widetilde{\delta}_\Gamma} (Dh) \\
&= - (Dh)^{i-1} D (1 - h {\widetilde{\delta}_\Gamma}) \\
&= - (Dh)^{i-1} D + (Dh)^i {\widetilde{\delta}_\Gamma}
\end{align*}
\end{proof}
\end{lem}

The first step in proving $\Psi_1$ is a morphism is to check that $\Psi_1 (c)$ does in fact land in the complex of $\Gamma$-invariant cyclic cochains $CC^{n+k} (\Gamma^\infty_c ( \End \E)  [u^{-1} , u] , b+uB)^\Gamma$.  To check this we compute ${\widetilde{\delta}_\Gamma} \Psi_1 (\alpha) $. By the above lemma
\begin{align*}
{\widetilde{\delta}_\Gamma} \Psi_1^k (c) &= {\widetilde{\delta}_\Gamma} (-Dh)^k c - {\widetilde{\delta}_\Gamma} h (-D h)^{k-1} D c - {\widetilde{\delta}_\Gamma} h (-D h)^k {\widetilde{\delta}_\Gamma} c \\
&= (-Dh)^k {\widetilde{\delta}_\Gamma} c + (-Dh)^{k-1} D c - {\widetilde{\delta}_\Gamma} h (-Dh)^{k-1} D c - {\widetilde{\delta}_\Gamma} h (-Dh)^k {\widetilde{\delta}_\Gamma} c \\
&= (1- {\widetilde{\delta}_\Gamma} h) (-Dh)^k {\widetilde{\delta}_\Gamma} c + (1- {\widetilde{\delta}_\Gamma} h) (-Dh)^{k-1} D c \\
&= (1- {\widetilde{\delta}_\Gamma} h) (-Dh)^{k-1} (-Dh {\widetilde{\delta}_\Gamma} + D) c \\
&= (1- {\widetilde{\delta}_\Gamma} h) (-Dh)^{k-1} (D (1- h {\widetilde{\delta}_\Gamma} )) c \\
&= ( h {\widetilde{\delta}_\Gamma} ) (-Dh)^{k-1} (D  {\widetilde{\delta}_\Gamma} h ) c \\
\intertext{and by applying the lemma again}
&= h ( (-Dh)^{k-1} {\widetilde{\delta}_\Gamma} + (-Dh)^{k-2} D) (D {\widetilde{\delta}_\Gamma} h ) c \\
&= h ( (-Dh)^{k-1} {\widetilde{\delta}_\Gamma} (D {\widetilde{\delta}_\Gamma} h ) c  + h (-Dh)^{k-2} D) (D {\widetilde{\delta}_\Gamma} h ) c \\
&=0
\end{align*}
by $D^2 = 0$ and $({\widetilde{\delta}_\Gamma})^2 = 0$ as well as commuting $D$ and ${\widetilde{\delta}_\Gamma}$. The cases ${\widetilde{\delta}_\Gamma} \Psi_1^1 (c) = 0$ and ${\widetilde{\delta}_\Gamma} \Psi_1^0 (c) = 0$ are straightforward. Since ${\widetilde{\delta}_\Gamma} \Psi_1 (c) = 0$ we see that $\Psi_1 (c)$ is $\Gamma$-invariant.

Further, the boundary operators commute with $\Psi_1$. For any $c \in C^k ( \Gamma , CC^{n} (\Gamma^\infty_c (\End \E))  [u^{-1} , u] )$
\begin{align*}
\Psi_1 ( ({\widetilde{\delta}_\Gamma} + D) c) &= \Psi_1^k (D c) + \Psi_1^{k+1} ({\widetilde{\delta}_\Gamma} c) \\
&= (-Dh)^k (D c) - h (-Dh )^{k-1} D (D c) - h (-Dh)^k {\widetilde{\delta}_\Gamma} (D c) \\
&\quad + (-Dh)^{k+1} ({\widetilde{\delta}_\Gamma} c) - h (-Dh )^{k} D ({\widetilde{\delta}_\Gamma} c) - h (-Dh)^{k+1} {\widetilde{\delta}_\Gamma} ({\widetilde{\delta}_\Gamma} c) \\
\intertext{which, by the commutation relations for the boundary operators,}
&= (-Dh)^k D c + (-Dh)^{k+1} {\widetilde{\delta}_\Gamma} c \\
&= -Dh (-Dh)^{k-1} D c - Dh (-Dh)^k {\widetilde{\delta}_\Gamma} c \\
&= D ( (-Dh)^k c - h (-Dh )^{k-1} D c - h(-Dh)^k {\widetilde{\delta}_\Gamma} c) \\
&= D \Psi_1^k (c) = D \Psi_1 (c)
\end{align*}
Hence $\Psi_1$ is a morphism.

\end{proof}

\end{thm}

\begin{thm}
Given a gerbe datum $(L, \mu, \nabla)$ on $M \rtimes \Gamma$ with an $(L, \mu, \nabla)$-twisted vector bundle $(\E,\varphi,\nabla^\E)$, there is a morphism 
\begin{equation}
\Psi_2 : CC^\bullet (\Gamma^\infty_c ( \End \E)  [u^{-1} , u], b+uB)^\Gamma \to CC^\bullet (C^\infty_c (M \rtimes \Gamma , L)  [u^{-1} , u], b+uB)
\end{equation}
defined by the formula
\begin{equation}
\label{eq:alg_morphism2}
\Psi_2 (c) (\widetilde{a_{g_0}} , \ldots , a_{g_n} ) = c(E_{1,g_0} (\widetilde{a_{g_0}}) , \ldots , E_{g_0 \cdots g_{i-1} , g_0 \cdots g_i} (a_{g_i}^{g_0 \cdots g_{i-1}}) , \ldots , E_{g_0 \cdots g_{n-1}, 1} (a_{g_n}^{g_0 \cdots g_{n-1}}))
\end{equation}
for any $g_0 , \ldots g_n \in \Gamma$. Here $c$ is a cochain in $CC^n (\Gamma_c^\infty (\End \E))^\Gamma$ and $a_{g_i}$ is a section in $\Gamma^\infty_c ( L_{g_i} )$ for $0 \leq i \leq n$, and $\widetilde{a_{g_0}} = (a_{g_0} , \lambda)$ for $\lambda \in \C$.
\begin{proof}
Note that $E_{1, g_0} (\widetilde{a_{g_0}}) = \widetilde{E_{1,g_0} (a_{g_0})}$ in $\widetilde{\Gamma^\infty_c ( \End \E)}$.  Since $a \in C^\infty_c (M \rtimes \Gamma, L)$ may be written $a = \sum_{g \in \Gamma} a_g$ where $a_g \in \Gamma^\infty_c ( L_g)$ for $g \in \Gamma$, it is sufficient to define $\Psi_2$ on such sections.  Also note that $(a_{g} * a_{g'} ) \in \Gamma^\infty_c ( L_{gg'})$ since 
\begin{align*}
(a_g * a_{g'}) (x, h) &= \sum_{h_1 h_2 = h} a_g ( x , h_1) \cdot a_{g'} (xh_1 , h_2) \\
&=a_g ( x , g) \cdot a_{g'} (x g , g') \\
&= \mu_{g , g'} (a_g \otimes a_{g'}^g) (x,g g').
\end{align*}
As the coboundary $b^n = \sum_{i=1}^n d^i_n$, we will prove that $d^i \Psi_2 = \Psi_2 d^i$ for $1\leq i \leq n$.  For $1 \leq i \leq n$,
\begin{align*}
&(d^i \Psi_2) ( c ) (\widetilde{a_{g_0}} , \ldots , a_{g_{n+1}} ) \displaybreak[2] \\
&= \Psi_2 (c) (\widetilde{a_{g_0}} , \ldots , a_{g_i} * a_{g_{i+1}} , \ldots , a_{g_{n+1}}) \displaybreak[2] \\
&=\Psi_2 (c) \left(\widetilde{a_{g_0}} , \ldots , \mu_{g_i , g_{i+1}} (a_{g_i} \otimes a_{g_{i+1}}^{g_i}) , \ldots , a_{g_{n+1}}\right) \displaybreak[2] \\
&=c \left( E_{1, g_0} (\widetilde{a_{g_0}})  , \ldots , E_{g_0 \cdots g_{i-1} , g_0 \cdots g_{i+1}} ((\mu_{g_i , g_{i+1}} (a_{g_i} \otimes a_{g_{i+1}}^{g_i} ))^{g_0 \cdots g_{i-1} } ) , \ldots , E_{g_0 \cdots g_n , 1} (a_{g_{n+1}}^{g_0 \cdots g_n}) \right)\displaybreak[2] \\
&=c \left( E_{1, g_0} (\widetilde{a_{g_0}})  , \ldots , E_{g_0 \cdots g_{i-1} , g_0 \cdots g_{i+1}} (\mu_{g_i , g_{i+1}}^{g_0 \cdots g_{i-1}} (a_{g_i}^{g_0 \cdots g_{i-1} } \otimes a_{g_{i+1}}^{g_0 \cdots g_i } ) ) , \ldots, E_{g_0 \cdots g_n , 1} (a_{g_{n+1}}^{g_0 \cdots g_n}) \right) \displaybreak[2] \\
&= c \left( E_{1, g_0} (\widetilde{a_{g_0}})  , \ldots , E_{g_0 \cdots g_{i-1} , g_0 \cdots g_i} (a_{g_i}^{g_0 \cdots g_{i-1}} ) E_{g_0 \cdots g_i , g_0 \cdots g_{i+1}} (a_{g_{i+1}}^{g_0 \cdots g_i} ) , \ldots, E_{g_0 \cdots g_n , 1} (a_{g_{n+1}}^{g_0 \cdots g_n}) \right) \displaybreak[2] \\
&= (d^i c) \left(E_{1, g_0} (\widetilde{a_{g_0}})  , \ldots , E_{g_0 \cdots g_{i-1} , g_0 \cdots g_i} (a_{g_i}^{g_0 \cdots g_{i-1}} ) , E_{g_0 \cdots g_i , g_0 \cdots g_{i+1}} (a_{g_{i+1}}^{g_0 \cdots g_i} ) , \ldots, E_{g_0 \cdots g_n , 1} (a_{g_{n+1}}^{g_0 \cdots g_n}) \right) \displaybreak[2] \\
&= (\Psi_2 d^i) (c) (\widetilde{a_{g_0}} , \ldots , a_{g_{n+1}} ).
\end{align*}
For the case of $d^0 \Psi_2 = \Psi_2 d^0$ note that $\widetilde{a_{g_0}} * a_{g_1} = a_{g_0} * a_{g_1} + \lambda a_{g_1}$ 
\begin{align*}
&(d^0 \Psi_2) ( c ) (\widetilde{a_{g_0}} , \ldots , a_{g_{n+1}} ) \displaybreak[2] \\
&= \Psi_2 (c) (\widetilde{a_{g_0}} * a_{g_1}, a_{g_2} , \ldots  , a_{g_{n+1}}) \displaybreak[2] \\
&= \Psi_2 (c) (a_{g_0} * a_{g_1} + \lambda a_{g_1}, a_{g_2} , \ldots  , a_{g_{n+1}}) \displaybreak[2] \\
&= \Psi_2 (c) (a_{g_0} * a_{g_1} , a_{g_2} , \ldots  , a_{g_{n+1}}) + \Psi_2 (c) ( \lambda a_{g_1}, a_{g_2} , \ldots  , a_{g_{n+1}}) \displaybreak[2] \\
&= \Psi_2 (c) (\mu_{g_0 , g_1} (a_{g_0} \otimes a_{g_1}^{g_0} ) , a_{g_2} , \ldots  , a_{g_{n+1}}) + \Psi_2 (c) ( \lambda a_{g_1},  a_{g_2} , \ldots  , a_{g_{n+1}}) \displaybreak[2] \\
&= c ( E_{1, g_0 g_1} (\mu_{g_0 , g_1} (a_{g_0} \otimes a_{g_1}^{g_0} )) , E_{g_0 g_1 , g_0 g_1 g_2} (a_{g_2}^{g_0 g_1}) , \ldots , E_{g_0 \dots g_n , 1} (a_{g_{n+1}}^{g_0 \cdots g_n})) \\
&\quad + c (\lambda E_{1, g_1} (a_{g_1}) , E_{g_1 , g_1 g_2} (a_{g_2}^{g_1}) , \ldots , E_{g_1 \cdots g_n ,  g_0^{-1} } (a_{g_{n+1}}^{g_1 \cdots g_n})) \displaybreak[2] \\
&= c ( E_{1, g_0 g_1} (\mu_{g_0 , g_1} (a_{g_0} \otimes a_{g_1}^{g_0} ) , E_{g_0 g_1 , g_0 g_1 g_2} (a_{g_2}^{g_0 g_1}) , \ldots , E_{g_0 \dots g_n , 1} (a_{g_{n+1}}^{g_0 \cdots g_n})) \\
&\quad + c (\lambda E_{1, g_1} (a_{g_1}) , E_{g_1 , g_1 g_2} (a_{g_2}^{g_1}) , \ldots , E_{g_1 \cdots g_n , 1} (a_{g_{n+1}}^{g_1 \cdots g_n})) \cdot g_0 \\
&= c ( E_{1, g_0 g_1} (\mu_{g_0 , g_1} (a_{g_0} \otimes a_{g_1}^{g_0} ) , E_{g_0 g_1 , g_0 g_1 g_2} (a_{g_2}^{g_0 g_1}) , \ldots , E_{g_0 \dots g_n , 1} (a_{g_{n+1}}^{g_0 \cdots g_n})) \\
&\quad + c (\lambda E_{g_0, g_0 g_1} (a_{g_1}^{g_0}) , E_{g_0 g_1 , g_0 g_1 g_2} (a_{g_2}^{g_0 g_1}) , \ldots , E_{g_0 g_1 \cdots g_n , g_0} (a_{g_{n+1}}^{g_0  \cdots g_n})) \\
\intertext{as the second term is nonzero only if $g_0 \cdots g_{n+1} =1$ and $g_1 \cdots g_{n+1} = 1$, i.e.\ $g_0=1$ }
&= c \left( E_{1,g_0} (a_{g_0}) E_{g_0 , g_0 g_1} (a_{g_1}^{g_0}) , E_{g_0 g_1 , g_0 g_1 g_2} (a_{g_2}^{g_0 g_1}),  \ldots  , E_{g_0 \cdots g_n , 1} (a_{g_{n+1}}^{g_0 \cdots g_n})  \right)\\
&\quad + c \left( \lambda E_{g_0 , g_0 g_1} (a_{g_1}^{g_0}) , E_{g_0 g_1 , g_0 g_1 g_2} (a_{g_2}^{g_0 g_1}),  \ldots , E_{g_0 \cdots g_n , 1} (a_{g_{n+1}}^{g_0 \cdots g_n})  \right) \displaybreak[2] \\
&= c \left( E_{1,g_0} (a_{g_0}) E_{g_0 , g_0 g_1} (a_{g_1}^{g_0}) + \lambda E_{g_0 , g_0 g_1} (a_{g_1}^{g_0}), E_{g_0 g_1 , g_0 g_1 g_2} (a_{g_2}^{g_0 g_1}),  \ldots , E_{g_0 \cdots g_n , 1} (a_{g_{n+1}}^{g_0 \cdots g_n})  \right) \displaybreak[2] \\
&= c \left( E_{1,g_0} (\widetilde{a_{g_0}}) E_{g_0 , g_0 g_1} (a_{g_1}^{g_0}), E_{g_0 g_1 , g_0 g_1 g_2} (a_{g_2}^{g_0 g_1}),  \ldots , E_{g_0 \cdots g_n , 1} (a_{g_{n+1}}^{g_0 \cdots g_n}) \right) \displaybreak[2] \\
&= (d^0 c) \left( E_{1,g_0} (\widetilde{a_{g_0}}) , E_{g_0 , g_0 g_1} (a_{g_0}^{g_1}),  \ldots , E_{g_0 \cdots g_n , 1} (a_{g_{n+1}}^{g_0 \cdots g_n}) \right) \displaybreak[2] \\
&= (\Psi_2 d^0) (c) (\widetilde{a_{g_0}} , \ldots , a_{g_{n+1}} ) 
\end{align*}
We also have $d^{n+1} \Psi_2 = \Psi_2 d^{n+1}$ as
\begin{align*}
&(d^{n+1} \Psi_2) ( c ) (a_{g_0} , \ldots , a_{g_{n+1}} ) \\
&= \Psi_2 (c) (a_{g_{n+1}} * a_{g_0}  , \ldots , a_{g_{n}}) \\
&= c (E_{1 , g_{n+1} g_0} (\mu (a_{g_{n+1}} \otimes a_{g_0}^{g_{n+1}}) ) , E_{g_{n+1} g_0 , g_{n+1} g_0 g_1} (a_{g_1}^{g_{n+1} g_0}) , \ldots , E_{g_{n+1} g_0 \cdots g_{n-1} , 1} (a_{g_n}^{g_{n+1} g_0 \cdots g_{n-1}})) \\
&=  c(E_{1 , g_{n+1} g_0} (\mu (a_{g_{n+1}} \otimes a_{g_0}^{g_{n+1}}) ) , E_{g_{n+1} g_0 , g_{n+1} g_0 g_1} (a_{g_1}^{g_{n+1} g_0}) , \ldots, E_{g_{n+1} g_0 \cdots g_{n-1} , 1} (a_{g_n}^{g_{n+1} g_0 \cdots g_{n-1}})) \cdot g_{n+1}^{-1}  \\
&= c(E_{g_0 \cdots g_n , g_0} (\mu (a_{g_{n+1}} \otimes a_{g_0}^{g_{n+1}}) ) , E_{g_0 , g_0 g_1} (a_{g_1}^{g_0}) , \ldots , E_{g_0 \cdots g_{n-1} , g_0 \cdots g_n} (a_{g_n}^{g_0 \cdots g_{n-1}}) ) \\
&=  c(E_{g_0 \cdots g_n , 1} (a_{g_{n+1}}) E_{1, g_0} ( a_{g_0}^{g_{n+1}} ) , E_{g_0 , g_0 g_1} (a_{g_1}^{g_0}) , \ldots , E_{g_0 \cdots g_{n-1} , g_0 \cdots g_n} (a_{g_n}^{g_0 \cdots g_{n-1}}) ) \\
&= (d^{n+1} c) ( E_{1, g_0} ( a_{g_0}^{g_{n+1}} )) , E_{g_0 , g_0 g_1} (a_{g_1}^{g_0}) , \ldots , E_{g_0 \cdots g_n , 1} (a_{g_{n+1}}^{g_0 \cdots g_n}) ) \\
&= (\Psi_2  d^{n+1}) (c) (a_{g_0} , \ldots , a_{g_{n+1}} )
\end{align*}
and the extension to $\widetilde{a_{g_0}} = (a_{g_0} , \lambda)$ is not much different from the previous calculation.  This shows that $b \Psi_2 = \Psi_2 b$.

Next we show that $B \Psi_2 = \Psi_2 B$.  As $B = \sum_i \bar{s}^i$ we will show that $\bar{s}^i \Psi_2 = \Psi_2 \bar{s}^i$ for $0 \leq i \leq n$.  The identity element in $C^\infty_c (M \rtimes \Gamma , L)$ is $(0,1)$ where $0$ may be considered as the sum of zero sections $0 : M \rightarrow L_g$ for any $g$, so we consider $0$ as $0 : M \rightarrow L_1$ below.  
\begin{align*}
& (\Psi_2 \bar{s}^i) (c) (\tilde{a}_{g_0} , \ldots , a_{g_n}) \\
&= (\bar{s}^i c) \left(E_{1, g_0} (\tilde{a}_{g_0}) , \ldots , E_{g_0 \cdots g_{i-1} , g_0 \cdots g_i} (a_{g_i}^{g_0 \cdots g_{i-1}}) , \ldots , E_{g_0 \cdots g_{n-1} , 1} (a_{g_n}^{g_0 \cdots g_{n-1}}) \right) \\
&= (\bar{s}^i c) \left(E_{1, g_0} (\tilde{a}_{g_0}) , \ldots , E_{g_0 \cdots g_{i-1} , g_0 \cdots g_i} (a_{g_i}^{g_0 \cdots g_{i-1}}) , \ldots , E_{g_0 \cdots g_{n-1} , 1} (a_{g_n}^{g_0 \cdots g_{n-1}}) \right) \cdot (g_i \cdots g_n) \\
&= (\bar{s}^i c) \left(E_{g_i \cdots g_n, g_i \cdots g_n g_0} (\tilde{a}_{g_0}^{g_i \cdots g_n}) , \ldots , E_{g_i \cdots g_n g_0 \cdots g_{i-2} , 1} (a_{g_{i-1}}^{g_i \cdots g_n g_0 \cdots g_{i-2}}), \right. \\
&\qquad \qquad \left. E_{1 , g_i } (a_{g_{i}}) , \ldots , E_{g_i \cdots g_{n-1} , g_i \cdots g_n} (a_{g_n}^{g_i \cdots g_{n-1}}) \right)  \\
&= c \left(1, E_{1,g_i} (a_{g_i}) , E_{g_i , g_i g_{i+1}} (a_{g_{i+1}}^{g_i}) , \ldots , E_{g_i \cdots g_{n-1} , g_i \cdots g_n} (a_{g_n}^{g_i \cdots g_n}), \right. \\ 
&\qquad \qquad \left. E_{g_i \cdots g_n , g_i \cdots g_n g_0} (a_{g_0}^{g_i \cdots g_n}) , \ldots , E_{g_i \cdots g_n g_0 \cdots g_{i-2} , 1} (a_{g_{i-1}}^{g_i \cdots g_n g_0 \cdots g_{i-2}})    \right) \\
&= (\Psi_2 c) (1, a_{g_i} , \ldots , a_{g_n} , a_{g_0} , \ldots , a_{g_i}) \\
&= (\bar{s}^i \Psi_2) (c) (\tilde{a}_{g_0} , \ldots , a_{g_n}). 
\end{align*}
Hence $\Psi_2$ is a morphism.
\end{proof}
\end{thm}

This completes the proof of the main theorem, Theorem \ref{thm:main_thm}.

\bibliographystyle{plain}

\bibliography{refs}

\begin{thebibliography}{10}

\bibitem{bott76}
Raoul Bott.
\newblock On characteristic classes in the framework of {G}elfand-{F}uchs
  cohomology.
\newblock {\em Ast{\'e}risque}, (32-33):113--139, 1976.

\bibitem{bott_tu}
Raoul Bott and Loring Tu.
\newblock {\em Differential {F}orms in {A}lgebraic {T}opology}.
\newblock Springer, New York, 1982.

\bibitem{bcmms02}
Peter Bouwknegt, Alan~L. Carey, Varghese Mathai, Michael~K. Murray, and Danny
  Stevenson.
\newblock Twisted {$K$}-theory and {$K$}-theory of bundle gerbes.
\newblock {\em Communications in Mathematical Physics}, 228(1):17--49, 2002.

\bibitem{bgnt07}
Paul Bressler, Alexander Gorokhovsky, Ryszard Nest, and Boris Tsygan.
\newblock Deformations of gerbes on smooth manifolds.
\newblock In {\em K-Theory and Noncommutative Geometry}, EMS Series of Congress
  Reports, pages 349--393. European Mathematical Society, 2008.

\bibitem{brylinski93}
Jean-Luc Brylinski.
\newblock {\em Loop Spaces, Characteristic Classes, and Geometric
  Quantization}.
\newblock Birkh{\"a}user, Boston-Basel-Berlin, 1993.

\bibitem{cmm00}
Alan~L. Carey, Jouko Mickellson, and Michael~K. Murray.
\newblock Bundle gerbes applied to quantum field theory.
\newblock {\em Reviews in Mathematical Physics}, 12(1):65--90, 2000.

\bibitem{connes94}
Alain Connes.
\newblock {\em Noncommutative Geometry}.
\newblock Academic Press, 1994.

\bibitem{crainic99}
Marius Crainic.
\newblock Cyclic {C}ohomology of \'{E}tale {G}roupoids; {T}he {G}eneral {C}ase.
\newblock {\em K-{T}heory}, 17, 1999.

\bibitem{dupont}
Johan~L. Dupont.
\newblock {\em Curvature and characteristic classes}.
\newblock Number 640 in Lecture Notes in Mathematics. Springer-Verlag,
  Berlin-Heidelberg-New York, 1978.

\bibitem{felisatti_neumann}
Marcello Felisatti and Frank Neumann.
\newblock Secondary theories for simplicial manifolds and classifying spaces.
\newblock {\em Geometry \& {T}opology {M}onographs}, 11:33--58, 2007.

\bibitem{gawedzki_reis02}
Krzysztof Gaw{\c{e}}dzki and Nuno Reis.
\newblock {WZW} branes and gerbes.
\newblock {\em Reviews in Mathematical Physics}, 14(12):1281--1334, 2002.

\bibitem{giraud71}
Jean Giraud.
\newblock {\em Cohomologie non ab{\'e}lliene}.
\newblock Number 179 in Die Grundlehren der mathematischen Wissenschaften.
  Springer-Verlag, Berlin-New York, 1971.

\bibitem{gorokhovsky99}
Alexander Gorokhovsky.
\newblock Characters of cycles, equivariant characteristic classes, and
  {F}redholm modules.
\newblock {\em Communications in Mathematical Physics}, 208(1):1--23, 1999.

\bibitem{hitchin01}
Nigel Hitchin.
\newblock Lectures on special {L}agrangian submanifolds.
\newblock In Cumrun Vafa and S.-T. Yau, editors, {\em Winter School on Mirror
  Symmetry, Vector Bundles and Lagrangian Submanifolds}, number~23 in Studies
  in Advanced Mathematics, pages 151--182, Providence, 2001. AMS/International
  Press.

\bibitem{jlo88}
Arthur Jaffe, Andrzej Lesniewski, and Konrad Osterwalder.
\newblock Quantum {$K$}-theory, {I}. {C}hern character.
\newblock {\em Communications in Mathematical Physics}, 118(1):1--14, 1988.

\bibitem{loday98}
Jean-Louis Loday.
\newblock {\em Cyclic Homology}.
\newblock Number 301 in Die Grundlehren der mathematischen Wissenschaften.
  Springer-Verlag, Berlin-Heidelberg-New York, 2nd edition, 1998.

\bibitem{mathai_stevenson06}
Varghese Mathai and Danny Stevenson.
\newblock On a generalised {C}onnes-{H}ochschild-{K}ostant-{R}osenberg theorem.
\newblock {\em Advances in Mathematics}, 200(2):303--335, 2006.

\bibitem{mathai_stevenson07}
Varghese Mathai and Danny Stevenson.
\newblock Entire cyclic homology of stable continuous trace algebras.
\newblock {\em Bulletin of the London Mathematical Society}, 39(1):71--75,
  2007.

\bibitem{moerdijk02}
Ieke Moerdijk.
\newblock Orbifolds as {G}roupoids: an {I}ntroduction.
\newblock In Alejandro Adem, Jack Morava, and Youngbim Ruan, editors, {\em
  Orbifolds in Mathematics and Physics}, number 310 in Contemporary
  Mathematics, pages 205--222, Providence, 2002. AMS.

\bibitem{moerdijk_mrcun}
Ieke Moerdijk and Janez Mr{\v{c}}un.
\newblock {\em Introduction to Foliations and Lie Groupoids}.
\newblock Cambridge Univesity Press, Cambridge, UK, 2003.

\bibitem{quillen}
Dan Quillen.
\newblock Algebra cochains and cyclic cohomology.
\newblock {\em Publ. Math. IHES}, 68:139--174, 1989.

\bibitem{tu_xu06}
Jean-Louis Tu and Ping Xu.
\newblock Chern character for twisted {$K$}-theory of orbifolds.
\newblock {\em Advances in Mathematics}, 207(2):455--483, 2006.

\end{thebibliography}

 \end{document}